\documentclass{elsart}

\usepackage{latexsym}
\usepackage{amsfonts}
\usepackage{graphicx,amssymb,enumerate,cite}
\usepackage{amsmath}

\usepackage{subfig}
 \usepackage[misc]{ifsym}

\usepackage{balance}

\newtheorem{tm}{Theorem}[section]   
\newtheorem{pp}[tm]{Proposition}   
\newtheorem{lm}[tm]{Lemma}

\newtheorem{tl}[tm]{Corollary}
\newtheorem{re}[tm]{Remark}
\newtheorem{eap}[tm]{Example}

\newcommand{\pof}{\noindent {\bf Proof} }
\newcommand{\ep}{$\quad \Box$}

\allowdisplaybreaks

\begin{document}

\newcommand{\al}{\alpha}
\newcommand{\be}{\beta}
\newcommand{\var}{\varepsilon}
\newcommand{\la}{\lambda}
\newcommand{\de}{\delta}
\newcommand{\st}{\stackrel}
\newcommand{\f}{F^{*}(R)}
\newcommand{\ce}{\centerline}
\newcommand{\Llr}{\Longleftrightarrow}
\newcommand{\Lr}{\Longrightarrow}

\begin{frontmatter}

\title{Approximation of fuzzy numbers by convolution method}

\author{Huan Huang  $^{a}$\corauthref{cor}  }
\author{Congxin Wu $^{b}$}
\author{}\ead{hhuangjy@126.com (H. Huang), wucongxin@hit.edu.cn (C. Wu)}

\address{
$^{a}$Department of Mathematics, Jimei
University, Xiamen 361021, China
\newline
 $^{b}$Department of
Mathematics, Harbin Institute of Technology, Harbin 150001, China
}
\corauth[cor]{Corresponding author.}
\date{}


%
%
%

%
%
%


\maketitle 

\begin{abstract}
In this paper we consider how to use the convolution method
to
construct approximations, which consist of fuzzy numbers sequences with good properties,
for a general fuzzy number.
It
shows
 that this convolution method can generate
differentiable approximations in finite steps for fuzzy numbers which have finite non-differentiable points.
 In the previous work,
this convolution method
only can be used to construct differentiable approximations
for continuous fuzzy numbers whose possible non-differentiable
points are the two endpoints of 1-cut.
The constructing of smoothers   is a key step in the construction process of approximations.
It
further points out
that,
if appropriately choose the smoothers,
then one can use the convolution method to provide approximations
which are differentiable, Lipschitz and preserve the core at the same time.

\end{abstract}


\end{frontmatter}


\section{Instructions}

  The approximations of fuzzy numbers attract many people's attention.
Mostly,
the researches can be grouped into two classes. One class is to use a given shape
fuzzy number to approximate the original fuzzy number.
There exist many important works which include but not limited to the following.
  Chanas \cite{cas}
and
Grzegorzewski \cite{grz1}
independently  presented the   interval approximations.
Ma et al. \cite{ma}
presented the
symmetric triangular approximations.
 Abbasbandy and Asady \cite{aay} presented
the
trapezoidal approximations.
 Grzegorzewski and Mr\'{o}wka \cite{grz2,  grz3}
presented the
trapezoidal approximations preserving the expected interval.
Zeng and Li \cite{zeng} presented the
weighted triangular approximations.
Nasibov and Peker \cite{nap} presented
the semi-trapezoidal approximations which is improved by Ban \cite{ban,ban2}.
 Yeh \cite{yeh}
presented the weighted semi-trapezoidal approximations.
 Yeh and Chu \cite{yeh2} presented a unified method to solve
the LR-type
approximation problems
without constraints
according to the weighted $L_2$-metric.
 Coroianu \cite{cor2}
discussed how to find the best Lipschitz constant of the trapezoidal approximation operator preserving the
value and ambiguity.
  Ban and Coroianu \cite{ban3}
proposed simpler methods to compute the parametric approximation
of a fuzzy number  preserving some important characteristics.
The works of this class of fuzzy numbers approximation
provide various methods to  approximate an arbitrary
fuzzy number
according to some metrics by  a  special type of fuzzy number
which is much more convenient to be calculated.
At the same time, since it finds the fuzzy number which has the
minimal distance to the original fuzzy number among all the given type fuzzy numbers,
it  minimizes
the loss of the information   to a certain extent.

But
there exist many   situations    in which the smaller the distance between approximated fuzzy number and
original fuzzy number becomes, the better the   effect appears.
So
it is also important to consider the problem whether one can approximate a fuzzy number
arbitrary well by fuzzy numbers with some good properties such as continuous, differentiable, etc.
This
is the topic of
another class of researches for fuzzy number approximation
which discuss how to construct a fuzzy numbers sequence with some properties
to
approximate a general fuzzy number.
There also exist many important contributions including the following works.
 Colling and Kloeden \cite{col} used the continuous fuzzy numbers sequence to
approximate an arbitrary fuzzy number.
 Coroianu et al. \cite{cor2016, cor}    constructed approximations which is made up of fuzzy numbers sequence
by using the F-transform and the max-product Bernstein operators, respectively.
Rom\'{a}n-Flores et al. \cite{ro}
pointed out a fact that
 the Lipschitzian
fuzzy numbers sequences can approximate any fuzzy numbers.
To demonstrate this fact,
they presented a method based on the convolution of two fuzzy numbers to
construct approximations for fuzzy numbers. For writing convenience, we call this method convolution method
in the sequel.
This
convolution method
                traces back to
the work of
 Seeger and   Volle \cite{sge}.

Differentiable fuzzy numbers play an important role in the
implementation of fuzzy intelligent systems and their applications
(see \cite{co,fs}).
For instance,
to use the well-known gradient descent algorithm,
it needs the fuzzy numbers be differentiable.
So it is an important and vital question that
whether one can use the differentiable
fuzzy numbers sequence to approximate
a general fuzzy number.
Chalco-Cano
et al. \cite{yo, yo2}
 used  the convolution method
to construct
  differentiable fuzzy numbers sequences to approximate a type of   non-differentiable fuzzy numbers
under supremum metric.
Since the convergence induced by supremum metric is stronger than the convergences induced by $L_p$-metric, sendograph metric, endograph metric, and  level-convergence (see \cite{da, wu, huang, huang3}),
it follows that
the constructed approximation is also an approximation for the original fuzzy number under the above mentioned convergences.
 This method is easy to be implemented since its operation is on level-cut-sets of the fuzzy numbers.
In
 the
 sequel, a differentiable fuzzy
number is also called a smooth fuzzy number,
  a non-differentiable fuzzy
number is called a non-smooth fuzzy number,
and
an approximation which consists of   differentiable fuzzy
numbers sequence is called a smooth approximation.

To construct smooth approximations for a type of non-smooth fuzzy numbers,
  Chalco-Cano et al.\cite{yo,yo2}   showed  an interesting fact that
the convolution transform can be used to smooth this type of fuzzy numbers, i.e. it can transfer
a non-smooth fuzzy number in this type to a smooth fuzzy number.
In fact, they constructed a class of
`smoothers' which are fuzzy numbers satisfying
some conditions.
Given a non-smooth fuzzy number of this type,   it can obtain a smooth fuzzy number
via
 convolution   of the original fuzzy number and the smoother.
The construction of
smoothers
is an important step in the construction of approximations.
The distance
 between the smooth fuzzy number and the
original fuzzy number
can be controlled
by the smoother.
Thus, by appropriately choosing the  smoother, it can get a smooth fuzzy number
such that the distance between which and the original fuzzy number is less than
an arbitrarily small  positive number given in advance.
So it can produce a sequence of smooth fuzzy numbers which constitute a smooth approximation of the original fuzzy numbers.

However, in the previous work, only a given type of fuzzy numbers
can be smoothed   by the convolution method.
Hence   only this type of fuzzy numbers can be smoothly approximated by using the convolution method.
 This type of fuzzy numbers have at most two possible non-differentiable points which are
the  endpoints of 1-cut.
 Whereas, an arbitrary fuzzy number may have other non-differentiable points,
or even non-continuous points.
So
it is  natural to consider
the question
 whether one can use the convolution method to smooth a general fuzzy number
and then the question
  whether one can give a smooth approximation to the original fuzzy number.

In this paper, we want to answer these questions. For this purpose, we first discuss the properties of
fuzzy numbers and convolution of fuzzy numbers.
Based on these discussions, we give partial positive answers to above questions.
The key is how to construct smoothers for a general fuzzy number so that it can be smoothed.
We do this step by step. It first shows how to construct smoothers for a subtype of continuous fuzzy numbers.
Then
it investigates how to construct smoothers
for the continuous fuzzy numbers.
At last,
it explores how to construct smoothers
for an arbitrary fuzzy number so that it can be transformed into a smooth fuzzy number.
On the basis of above
   results,
 it shows that how to construct smooth approximations
for
 fuzzy numbers which have finite non-differentiable points.
This type of fuzzy numbers are quite general in real world applications.
It
further
finds that,
by appropriately choosing the smoothers,
the smooth approximations can be Lipschitz approximations and can preserve the core
at the same time.
 We
 give simulation examples to validate and to illustrate the theoretical results.

The remainder of this paper is organized as follows.
Section 2
presents preliminaries about
fuzzy numbers and the convolution method for approximating fuzzy numbers.
Section 3 gives properties on the continuity of fuzzy numbers.
In Section 4, it shows  that the   convolution transform can keep the differentiability of fuzzy numbers
which
is
the key property to  ensure that the convolution method  can   be used to smooth fuzzy numbers.
On the basis of the results in Sections 3 and 4,
it
discusses how to smooth and approximate a general fuzzy number in Section 5.
In Section 6, it investigates advantages of constructing
approximations
 by the convolution method.
In Section 7, we draw conclusions.

%
%
%
%
%
%

%
%
%
%
%
%

\section{Preliminaries}

\subsection{Fuzzy numbers}

In this subsection, we introduce some basic and important notations and properties
about fuzzy numbers
               which will be used in the sequel.
  For details, we refer the reader to
references \cite{da,wu}.

Let $\mathbb{N}$ be the set of all natural numbers, $ \mathbb{R}$ be the
set of all real numbers.
A fuzzy subsets $u$ on $\mathbb{R}$ can be seen as a mapping  from
$\mathbb{R}$ to [0,1]. For $\al\in (0,1]$, let $[u]_{\al}$ denote
the $\al$-cut of $u$; i.e., \ $[u]_{\al}\equiv\{x\in \mathbb
R:u(x)\geq \al \}$ and $[u]_0$ denotes $\overline{\{x\in \mathbb
R:u(x)>0\}}$. We call $u$ a fuzzy number if $u$ has the following
properties:
\\
(\romannumeral1) $[u]_1\neq \emptyset$; and
\\
(\romannumeral2) $[u]_\al=[u^{-}(\al), u^{+}(\al)]$ are compact intervals of
$\mathbb{R}$ for all $\al\in [0,1]$.
\\
The set of all fuzzy numbers is denoted by
$\mathcal{F}(\mathbb{R})$. In \cite{yo}, a fuzzy number is also called a fuzzy interval.

Suppose that $u$ is
a fuzzy number. The 1-cut of $u$ is also called the core of $u$, which is denoted by Core($u$),
i.e.
Core($u)= [u]_1$.
 $u$
is said to be Lipschitz if $u$ is a Lipschitz function on $[u]_0$,
i.e.
$|u(x) - u(y)| \leq K |x-y|$ for all $x,y \in [u]_0$,
where $K$ is a constant which is called the Lipschitz contant.

The following  is a widely used representation theorem of
 fuzzy numbers.
\begin{pp}
\label{ulc}(Goetschel and Voxman \cite{gn})
 Given $u\in \mathcal{F}(\mathbb{R}),$ then

(\romannumeral1) \ $u^{-}(\cdot)$ is a left-continuous nondecreasing bounded
function on $(0,1]$;\\
(\romannumeral2) \ $u^{+}(\cdot)$ is a left-continuous nonincreasing
bounded function on $(0,1]$;\\
(\romannumeral3) \ $u^{-}(\cdot)$ and  $u^{+}(\cdot)$ are right continuous at
$\al=0$;\\
(\romannumeral4) \ $u^{-}(1)\leq u^{+}(1).$

Moreover, if the pair of functions $a(\la)$ and $b(\la)$ satisfy
conditions (\romannumeral1) through (\romannumeral4), then there exists a unique $u\in \mathcal{F}(\mathbb{R})$ such
that $[u]_{\alpha}=[a(\la),b(\la)]$ for each $\alpha\in (0,1].$
\end{pp}

From the definition of fuzzy numbers,
we know that,
given $x<u^-(1)$,
 then
  $u(x)=\lim_{z\to x+}u(z)$, i.e. $u$ is right-continuous at $x< u^-(1)$.
 Similarly,
 $u(x)=\lim_{z\to x-}u(z)$ for each $x>u^+(1)$, i.e., $u$ is left-continuous at each $x>u^+(1)$.

The algebraic operations
 on $\mathcal{F}(\mathbb{R})$ are defined  as follows: given $u,\
v\in \mathcal{F}(\mathbb{R})$, $\al\in [0,1]$,
  \begin{gather}
[u+v]_\al= [u]_\al+[v]_\al=[u^-(\alpha) + v^-(\alpha),\ u^+(\alpha) + v^+(\alpha)],\nonumber
\\
[u-v]_\al= [u]_\al-[v]_\al=[u^-(\alpha) - v^+(\alpha),\ u^+(\alpha) - v^-(\alpha)], \nonumber
\\
[u\cdot v]_\al=[u]_\al\cdot[v]_\al
=
[
\min  \{    xy:     \;     x\in [u]_\al,\, y\in[v]_\al   \}
,\
 \max  \{     xy:    \;     x\in [u]_\al,\, y\in [v]_\al   \}
 ].\label{fnm}
\end{gather}
From \eqref{fnm}, we know that if $r$ is a real number and $v$ is a fuzzy number, then
\[
(r\cdot v)(t) =
\left\{
\begin{array}{ll}
v(t/r), & \   r\not=0,
\\
\chi_{ \{0\}   }(t), & \   r=0,
\end{array}
\right.
\]
where $\chi_{ \{0\}   }$ is the characterization of $\{0\}$.

Suppose
that
$u$ is a fuzzy number. Its strong-$\alpha$-cuts $[u]_\al^s$, $\al\in [0,1]$,    are defined by:
\[[u]_\al^s=[u_s^-(\al), u_s^+(\al)] =
\left\{
\begin{array}{ll}
\overline{\bigcup_{\beta > \alpha} [u]_\beta } =   [\lim_{\beta \to \al+}  u^-(\beta), \lim_{\beta \to \al+}  u^+(\beta)], & \   \al<1,
\\
\mbox{}[u]_1=[u^-(1), u^+(1)], & \   \al=1.
\end{array}
\right.
\]
Clearly,
$[u]_1=[u]^s_1$, $[u]_0=[u]^s_0$, and $[u]^s_\al\subseteq    [u]_\al$
 for all $\al\in [0,1]$.
It is easy to show that
 \begin{gather*}
[u+v]^s_\al= [u]^s_\al+[v]^s_\al=[u^-_s(\alpha) + v^-_s(\alpha),\ u^+_s(\alpha) + v^+_s(\alpha)],  \\
[u-v]^s_\al= [u]^s_\al-[v]^s_\al=[u^-_s(\alpha) - v^+_s(\alpha),\ u^+_s(\alpha) - v^-_s(\alpha)].
\end{gather*}
We call $u^-(\cdot), \ u^+(\cdot), \ u_s^-(\cdot), \ u_s^+(\cdot)$ cut-functions.
The $\alpha$-cut and strong-$\alpha$-cut
are also called level-cut-set or strong-level-cut-set, respectively.

The supremum metric
on $\mathcal{F}(\mathbb{R})$
is defined by
$$
d_{\infty}(u,v)=\sup_{\al\in [0,1]} \max\{|u^{-}(\al)-v^{-}(\al)|,\
|u^{+}(\al)-v^{+}(\al)|\},
$$
where $u,v\in
\mathcal{F}(\mathbb{R}).$

\subsection{Convolution method for approximating fuzzy numbers}

This subsection describes a method based on the convolution transform
to   approximate a fuzzy number.
This
convolution method
                was first putted forward by
 Rom\'{a}n-Flores \cite{ro},
and
 traced back to
the work of
 Seeger and   Volle \cite{sge}.
Chalco-Cano
et al. \cite{yo, yo2}
gave important contributions to this method. They used  this convolution method
to produce
 smooth approximations for a class of   non-smooth fuzzy numbers.

The
 sup-min convolution $u\nabla v$ of fuzzy numbers $u$ and $v$ is defined by
$$(u\nabla v)(x)=\sup_{y\in \mathbb{R}}\{u(y)\wedge v(x-y)\}.$$

\begin{re}

In fact $u\nabla v=u+v$ for all $u,v\in \mathcal{F}(\mathbb{R})$. For details, see \cite{da, wu}.

\end{re}

%
%
%
%
%
%

The following is some symbols which are used to denote subsets of
$\mathcal{F}(\mathbb{R})$.
\begin{itemize}

\item \
$\mathcal{F}_\mathrm{T}(\mathbb{R}) $ is denoted
the family of all fuzzy numbers $u$ such that $u$ is strictly
increasing on
$[u^-(0), u^-(1)]$,
strictly decreasing on
$[u^+(1), u^+(0)]$,
 and
  differentiable on
   $(  u^-(0), u^-(1)   ) \cup  (   u^+(1), u^+(0)   )$.

\vspace{1mm}

\item \
$\mathcal{F}_\mathrm{N}(\mathbb{R})$ is denoted the family of all
fuzzy numbers $u$ such that $u$ is differentiable on $(u^-(0), u^-(1)) \cup  (u^+(1), u^+(0))$.

\vspace{1mm}

\item \
$\mathcal{F}_\mathrm{C}(\mathbb{R})$ is denoted the family of all
 fuzzy
numbers $u$ such that $u:\mathbb{R} \to [0, 1]$ is continuous
on $[u]_0   =  [u^-(0), u^+(0)]$. In other words,
given $u\in \mathcal{F} (\mathbb{R})$
 with $[u]_0$ is not a singleton,
then
$u\in \mathcal{F}_\mathrm{C}(\mathbb{R}) $ if and only if $u$ is continuous
on $(u^-(0), u^+(0))$, right-continuous on $u^-(0)$ and left-continuous on $u^+(0)$.

\vspace{1mm}

\item \
$\mathcal{F}_\mathrm{D}(\mathbb{R})$ is denoted the family of all
differentiable fuzzy numbers, i.e., the family of all fuzzy
numbers $u \in \mathcal{F}_\mathrm{C}(\mathbb{R})$ such that $u:\mathbb{R} \to [0, 1]$ is differentiable
on $(u^-(0), u^+(0))$.

\end{itemize}
Given
 a
 fuzzy number $u$ in $\mathcal{F}_\mathrm{N}(\mathbb{R})$,
$u$ need not be strictly increasing on $(u^-(0), u^-(1))$  and
strictly decreasing on $(u^+(1), u^+(0))$. So
$$\mathcal{F}_\mathrm{T}(\mathbb{R})
 \subsetneq
\mathcal{F}_\mathrm{N}(\mathbb{R}).$$
Observe
 that, for each $v\in \mathcal{F}(\mathbb{R})$,
  $v$ is differentiable on $(v^-(1),   v^+(1))$.
Thus,
 for every $u\in \mathcal{F}_\mathrm{N}(\mathbb{R})$,
 its   possible non-differentiable points in $(u^-(0),   u^+(0))$ are $u^-(1)$ and $u^+(1)$.
It
is easy to check that
\begin{gather*}
 \mathcal{F}_\mathrm{D}(\mathbb{R})
 \subsetneq
  \mathcal{F}_\mathrm{C}(\mathbb{R})
 \cap
      \mathcal{F}_\mathrm{N}(\mathbb{R}).
\end{gather*}
Clearly, given  $u\in \mathcal{F}_\mathrm{D}(\mathbb{R})$, if $x\in [u]_1$  and $x$ is an inner point of $[u]_0$, then
 $u'(x)=0$.
We also call $u\in \mathcal{F}_\mathrm{D}(\mathbb{R})$ a smooth fuzzy number.

Suppose that $u\in \mathcal{F}(\mathbb{R})$ and $ v\in \mathcal{F}_\mathrm{D}(\mathbb{R})$,
then $v$ is said to be a smoother of $u$ if $u\nabla v \in   \mathcal{F}_\mathrm{D}(\mathbb{R})$.

 Chalco-Cano et al. \cite{yo} pointed out that each fuzzy number in $\mathcal{F}_\mathrm{T}(\mathbb{R})$
can be approximated by a smooth fuzzy numbers sequence which is constructed by using the convolution method.
They
constructed
fuzzy numbers $w_p$, $p>0$, as follows:
\begin{equation} \label{wpc}
 w_p(x)=\left\{
\begin{array}{ll}
1-{\left(\frac{x}{p}\right)}^2,& \mbox{if}\ x\in [-p,p],\\
0, & \mbox{if}\ x\notin [-p,p].
\end{array}
\right.
\end{equation}
Obviously, $w_p\in \mathcal{F}_\mathrm{D}(\mathbb{R})$ for all  $p>0$.
They
 presented the
 following
 result.
\begin{pp}
\cite{yo} \label{1}
If
$u\in \mathcal{F}_\mathrm{T}(\mathbb{R})$,
then
$u\nabla w_p \in \mathcal{F}_\mathrm{D}(\mathbb{R})$.
\end{pp}
Notice
that
 $d_\infty (  u, u\nabla w_{p}   )   \to    0$
as $p\to 0$.
Thus Proposition \ref{1}
indicates
that
every fuzzy number in $\mathcal{F}_\mathrm{T}(\mathbb{R})$ can be approximated by
fuzzy numbers sequences contained in $\mathcal{F}_{\mathrm{D}}(\mathbb{R})$.

We can see that
the fuzzy numbers $w_p$,  $p>0$, work as  smoothers, which transfer each fuzzy number $u$
to a smooth fuzzy number $u \nabla w_p$. The  smooth fuzzy numbers sequence $\{ u \nabla w_p \}$ construct a smooth approximation
 of the original fuzzy number $u$,
i.e.,
 $u\nabla w_{p}    \to  u$ as $p\to 0$.

 Chalco-Cano et al. \cite{yo2}
  further presented an approach to   produce a more large class of  smoothers.
       A class of fuzzy numbers $Z_p^f$ are defined by
\begin{equation}\label{zpf}
  Z_p^f(x)=
\begin{cases}
f^{-1} ( \|x\| / p), &   \|x\|  \leq   p,
\\
0,      &     \|x\|  >   p,
\end{cases}
\end{equation}
where $p>0$ is a real number
 and
     that   $f: [0,1]   \to  [0,1]$ is a continuous and strictly decreasing function with $f(0)=1, f(1)=0$.
It
 is
 easy to see that $Z_p^f=w_p$ when $f=\sqrt{1-t}$.
They gave the following
result.
\begin{pp} \cite{yo2}\label{2}
Suppose that $Z_p^f$ is defined by \eqref{zpf}.
If $f$ is differentiable and $\lim_{\alpha \to 1-} f'(\al)   =  -\infty$,
then
$u \nabla Z_p^f    \in   \mathcal{F}_\mathrm{D}(\mathbb{R})$
for each $u\in   \mathcal{F}_\mathrm{T} (\mathbb{R})$.
\end{pp}
Notice
    that $d_\infty ( u \nabla Z_p^f, \   u  ) \to  0 $ as $p\to 0$.
This means that
given $f$ satisfies the above conditions,
it produces a class of smoothers $\{Z_p^f:  \  p>0 \}$
and
a
smooth approximation
$\{  u \nabla Z_p^f:  \  p>0 \}$ of the fuzzy number $u$.
Different $f$ corresponds to different class of smoothers and then corresponds to
different    smooth approximation.

\section{Properties of fuzzy numbers}

In this section, we investigate some properties on the continuity of fuzzy numbers which will be used in the sequel.
It
 first
 lists some conclusions on the values of membership functions of fuzzy numbers,
and
conclusions
on
 characterizations of
         continuous points of  fuzzy numbers.
Based on this,
it
gives some characterizations of
         continuous intervals of  fuzzy numbers.
At last,
it considers the properties of
              continuous points of the cut-functions of  fuzzy numbers.

We list some
 propositions and corollary which
can
be found in \cite{da,wu} or as direct consequences of the conclusions therein.
The following four
conclusions
discuss values of fuzzy numbers' membership functions.

\begin{pp}  \label{bse}

Suppose that $u\in F(\mathbb{R})$ and that $x\in    \mathbb{R}$, then the following statements hold.
\\
(\romannumeral1) \ If $u^-(1) \geq x> u_s^-(\al)$, then $u(x)>\al$.
\\
(\romannumeral2) \ If $u^+(1) \leq x < u_s^+(\al)$, then $u(x)>\al$.
\\
(\romannumeral3) \ If $x < u_s^-(\al)$ or $x > u_s^+(\al)$, then $u(x)\leq   \al$.

\end{pp}

\begin{pp}\label{nve} \ Suppose that $u\in F(\mathbb{R})$ and that $x\in    \mathbb{R}$, then the following statements hold.
\\
(\romannumeral1) \ If $u^-(\al)< u^-_s(\al)$, then $u(x)= \alpha$ when $x\in (u^-(\alpha), u^-_s(\al))$.
\\
(\romannumeral2) \ If $u^+(\al)> u^+_s(\al)$, then $u(x)=\alpha$ when $x \in ( u^+_s(\al), u^+(\alpha))$.
\end{pp}

\begin{pp}\label{vec}
Suppose that $u\in \mathcal{F}(\mathbb{R})$.
If $u$ is continuous at
a point $x\in (u^-(0), u^+(0))$ which
is equal to $u^-(\alpha)$ or $u^+(\alpha)$ or  $u_s^-(\alpha)$ or $u_s^+(\alpha)$,
then $u(x)=\al$.
\end{pp}

\begin{tl}
\label{ncp}
If $u\in \mathcal{F}_\mathrm{C}(\mathbb{R})$, then
\begin{gather*}
u(u^-(\al))=u(u_s^-(\al))=\al,\\
u(u^+(\beta))=u(u_s^+(\beta))=\beta
\end{gather*}
for all $\al  \geq  u(u^-(0))$ and $\beta \geq u(u^+(0))$.
\end{tl}

Propositions \ref{cae} and \ref{pce} consider characterizations of continuous points of a fuzzy number.

\begin{pp} \label{cae}
Suppose
that $u\in \mathcal{F}(\mathbb{R})$,
then the following statements hold.
\\
(\romannumeral1) \  Given $x \in  (  u^-(0),  u^+(1)  ]$,
 then
$u $ is
 left-continuous at $x$ if and only if $u^-(\beta)   <   x$ for each $\beta < u(x)$.
\\
 (\romannumeral2) \  Given $x \in  [  u^-(1),  u^+(0)  )$,
 then
$u $ is
 right-continuous at $x$ if and only if $u^+(\beta)   >   x$ for each $\beta < u(x)$.

\end{pp}

\begin{pp} \label{pce}
Suppose that $u\in \mathcal{F}(\mathbb{R})$.
 Then the following statements hold.
 \\
(\romannumeral1) \ Given $ u^-(\al) \in (  u^-(0),  u^-(1)   )$, then $u$ is continuous at $u^-(\al)$ if and only if
$u( u^-(\al)  ) = \alpha$, and
 $u^-(\beta)   <    u^-(\al)$ for each $\beta<\alpha$.
\\
(\romannumeral2) \ Given $ u_s^-(\al) \in (  u^-(0),  u^-(1)   )   $, then $u$ is continuous at $u_s^-(\al)$ if and only if
$u( u_s^-(\al)  ) = \alpha$, and
 $u^-(\beta)   <    u_s^-(\al)$ for each $\beta<\alpha$.
\\
(\romannumeral3) \ Given $ u^+(\al) \in (  u^+(1),  u^+(0)   )$, then $u$ is continuous at $u^+(\al)$ if and only if
$u( u^+(\al)  ) = \alpha$, and
 $u^+(\beta)   >    u^+(\al)$ for each $\beta<\alpha$.
\\
(\romannumeral4) \ Given $ u^+_s(\al) \in (  u^+(1),  u^+(0)   )$, then $u$ is continuous at $u_s^+(\al)$ if and only if
$u( u_s^+(\al)  ) = \alpha$, and
 $u^+(\beta)   >    u_s^+(\al)$ for each $\beta<\alpha$.

\end{pp}

The following lemmas and theorems give characterizations of
 continuous intervals of a fuzzy number.

\begin{lm} \label{lec}
Suppose that $u\in \mathcal{F}(\mathbb{R})$.
Given $a,b \in   [u^-(0),   u^-(1)] $ with $a<b$,
 then the following statements are equivalent.
\\
(\romannumeral1) There exists $x$ in $(a, b]$
such that $u$ is not left-continuous at $x$.
\\
(\romannumeral2) There exists $\alpha,   \beta$ in $[u(a), u(b)]$
such that $\alpha \not= \beta$
and
$u^-(\al) = u^-(\beta)$.
\\
(\romannumeral3)  $u^-(\cdot)$ is not strictly increasing on $[u(a),  u(b)]$.

\end{lm}

\pof \ If statement (\romannumeral1) holds, then there exists $x \in (a, b]$
such that
$u$ is not left-continuous at $x$.
Hence
$\alpha=u(x)> \lim_{z\to x-}u(z)=\beta$,
and thus
 $u^-(\al)=u^-(\frac{\alpha+\beta}{2})=x$.
 Note that
  $\alpha \not= \frac{\alpha+\beta}{2}$ and
 $\alpha,\beta\in [u(a), u(b)]$.
 So
statement (\romannumeral2)
holds.

 If statement (\romannumeral2) holds,
 then there exist $\al>\beta$ such that $u^-(\al) = u^-(\beta)=x$ and $\alpha,\beta\in [u(a), u(b)]$. So  $\lim_{z\to x-}u(z)
\leq \beta < \alpha\leq u(x)$. This means that $u$ is not left-continuous at $x\in (a, b]$,
i.e.,
statement (\romannumeral1)
holds.

The equivalence of statement (\romannumeral2) and statement (\romannumeral3)
follows immediately from
the monotonicity of $u^-(\cdot)$.
\ep

\begin{lm} \label{rec}
Suppose that $u\in \mathcal{F}(\mathbb{R})$.
Given $c,d \in   [u^+(1),   u^+(0)] $ with $c<d$,
 then the following statements are equivalent.
\\
(\romannumeral1) There exists $x$ in $[ c,d)$
such that $u$ is not right-continuous at $x$.
\\
(\romannumeral2) There exists $\alpha,   \beta$ in $[u(c), u(d)]$
such that $\alpha \not= \beta$
and
$u^+(\al) = u^+(\beta)$.
\\
(\romannumeral3)  $u^+(\cdot)$ is not strictly decreasing on $[u(c),  u(d)]$.
\end{lm}

\pof \ The proof is similar to the proof of Lemma \ref{lec}. \ep

\begin{tm} \label{cfn}
Suppose that $u\in \mathcal{F}(\mathbb{R})$. Then
$u\in \mathcal{F}_\mathrm{C}(\mathbb{R})$ if and only if the following
two conditions are satisfied.
\\
(\romannumeral1) \
$u^-(\al) \not= u^-(\beta)$ for all $\alpha,\beta\in [u(u^-(0)), 1]$ with $\alpha \not= \beta$.
\\
(\romannumeral2) \
 $u^+(\al) \not= u^+(\beta)$
for all $\alpha,\beta\in [ u(u^+(0)),1]$ with $\al\not = \beta$.
\end{tm}

\pof \ If $[u]_0$ is a singleton, the the conclusion holds obviously.
If $[u]_0$ is not a singleton, the desired results follow from Lemmas \ref{lec}
and \ref{rec}. \ep

%
%
%
%
%
%

\begin{tm}\label{smn}
Suppose that $u\in \mathcal{F}(\mathbb{R})$
  and that $u(u^-(0))=\alpha_0$ and $u(u^+(0))=\beta_0$.
  Then
   $u\in \mathcal{F}_\mathrm{C}(\mathbb{R})$
   if and only if
  $u^-(\cdot)$ is strictly increasing on $[\alpha_0,1]$ and $u^+(\cdot)$ is strictly decreasing on $[\beta_0,1]$.
\end{tm}

\pof \      Suppose that $u\in \mathcal{F}_\mathrm{C}(\mathbb{R})$, it then follows
from Lemmas \ref{lec} and \ref{rec}
that
  $u^-(\cdot)$ is strictly increasing on $[\alpha_0,1]$ and $u^+(\cdot)$ is strictly decreasing on $[\beta_0,1]$.

Suppose
that
$u^-(\cdot)$ is strictly increasing on $[\alpha_0,1]$ and $u^+(\cdot)$ is strictly decreasing on $[\beta_0,1]$.
Given
$u\in \mathcal{F}(\mathbb{R})$ with $[u]_0$ is not a singleton, then
$u^-(0) < u^-(1)$
or
$u^-(0) = u^-(1) <  u^+(1) $
or
$u^-(0) = u^-(1) =  u^+(1)< u^+(0) $.
We can check that in all above cases, $u$ is right-continuous on $u^-(0)$.
Similarly, we know
that
$u$ is left-continuous on $u^+(0)$.
By using Lemmas \ref{lec} and \ref{rec},
we
can
also
check that $u$ is continuous on $(u^-(0),  u^+(0))$.
Combined with above conclusions,
we know
that
  $u\in \mathcal{F}_\mathrm{C}(\mathbb{R})$. \ep

\begin{lm} \label{evf}
Suppose that $u,v \in \mathcal{F}(\mathbb{R})$ and that $\al\in [0,1]$, then
\begin{gather*}
(u\nabla v) ((u\nabla v)^-(\alpha))=u(u^-(\alpha))\wedge v(v^-(\alpha)),
\\
(u\nabla v) ((u\nabla v)^+(\alpha))=u(u^+(\alpha))\wedge v(v^+(\alpha)).
\end{gather*}
\end{lm}

\pof \ From Zadeh's extension principle, we know that
\begin{align*}
&(u\nabla v) ((u\nabla v)^-(\alpha))
\\
&=\sup\{u(x) \wedge v(y): x+y=(u\nabla v)^-(\alpha)=u^-(\al)+v^-(\al) \}\\
&\geq u(u^-(\al))   \wedge  v(v^-(\al))
\end{align*}
 for all $\al\in [0,1]$.

  On the other hand, given $x+y=u^-(\al)+v^-(\al)$, if $x>u^-(\al)$, then $y<v^-(\al)$, hence $v(y) \leq \al$, and thus
$u(x) \wedge v(y) \leq \al \leq u(u^-(\al))   \wedge  v(v^-(\al))$. Similarly we can show that if $x<u^-(\al)$, then
 $u(x) \wedge v(y) \leq \al \leq u(u^-(\al))   \wedge  v(v^-(\al))$. So we know that
 $$(u\nabla v) ((u\nabla v)^-(\alpha))=u(u^-(\al))   \wedge  v(v^-(\al)).$$
In the same way, we can prove that
$$(u \nabla v) ((u\nabla v)^+(\alpha))=u(u^+(\al))   \wedge  v(v^+(\al)). \quad \Box$$

\begin{tm}
\label{cst}
Let $u\in \mathcal{F}_\mathrm{C}(\mathbb{R})$ and let $v\in  \mathcal{F}(\mathbb{R})$.
Suppose that
 $u(u^-(0))=\alpha_0$ and $u(u^+(0))=\beta_0$,
 then
$u\nabla v$ is continuous on $[   (u\nabla v)^-(\alpha_0),     (u\nabla v)^+(\beta_0)    ]$.

\end{tm}

\pof \ Since $u(u^-(0))=\alpha_0$ and $u(u^+(0))=\beta_0$, we know that
$u^-(0) = u^-(\al_0)$ and $u^+(0) = u^+(\beta_0)$.
Thus, by Lemma \ref{evf},
\begin{gather*}
(u\nabla v)  (   (u\nabla v)^-(\alpha_0)   )= u(u^-(\al_0)) \wedge  v(v^-(\al_0)) = \alpha_0,
\\
(u\nabla v)  (   (u\nabla v)^+(\beta_0)   )= u(u^+(\beta_0)) \wedge  v(v^+(\beta_0)) = \beta_0.
\end{gather*}
Since
 $u\in \mathcal{F}_\mathrm{C}(\mathbb{R})$,
  by Theorem \ref{smn},
we know
   $u^-(\cdot)$ is strictly increasing on $[\alpha_0,1]$
   and
    $u^+(\cdot)$ is strictly decreasing on $[\beta_0,1]$.
 Note that
 \begin{gather*}
(u\nabla v)^-(\cdot)= u^-(\cdot) +  v^-(\cdot),
\\
(u\nabla v)^+(\cdot)= u^+(\cdot) +  v^+(\cdot),
\end{gather*}
  and hence
    $(u\nabla v)^-(\cdot)$ is strictly increasing on $[\alpha_0,1]$
    and
    $(u\nabla v)^+(\cdot)$ is strictly decreasing on $[\beta_0,1]$.
    Then, using Lemmas \ref{lec} and \ref{rec}
and
 reasoning as in
   the proof Theorem \ref{smn},
    we can prove that
if $[   (u\nabla v)^-(\alpha_0),     (u\nabla v)^+(\beta_0)]$ is not a singleton,
then
$u\nabla w$ is right-continuous on $ (u\nabla v)^-(\alpha_0)$,
 left-continuous on $(u\nabla v)^+(\beta_0)$
and
continuous on $( (u\nabla v)^-(\alpha_0),  (u\nabla v)^+(\beta_0))$.
So $u\nabla w$ is continuous
on
$[   (u\nabla v)^-(\alpha_0),     (u\nabla v)^+(\beta_0)]$.
\ep

\begin{tm}\label{arp}
Suppose that $u\in \mathcal{F}_\mathrm{C}(\mathbb{R})$ and that
 $v\in \mathcal{F}(\mathbb{R})$.
 If $u$ satisfies that
 $u(u^-(0)) \leq v(v^-(0))$ and $u(u^+(0)) \leq v(v^+(0))$,
 then
$u \nabla v\in \mathcal{F}_\mathrm{C}(\mathbb{R})$.
\end{tm}

\pof \ From Lemma \ref{evf}, we know that
\begin{gather*}
(u\nabla v)((u\nabla v)^-(0))=u(u^-(0)),
\\
(u\nabla v)((u\nabla v)^+(0))=u(u^+(0)).
\end{gather*}
So, by Theorem \ref{cst}, we know that
$u\nabla v$
is continuous on
$[   (u\nabla v)^-(0),  (u\nabla v)^+(0)    ]$,
i.e.
$u\nabla v \in \mathcal{F}_\mathrm{C}(\mathbb{R})$. \ep

The following theorem considers the properties of continuous points of cut-functions.

\begin{tm} \label{fcn}
Suppose that $u\in \mathcal{F}(\mathbb{R})$, then the following statements hold.
\\
(\romannumeral1)\ $u^-(\cdot)$ is continuous at $\al$, if and only if, $u^-(\al)=u_s^-(\al)$.
\\
(\romannumeral2)\ $u^+(\cdot)$ is continuous at $\al$, if and only if, $u^+(\al)=u_s^+(\al)$.
\end{tm}

\pof \
From Proposition \ref{ulc}, we know that $u^-(\cdot)$ is discontinuous at $\al$, if and only if,
 $u^-(\cdot)$ is not right continuous at $\al$, i.e. $u^-(\al)< \lim_{\beta \rightarrow \alpha+}   u^-(\beta)=u_s^-(\al)$. So statement (\romannumeral1) is true. Statement (\romannumeral2)
can be proved similarly.    \ep

\begin{re}
The assumption that
$u\in \mathcal{F}_\mathrm{C}(\mathbb{R})$, even the stronger assumption that $u\in \mathcal{F}_\mathrm{D}(\mathbb{R})$, cannot imply that $u^-(\cdot)$ and $u^+(\cdot)$ are continuous on [0,1]. Also, the assumption that  $u^-(\cdot)$ and $u^+(\cdot)$ are continuous on [0,1] cannot imply that $u\in \mathcal{F}_\mathrm{C}(\mathbb{R})$.
\end{re}

%
%
%
%
%
%

\section{Properties of convolution of fuzzy numbers} \label{pcm}

In this section, we investigate some properties of convolution of two fuzzy numbers.
It finds  that the   convolution transform can keep the differentiability of fuzzy numbers.
This
is
the key property which ensures that the convolution method  can   be used to smooth and to approximate fuzzy numbers.

The following theorem states that convolution  transform can retain the differentiability when
the derivative is zero.
The
 symbols
  $f'_-(\cdot)$ and $f'_+(\cdot)$ are used to denote the left derivative and
the right derivative of $f$, respectively.

\begin{tm}    \label{dzeron}
Let $u,v\in \mathcal{F}(\mathbb{R}) $, and let
 $\al\in [0,1]$,
then
 the following statements hold.
\\
(\romannumeral1)
If $v_-'(v^-(\al))=0$, then
$(u\nabla v)_-'((u\nabla v)^-(\al))=0$.
\\
(\romannumeral2)
If $v_+'(v^-(\al))=0$, and $v(v^-(\alpha))=  \beta$, then
$(u\nabla v)_+'((u\nabla v)^-(\beta))=0$.
\\
(\romannumeral3)
If $v_-'(v_s^-(\al))=0$, then
$(u\nabla v)_-'((u\nabla v)_s^-(\al))=0$.
\\
(\romannumeral4)
If $v_+'(v_s^-(\al))=0$, and $v(v_s^-(\al))=\beta$, then
$(u\nabla v)_+'((u\nabla v)_s^-(\beta))=0$.
\\
(\romannumeral5)
If $v_-'(v^+(\alpha))=0$,  and $v(v^+(\al))=\beta$, then
$(u\nabla v)_-'((u\nabla v)^+(\beta))=0$.
\\
(\romannumeral6)
If $v_+'(v^+(\alpha))=0$, then
$(u\nabla v)_+'((u\nabla v)^+(\alpha))=0$.
\\
(\romannumeral7)
If $v_-'(v_s^+(\alpha))=0$,  and $v(v_s^+(\alpha))=\beta$, then
$(u\nabla v)_-'((u\nabla v)_s^+(\beta))=0$.
\\
(\romannumeral8)
If $v_+'(v_s^+(\alpha))=0$,
then
$(u\nabla v)_+'((u\nabla v)_s^+(\alpha))=0$.

\end{tm}

\pof \ See Appendix A.   \ep

The following theorem expresses the fact that the differentiability at the left-endpoints of level-cut-sets
still   holds
 after convolution transform.  Furthermore, it gives the corresponding derivatives.

\begin{tm}
\label{dbzeron}

Let $u,v\in \mathcal{F}(\mathbb{R}) $, and let
 $\al\in [0,1]$,
then
\\
(\romannumeral1) \
If $u_-'(u^-(\al))=\varphi>0$ and $v_-'(v^-(\al))=\psi>0$,
then
$(u\nabla v)_-'(     (u\nabla v)  ^- (\al)    )=(\varphi^{-1} + \psi^{-1})^{-1}$.
\\
(\romannumeral2) \
If $u_+'(u^-(\al))=\varphi>0$,     $v_+'(v^-(\al))=\psi>0$, and
$u(u^-(\al))=v(v^-(\al))=\beta$, then
$(u\nabla v)_+'(     (u\nabla v)  ^- (\beta)    )=(\varphi^{-1} + \psi^{-1})^{-1}$.
\\
(\romannumeral3) \
If $u_+'(u^-(\al))=\varphi$,
$u(u^-(\al))=\beta$, and $v(v^-(\al))=\gamma>\beta$, then
$(u\nabla v)_+'(     (u\nabla v)  ^- (\beta)    )=\varphi$.
\\
(\romannumeral4) \
If $u_-'(u^-(\al))=\varphi$,
 and $\lim_{y\to v^-(\al)-}  v(y) = \lambda  <   \alpha$,
then
$(u\nabla v)_-'(     (u\nabla v)  ^- (\alpha)    )=\varphi$.
\end{tm}

\pof \ See   Appendix B. \ep

The following theorem shows that convolution transform keeps the differentiability at the right-endpoints of level-cut-sets. It also computes     the corresponding derivatives.

\begin{tm}
\label{rdbzeron}

Let $u,v\in \mathcal{F}(\mathbb{R}) $, and let
 $\al\in [0,1]$,
then
\\
(\romannumeral1) \
If $u_+'(u^+(\al))=\varphi>0$ and $v_+'(v^+(\al))=\psi>0$,
then
$(u\nabla v)_+'(     (u\nabla v)  ^+ (\al)    )=(\varphi^{-1} + \psi^{-1})^{-1}$.
\\
(\romannumeral2) \
If $u_-'(u^+(\al))=\varphi>0$,     $v_-'(v^+(\al))=\psi>0$, and
$u(u^+(\al))=v(v^+(\al))=\beta$, then
$(u\nabla v)_-'(     (u\nabla v)  ^+ (\beta)    )=(\varphi^{-1} + \psi^{-1})^{-1}$.
\\
(\romannumeral3) \
If $u_-'(u^+(\al))=\varphi$,
$u(u^+(\al))=\beta$, and $v(v^+(\al))=\gamma>\beta$, then
$(u\nabla v)_-'(     (u\nabla v)  ^+ (\beta)    )=\varphi$.
\\
(\romannumeral4) \
If $u_+'(u^+(\al))=\varphi$,
 and $\lim_{y\to v^+(\al)+}  v(y) = \lambda  <   \alpha$,
then
$(u\nabla v)_+'(     (u\nabla v)  ^+ (\alpha)    )=\varphi$.

\end{tm}

\pof \ The proof is similar to the proof of Theorem \ref{dbzeron}. \ep

\section{Smooth approximations of fuzzy numbers generated by the convolution method} \label{abn}

In this section, we consider how to use the convolution method
to
give smooth approximations for
an arbitrarily given fuzzy number.
The key step is the constructing of smoothers.
It discusses how to
   construct smoothers
for a general fuzzy number
so that it can be smoothed by
   the convolution method.
We show how to do this  step by step.
Firstly, it shows how to construct smoothers for
fuzzy numbers in $ \mathcal{F}_\mathrm{N}(\mathbb{R})\cap \mathcal{F}_\mathrm{C}(\mathbb{R})$.
Secondly,
it discusses how to construct smoothers for
fuzzy numbers in $ \mathcal{F}_\mathrm{C}(\mathbb{R})$.
At last,
 it investigates how to construct smoothers
for
 an arbitrary fuzzy number so that it can be smoothed.
Based on above results,
we then assert that,
given an arbitrary fuzzy number with finite non-differentiable points,
 one
can
use the convolution method
to
generate smooth approximations
 for it in finite steps.
The condition that the number of non-differentiable points is finite
is
quite general for fuzzy numbers used in real world applications.
Several simulation examples are given to validate and illustrate the theoretical results.
All
computations in this section are implemented by
Matlab.

In the following,
 we give some lemmas to discuss differentiability of the convolution of two fuzzy numbers at various types of its inner points.


\begin{lm}
\label{rapfes}
Suppose that $u\in \mathcal{F} (\mathbb{R})$   and that
 $w\in \mathcal{F} (\mathbb{R})$.
Then the following statements hold.
\begin{description}
  \item[A1] $(u\nabla w)'(x) =  (\varphi^{-1}   +    \psi^{-1})^{-1}$
when
$ x=(u\nabla w)^-(\al)$, $u'(u^-(\al)):=\varphi>0$, and $w'(w^-(\al)):=\psi>0$.

\item[A2] $(u\nabla w)'(x) = 0$
when
$ x=(u\nabla w)^-(\al)$ and $u'(u^-(\al)) \cdot w'(w^-(\al))=0$.

\item[A3] $(u\nabla w)'(x) =  (\varphi^{-1}   +    \psi^{-1})^{-1}$
when
$ x=(u\nabla w)_s^-(\al)$, $u'(u_s^-(\al))=:\varphi>0$, and $w'(w_s^-(\al))=:\psi>0$.

\item[A4] $(u\nabla w)'(x) = 0$
when
$ x=(u\nabla w)_s^-(\al)$ and $u'(u_s^-(\al)) \cdot w'(w_s^-(\al))=0$.

 \item[A5] $(u\nabla w)'(x) =  (\varphi^{-1}   +    \psi^{-1})^{-1}$
when
$ x=(u\nabla w)^+(\al)$, $u'(u^+(\al)):=\varphi< 0$, and $w'(w^+(\al)):=\psi<0$.

\item[A6] $(u\nabla w)'(x) = 0$
when
$ x=(u\nabla w)^+(\al)$ and $u'(u^+(\al)) \cdot w'(w^+(\al))=0$.

\item[A7] $(u\nabla w)'(x) =  (\varphi^{-1}   +    \psi^{-1})^{-1}$
when
$ x=(u\nabla w)_s^+(\al)$, $u'(u_s^+(\al))=:\varphi < 0$, and $w'(w_s^+(\al))=:\psi < 0$.

\item[A8] $(u\nabla w)'(x) = 0$
when
$ x=(u\nabla w)_s^+(\al)$ and $u'(u_s^+(\al)) \cdot w'(w_s^+(\al))=0$.

\item[A9] $(u\nabla w)'(x) = 0$
when
$ x \in (   (u\nabla w)^-(1),    (u\nabla w)^+(1)  )$   or  $x\in \bigcup_{0<\al<1}   ((u\nabla w)^-(\alpha), (u\nabla w)_s^-(\alpha) )$
 or
$x\in \bigcup_{0<\al<1}(  (u\nabla w)_s^+(\alpha), (u\nabla w)^+(\alpha)  )  $.

\end{description}

\end{lm}

\pof \ We only prove statements \textbf{A1}, \textbf{A2} and \textbf{A9}. Other statements can be proved similarly.

Suppose that $x$, $u^-(\al)$ and $w^-(\al)$
satisfy
 the   premise of statement \textbf{A1}.
By Proposition \ref{vec},
we know
that
$u(u^-(\al))=w(w^-(\al))=\al$.
Since
$u'(u^-(\al))=\varphi>0$
and
$w'(w^-(\al))=\psi>0$,
by Theorem
\ref{dbzeron} (\romannumeral1), (\romannumeral2),
we have that
$(u\nabla w)'(x) = (u\nabla w)' (   (u\nabla w)^-(\alpha)   )=   (\varphi^{-1}   +    \psi^{-1})^{-1}$.
So statement \textbf{A1} holds.

Assume that $x$, $u^-(\al)$ and $w^-(\al)$
meet
 the   premise of statement \textbf{A2}.
Since
$u'(u^-(\al))=0$
or
$w'(w^-(\al))=0$,
by Theorem
\ref{dzeron} (\romannumeral1), (\romannumeral2),
we have that
$ (u\nabla w)'(x) = (u\nabla w)' (   (u\nabla w)^-(\alpha)   )=  0$.
Hence statement \textbf{A2} holds.

Suppose that
$x\in ((u\nabla w)^-(1), (u\nabla w)^+(1) )$
or
$x\in \bigcup_{0<\al<1} (  ((u\nabla w)^-(\alpha), (u\nabla w)_s^-(\alpha) ) \cup ((u\nabla w)_s^+(\alpha), (u\nabla w)^+(\alpha) )  )$.
Then, clearly, $(u\nabla w)'(x) =0 $.
This is statement  \textbf{A9}.
\ep

\begin{re}
  {\rm
If $u'(u_s^-(\al))=:\varphi>0$ and  $w'(w_s^-(\al))=:\psi>0$, then
$u^-(\al)= u_s^-(\al)$
and
$w^-(\al)= w_s^-(\al)$. So statements \textbf{A1} and \textbf{A3} are exactly the same.
Similarly,
 statements \textbf{A5} and \textbf{A7} are same.

}
\end{re}

\begin{lm} \label{cu1r}
  Suppose that $u\in \mathcal{F} (\mathbb{R})$   and that
 $w\in \mathcal{F}_\mathrm{D}(\mathbb{R})$. If $w$ satisfies conditions (\romannumeral1) and (\romannumeral2) listed below:
\begin{description}
 \item[(\romannumeral1)] $w(w^-(0))=u(u^-(0))$ and $w(w^+(0))=u(u^+(0))$.

 \item[(\romannumeral2-1)] If $u^-(1)$ is an non-differential inner point of $[u]_0$, then
$w_-'(w^-(1))=0$.

  \item[(\romannumeral2-2)] If $u^+(1)$ is an  non-differential inner point of $[u]_0$, then
$w_+'(w^+(1))=0$.
\end{description}

Then the following statement holds.
\\
\textbf{A10}    $ (u\nabla w)' (x) =0$ for each $x\in (   (u\nabla w)^-(0),    (u\nabla w)^+(0)  ) $
with
$ (u\nabla w) (x) =1$.
\end{lm}

\pof \ Set $w(w^-(0))=u(u^-(0)):=\alpha_0$,
then by Lemma \ref{evf},
we know that
$(u\nabla w)   (   (u\nabla w)^-(0)    )   =   \alpha_0$.
Suppose that
$x=(u\nabla w)^-(1)$. Since $x$ is an inner point of
$[u\nabla w]_0$,
we
know
$\al_0 < 1$.
Now
 we prove $$(u\nabla w)' (x) =0.$$ The proof is divided into two cases.

Case (A) \  $w^-(1)$ is an inner point of $[w]_0$.

In this case, it follows from $w\in \mathcal{F}_\mathrm{D}(\mathbb{R})$ that    $w'(w^-(1))=0$.
Thus, from Theorem \ref{dzeron}(\romannumeral1)(\romannumeral2),
 $
   (u\nabla w)' (x) =  (u\nabla w)' (    (u\nabla w)^-(1)   )=0.
$

Case (B) \  $w^-(1)$ is not an inner point of $[w]_0$.

In this case $w^-(1)=w^+(0)$. So $w(w^+(0))=u(u^+(0))=1$.
Since
$x=(u\nabla w)^-(1)$ is an inner point of $[u\nabla w]_0$,
we
know
that
$u^-(0) < u^-(1) < u^+(1)=u^+(0)$. Hence $u^-(1)$
is
an inner point of $[u]_0$.
If
$u^-(1)$
is a differential point,
then
$u'(u^-(1))=0$.
So, from Theorem \ref{dzeron}(\romannumeral1)(\romannumeral2),
 $$
   (u\nabla w)' (x) =  (u\nabla w)' (    (u\nabla w)^-(1)   )=0.
$$
If
$u^-(1)$
is a non-differential point,
then
by condition (\romannumeral2-1) that
$w'_-(w^-(1))=0$,
 and thus, by Theorem \ref{dzeron}(\romannumeral1),
we
know that
\begin{equation}\label{uwrd0}
  (u\nabla w)'_- (    (u\nabla w)^-(1)   )=0.
\end{equation}
Since
 $u^-(1) < u^+(1)$, we obtain that $u'_+( u^-(1) ) = 0$, and hence by
Theorem \ref{dzeron}(\romannumeral2),
\begin{equation}\label{uwld0}
  (u\nabla w)'_+ (    (u\nabla w)^-(1)   )=0.
\end{equation}
  Combined
 with \eqref{uwrd0}
 and \eqref{uwld0}, we obtain
that
$ (u\nabla w)' (x) = (u\nabla w)' (    (u\nabla w)^-(1)   )=0$.

Similarly, we can show that  $ (u\nabla w)' (   x  )=0$ when $x= (u\nabla w)^+(1)   $.
From
statement \textbf{A9},
we know
that
 $ (u\nabla w)' (   x  )=0$
when
$ x \in (   (u\nabla w)^-(1),    (u\nabla w)^+(1)  )$.
In summary,
statement
\textbf{A10} is correct.
\ep


\begin{lm} \label{bpce}
  Suppose that $u\in \mathcal{F}(\mathbb{R})$ and $\al<1$.
\\
(\romannumeral1) \  If $ u_s^-(\al) = u^+(0)$, then $u$ is not left-continuous at $u_s^-(\al)$.
So we know that if $u^-(0)< u_s^-(\al) = u^+(0)$, then $u\notin \mathcal{F}_\mathrm{C}(\mathbb{R})$.
\\
(\romannumeral2) \  If $ u_s^+(\al) = u^-(0)$, then $u$ is not right-continuous at $u_s^+(\al)$.
So we know that if $u^-(0)= u_s^+(\al) < u^+(0)$, then $u\notin \mathcal{F}_\mathrm{C}(\mathbb{R})$.
\\
(\romannumeral3) \  If $ u^-(\al) = u^+(0)$, then $u$ is not left-continuous at $u^-(\al)$.
So we know that if $u^-(0)< u^-(\al) = u^+(0)$, then $u\notin \mathcal{F}_\mathrm{C}(\mathbb{R})$.
\\
(\romannumeral4) \  If $ u^+(\al) = u^-(0)$, then $u$ is not right-continuous at $u^+(\al)$.
So we know that if $u^-(0)= u^+(\al) < u^+(0)$, then $u\notin \mathcal{F}_\mathrm{C}(\mathbb{R})$.

\end{lm}

\pof \ Suppose that $u_s^-(\al) = u^+(0)$, $\al<1$, then $u_s^-(\al) = u^-(1) = u^+(1)= u^+(0)$.
Hence $u(   u_s^-(\al)  ) =1$.
  Note that $\lim_{x\to  u_s^-(\al)-  }    u(x) \leq  \al      $,
thus
$u$ is not left-continuous at $ u_s^-(\al)$.
If
 $u^-(0)< u_s^-(\al) = u^+(0)$, then, clearly,
$u\notin \mathcal{F}_\mathrm{C}(\mathbb{R})$.
This is statement
(\romannumeral1).
Similarly,
we can
prove
statements
(\romannumeral2)--(\romannumeral4).
\ep


\begin{lm}
  \label{s0e}
Suppose that $u\in \mathcal{F} (\mathbb{R})$
 and
that
 $w \in \mathcal{F}_\mathrm{D}(\mathbb{R})$ satisfies condition (\romannumeral1) in Lemma \ref{cu1r}.
Then,
 given $x \in (  (u\nabla w)^-(0), (u\nabla w)^+(0)   ) $ with $(u\nabla w)(x)<1$,
the following statements hold.
\begin{description}

  \item[B1] $ (u\nabla w)' (x)=0$ and $x=(u\nabla w)_s^-(0)$ when $x= (u\nabla w)_s^- (\alpha)$ and $u_s^-(\al)= u^-(0)$.

  \item[B2] $ (u\nabla w)' (x)=0$ and $x=(u\nabla w)_s^+(0)$ when $x= (u\nabla w)_s^+ (\alpha)$ and $u_s^+(\al)= u^+(0)$.

\end{description}

\end{lm}

\pof \  We only prove statement \textbf{B1}. Statement \textbf{B2} can be proved similarly.
 Set
 $w(w^-(0))=u(u^-(0))=(u\nabla w) (  (u\nabla w)^-(0)    ):=\alpha_0$.
Suppose that $x=(u\nabla w)_s^- (\alpha)$
and
 $u_s^-(\al)= u^-(0) $.
Then
$\al=\al_0< 1$.
Note
 that
$x=(u\nabla w)_s^- (\alpha_0)= u_s^-(\al_0) +   w_s^-(\al_0)$ is an inner point of $[u\nabla w]_0$,
hence
$w^-(0) = w^-(\al_0)< w_s^-(\al_0)$.
Since $w\in \mathcal{F}_\mathrm{D}(\mathbb{R})$, by Lemma \ref{bpce},
we know $w_s^-(\al_0)< w^+(0) $,
and therefore $w_s^-(\al_0)$ is an inner point of
$[w]_0$.
So
$w'(   w_s^-(\al_0)  ) =0$.
It
 thus
  follows from Theorem \ref{dzeron}(\romannumeral3), (\romannumeral4)
 that
$ (u\nabla w)' (x) = (u\nabla w)' (   (u\nabla w)_s^-(\al_0) )   = 0$.
 \ep


\begin{lm}  \label{s1e}
Suppose that $u\in \mathcal{F} (\mathbb{R})$
and
that
 $w \in \mathcal{F}_\mathrm{D}(\mathbb{R})$ satisfies condition (\romannumeral1) in Lemma \ref{cu1r} and condition (\romannumeral3) listed below
\begin{description}

\item[(\romannumeral3-1)] If $u_s^-(  u(u^-(0))     )$ is an non-differential inner point of $[u]_0$, then
$w_+'(  w^-(0)   )=0$.

\item[(\romannumeral3-2)] If $u_s^+(  u(u^+(0))     )$ is an non-differential inner point of $[u]_0$, then
$w_-'(  w^+(0)   )=0$.

\end{description}
  Then, given $x \in (  (u\nabla w)^-(0), (u\nabla w)^+(0)   ) $ with $(u\nabla w)(x)<1$,
the following statements hold.
\begin{description}
 \item[B3]  $ (u\nabla w)' (x)=0$ and $x=(u\nabla w)_s^-(0)$ when $x= (u\nabla w)_s^- (\alpha)$,
$w_s^-(\al)= w^-(0) $ and $u^-(0)<u_s^-(\al) < u^+(0) $.

\item[B4]  $ (u\nabla w)' (x)=0$ and $x=(u\nabla w)_s^+(0)$ when $x= (u\nabla w)_s^+ (\alpha)$,
$w_s^+(\al)= w^+(0) $ and $u^-(0)<u_s^+(\al) < u^+(0) $.
\end{description}

\end{lm}

\pof \ We only prove statement \textbf{B3}. Statement \textbf{B4} can be proved similarly.
 Set
 $w(w^-(0))=u(u^-(0))=(u\nabla w) (  (u\nabla w)^-(0)    ):=\alpha_0$.
From
 $w_s^-(\al)= w^-(0) $ and $u^-(0)<u_s^-(\al) < u^+(0) $, we know that $\al=\al_0<1$.

If
$u$ is differentiable at $u_s^-(\al)$,
it then follows from
$u^-(0)<u_s^-(\al_0) < u^+(0) $
that
$u'(   u_s^-(\al_0)   ) =   0$, and thus, by Theorem \ref{dzeron}(\romannumeral3), (\romannumeral4)
$$ (u\nabla w)' (x) = (u\nabla w)' (   (u\nabla w)_s^-(\al_0) )   = 0.$$
 If
$u$ is not differentiable at $u_s^-(\al)$, then from
condition (\romannumeral3-1),
we know that
$w_+'(  w_s^-(\al_0)  ) =  w_+'(  w^-(0)  )  =0$,
and
so, by Theorem \ref{dzeron}(\romannumeral4),
$$ (u\nabla w)_+' (x) = (u\nabla w)_+' (   (u\nabla w)_s^-(\al_0) )   = 0.$$
Note
that
 $x= (u\nabla w)_s^- (\alpha_0)$
is an inner point, hence
$ (u\nabla w)^-(\al_0) < (u\nabla w)_s^-(\al_0) $, and thus we know
$(u\nabla w)_-' (x) = (u\nabla w)_-' (   (u\nabla w)_s^-(\al_0) )   = 0. $
 So
$$ (u\nabla w)' (x) =(u\nabla w)' (   (u\nabla w)_s^-(\al_0) ) = 0. \quad \Box$$

The following theorem gives a method to search smoothers
for fuzzy numbers
 in
$ \mathcal{F}_\mathrm{N}(\mathbb{R})\cap \mathcal{F}_\mathrm{C}(\mathbb{R})$.

\begin{tm}
\label{rap}
Suppose that $u\in \mathcal{F}_\mathrm{N}(\mathbb{R})\cap \mathcal{F}_\mathrm{C}(\mathbb{R})$   and that
 $w\in \mathcal{F}_\mathrm{D}(\mathbb{R})$, then $w$ is a smoother of $u$, i.e.
$u\nabla
w\in \mathcal{F}_\mathrm{D}(\mathbb{R})$, when $w$ satisfies conditions (\romannumeral1) and (\romannumeral2).

\end{tm}

\begin{figure}
  \subfloat
 [ $ 0<\alpha_0<1$, $0<\beta_0<1$]
 { \label{fugal1}
\scalebox{0.33}{   \includegraphics{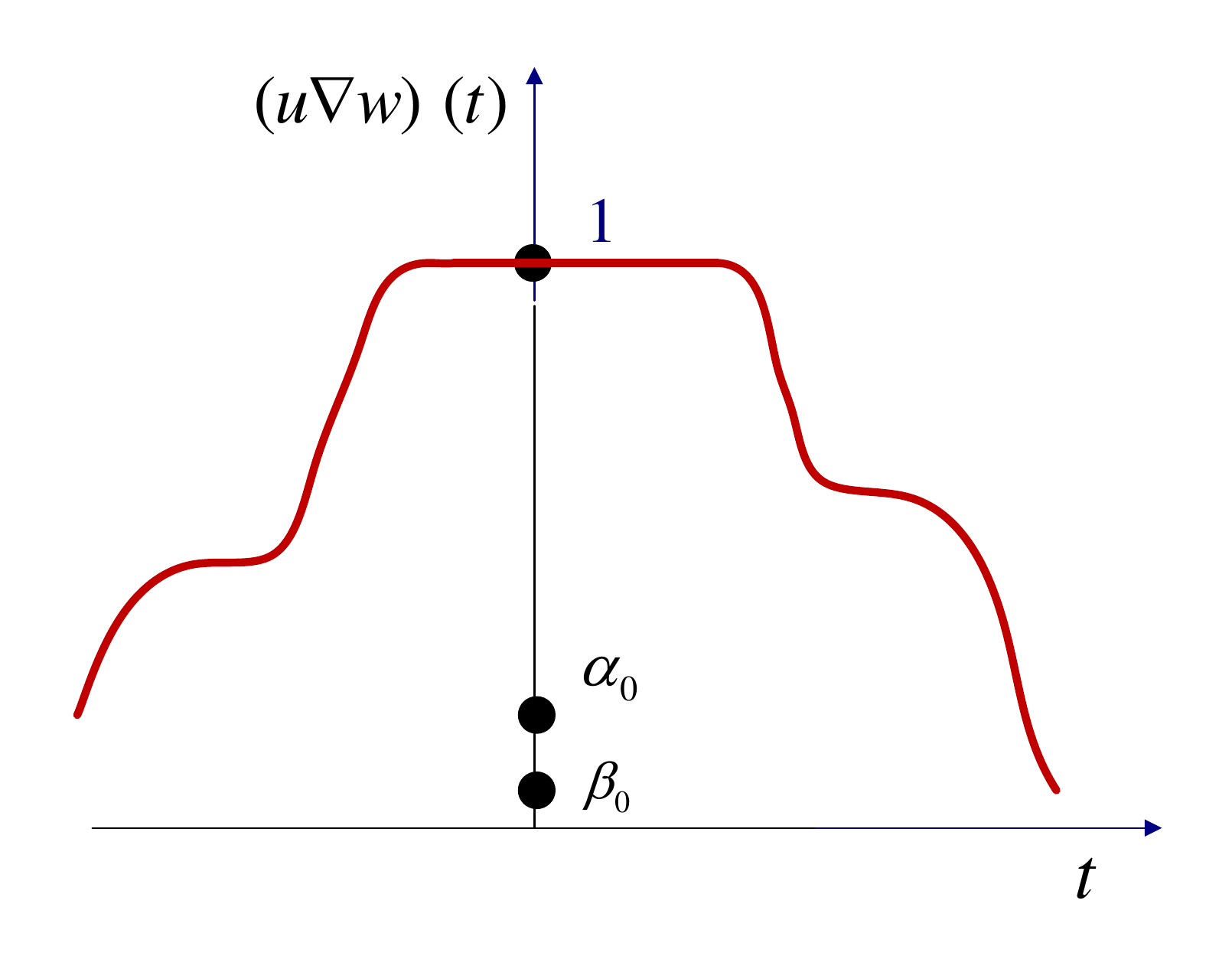}      }
}
\subfloat
[ $\alpha_0  =  \beta_0  = 0$]
 {\label{fugal0}
\scalebox{0.33}{   \includegraphics{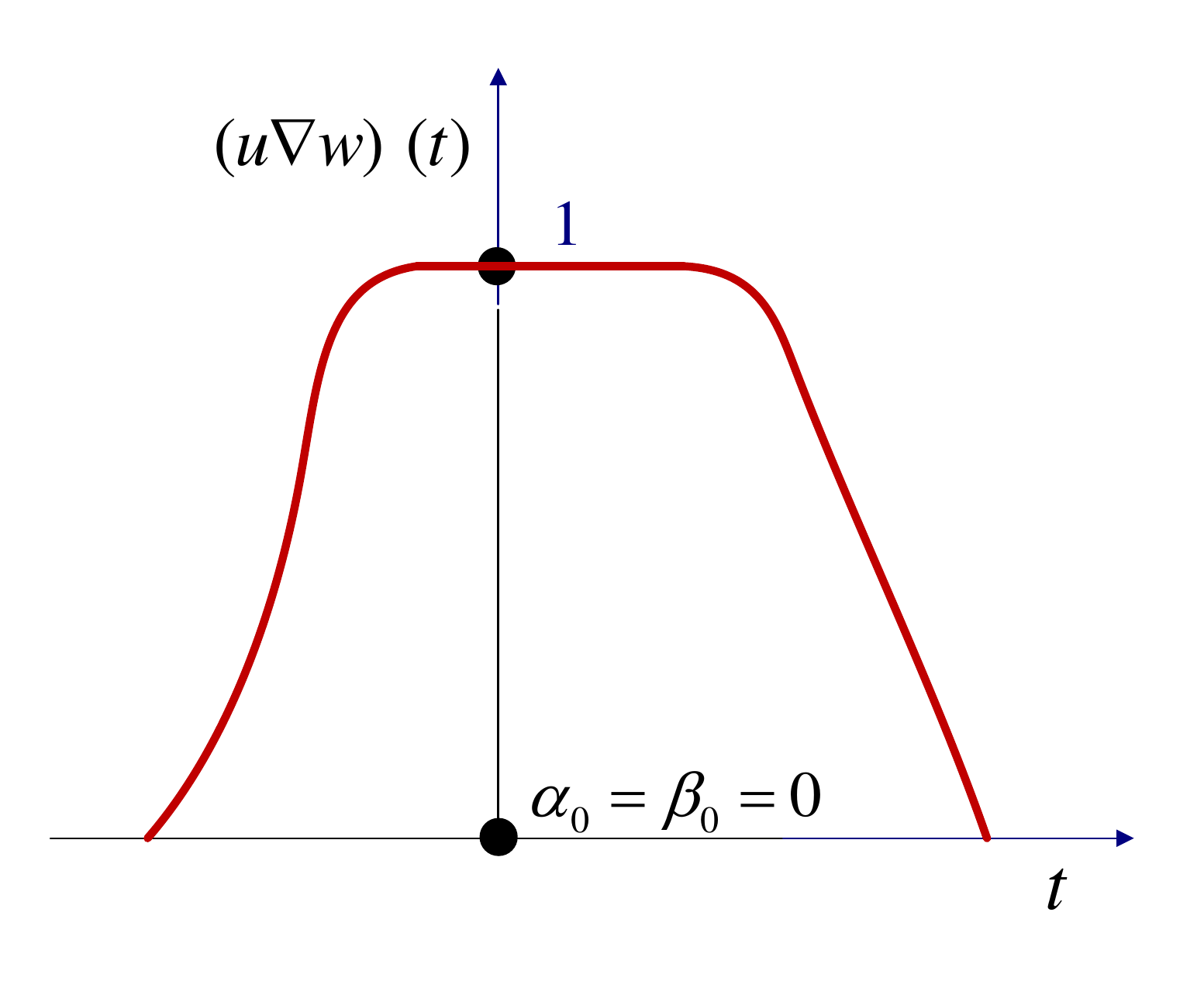}      }
}
\subfloat
[$\alpha_0=\beta_0=1$]
 { \label{fugde1}
\scalebox{0.33}{   \includegraphics{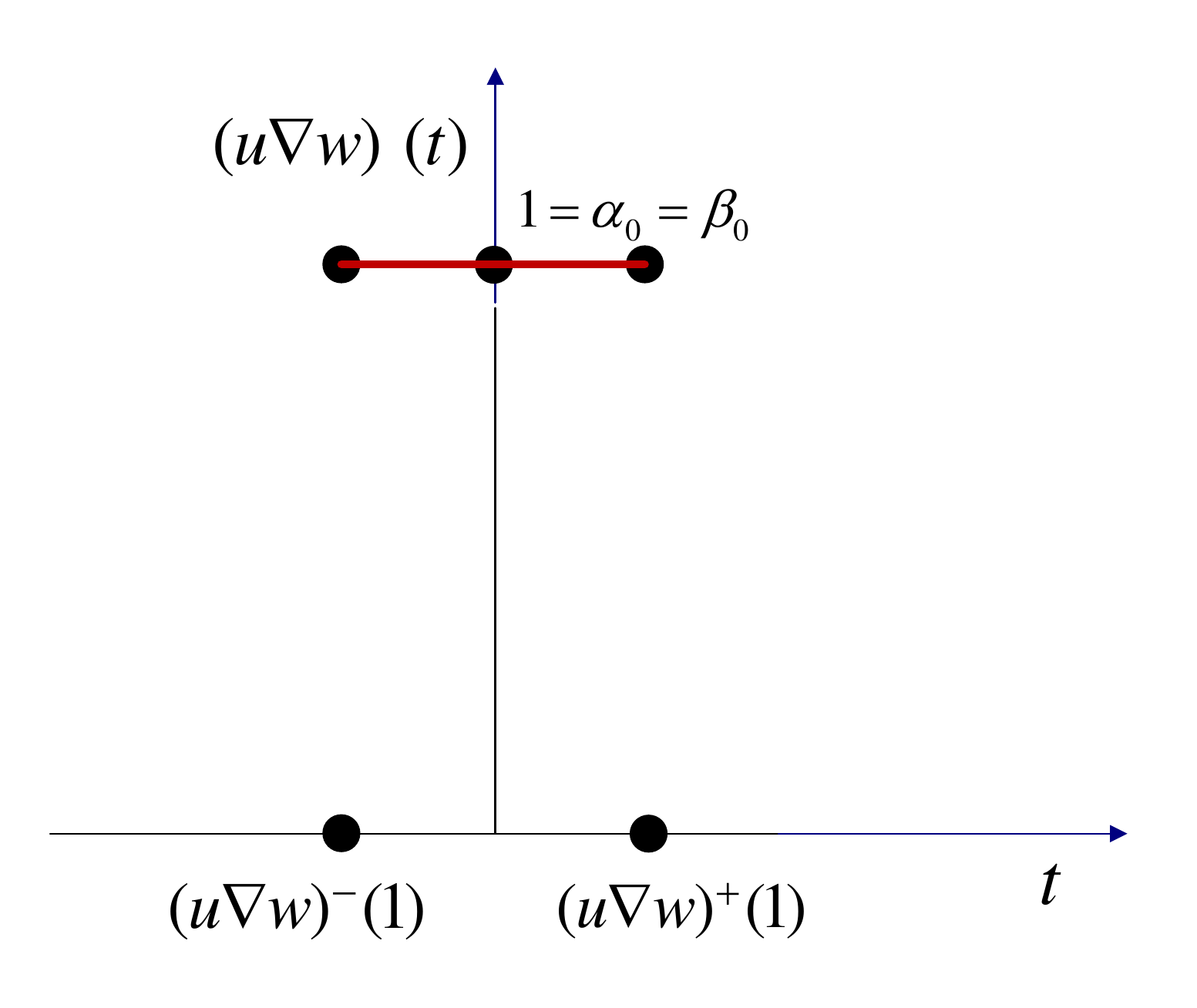}      }
}

 \caption{Some examples of $u\nabla w$  }
 \label{fuga}
 \end{figure}

\pof \ Set $w(w^-(0))=u(u^-(0)):=\alpha_0$ and $w(w^+(0))=u(u^+(0)):=\beta_0$,
then by Lemma \ref{evf},
we know that
$(u\nabla w)   (   (u\nabla w)^-(0)    )   =   \alpha_0$
and
$(u\nabla w)   (   (u\nabla w)^+(0)    )   =   \beta_0$.
Thus
$ [ u\nabla w ]_0 = [ (u\nabla w)^-(\alpha_0),   (u\nabla w)^+(\beta_0)    ] $.
In
 Fig \ref{fuga}, we give some examples of $u\nabla w$ with different $\al_0$ and $\beta_0$.

Since $w$ satisfies condition (\romannumeral1),
from Theorem \ref{arp},
we know
that
$u\nabla w   \in \mathcal{F}_\mathrm{C}(\mathbb{R})$.
To show
$u\nabla
w\in \mathcal{F}_\mathrm{D}(\mathbb{R})$
 is equivalent to
prove
that
$u\nabla w$ is differentiable
at
each
inner point
of
$[ u\nabla w ]_0 $.

Let $x$ be an inner point of $[ u\nabla w ]_0 $,
i.e.
$ x \in ( (u\nabla w)^-(\alpha_0),   (u\nabla w)^+(\beta_0)    )$.
From statement \textbf{A10}, we know that $(u\nabla w)'(x) =0$  when $(u\nabla w)(x) =1$.
By
statement \textbf{A9},
we obtain that
$(u\nabla w)'(x) =0$
when
 $x$ is neither an endpoint of an $\al$-cut nor an endpoint of a strong-$\al$-cut.

Next,
we consider the rest situation of $x$, i.e. $(u\nabla w)(x) < 1$
and
 $x$ is an endpoint of an $\al$-cut or an endpoint of a strong-$\al$-cut.
We split this situation into two cases.
It is clear that
$\al<1$ in these two cases.

Case (A)\ $x= (u\nabla w)^-(\al)$ ($(u\nabla w)_s^-(\al)$, $(u\nabla w)^+(\al)$,  $(u\nabla w)_s^+(\al)$)
with $(u\nabla w) (x)<1$,
and both
 $u^-(\al)$ and $w^-(\al)$ ($u_s^-(\al)$ and $w_s^-(\al)$,  $u^+(\al)$ and $w^+(\al)$, $u_s^+(\al)$ and $w_s^+(\al)$) being inner points of $[u]_0$ and $[w]_0$, respectively.

Note that  $u\in \mathcal{F}_\mathrm{N}(\mathbb{R})\cap \mathcal{F}_\mathrm{C}(\mathbb{R})$
and
 $w\in \mathcal{F}_\mathrm{D}(\mathbb{R})$,
hence
$w$ is differentiable at each inner point
and
$u$ is differentiable at each inner point $y$ when $u(y)<1$.
So,
by statements \textbf{A1}--\textbf{A8},
we can compute
$(u\nabla w)'(x)$.
For example,
suppose that $x= (u\nabla w)^-(\al)= u^-(\al) + w^-(\al)$, and
both
$u^-(\al)$ and $w^-(\al)$
are
 inner points of $[u]_0$ and $[w]_0$,
then,
by statements \textbf{A1} and \textbf{A2},
we
obtain
the differentiability
of $u\nabla w$
at
$x$.

Case (B)\ $x= (u\nabla w)^-(\al)$ ($(u\nabla w)_s^-(\al)$, $(u\nabla w)^+(\al)$, $(u\nabla w)_s^+(\al)$)
 with
$(u\nabla w)(x) < 1$
and $x$ is not in
Case (A).

We claim that the following situations will not occur.
\begin{description}
 \item[\romannumeral1]  $x=(u\nabla w)_s^-(\al)$, $u_s^-(\al) = u^-(0)$ and $w_s^-(\al) = w^-(0)$.
 \item[\romannumeral2]  $x=(u\nabla w)_s^+(\al)$, $u_s^+(\al) = u^+(0)$ and $w_s^+(\al) = w^+(0)$.
\item[\romannumeral3] $w^-(0) < w_s^-(\al) = w^+(0)$.
\item[\romannumeral4] $w^-(0) = w_s^+(\al) < w^+(0)$.
\end{description}
In fact, if $x=(u\nabla w)_s^-(\al)$, $u_s^-(\al) = u^-(0)$ and $w_s^-(\al) = w^-(0)$,
then
$x=(u\nabla w)_s^-(\al)= u^-(0) +  w^-(0) = (u\nabla w)^-(0) $
is not an inner point of $[u\nabla w]_0$.
Hence   situation (\textbf{\romannumeral1})
can
not
occur.
Similarly, situation (\textbf{\romannumeral2})
does
not
occur.
Since
$w \in \mathcal{F}_\mathrm{C}(\mathbb{R})$, by Lemma \ref{bpce},
we know
that
situations (\textbf{\romannumeral3}) and (\textbf{\romannumeral4}) will not occur.

If
 $x = (u\nabla w)^-(\al)$ is an inner point of $[u \nabla w]_0$,
we
obtain
$1>\al>\al_0$.
Note that
 $w$ is in $\mathcal{F}_\mathrm{D}(\mathbb{R})$,
hence
$w^-(\al_0)< w^-(\al) <w^-(1)$,
and therefore $u^-(\al_0)< u^-(\al) = u^+(0)$.
So case (B) can be divided into the following subcases.
\begin{description}
  \item[B\romannumeral1] $x= (u\nabla w)_s^- (\alpha)$ and $u_s^-(\al)= u^-(0)$.

  \item[B\romannumeral2] $x= (u\nabla w)_s^+ (\alpha)$ and $u_s^+(\al)= u^+(0)$.

 \item[B\romannumeral3] $x= (u\nabla w)_s^- (\alpha)$,
$w_s^-(\al)= w^-(0) $ and $u^-(0)<u_s^-(\al) < u^+(0) $.

\item[B\romannumeral4]  $x= (u\nabla w)_s^+ (\alpha)$,
$w_s^+(\al)= w^+(0) $ and $u^-(0) < u_s^+(\al) < u^+(0) $.

\item[B\romannumeral5] $x= (u\nabla w)_s^- (\alpha)$,
$w_s^-(\al)= w^-(0) $ and $u^-(0) < u_s^-(\al) = u^+(0) $.

\item[B\romannumeral6]  $x= (u\nabla w)_s^+ (\alpha)$,
$w_s^+(\al)= w^+(0) $ and $u^-(0)   =    u_s^+(\al) < u^+(0) $.

\item[B\romannumeral7] $x= (u\nabla w)_s^- (\alpha)$,
$ w^-(0)< w_s^-(\al) < w^+(0) $ and $u^-(0) < u_s^-(\al) = u^+(0) $.

\item[B\romannumeral8]  $x= (u\nabla w)_s^+ (\alpha)$,
$w^-(0) < w_s^+(\al) < w^+(0) $ and $u^-(0)   =    u_s^+(\al) < u^+(0) $.

 \item[B\romannumeral9]    $x = (u\nabla w)^-(\al)$,   $w^-(0)< w^-(\al) <w^-(1)$
and
 $u^-(0)< u^-(\al) = u^+(0)$.

  \item[B\romannumeral10]   $x = (u\nabla w)^+(\al)$,
$w^+(1) < w^+(\al) <w^+(0)$
and
$u^-(0)= u^+(\al) < u^+(0)$.
\end{description}
Since
$u\in \mathcal{F}_\mathrm{C}(\mathbb{R})$,   by Lemma \ref{bpce},
we
know
that
subcases \textbf{B\romannumeral5}--\textbf{B\romannumeral10} do not occur,
i.e.,
case (B)
 can be decomposed into
four subcases: \textbf{B\romannumeral1}--\textbf{B\romannumeral4}.
Notice that if $u \in \mathcal{F}_\mathrm{N}(\mathbb{R})$,
then
$w$ satisfies condition (\romannumeral3).
So,
by statements \textbf{B1}--\textbf{B4},
we can prove the differentiability of $u\nabla w$
at
$x$. \ep

Now we give some
 examples     to illustrate  Theorem \ref{rap}.
By
Example \ref{eapc1m},
we describe the role of
 condition (\romannumeral1) in Theorem \ref{rap}.

\begin{eap} \label{eapc1m}    {\rm

\begin{figure}
  \subfloat
 [$u$]
 { \label{eapc1u}
\scalebox{0.24}{   \includegraphics{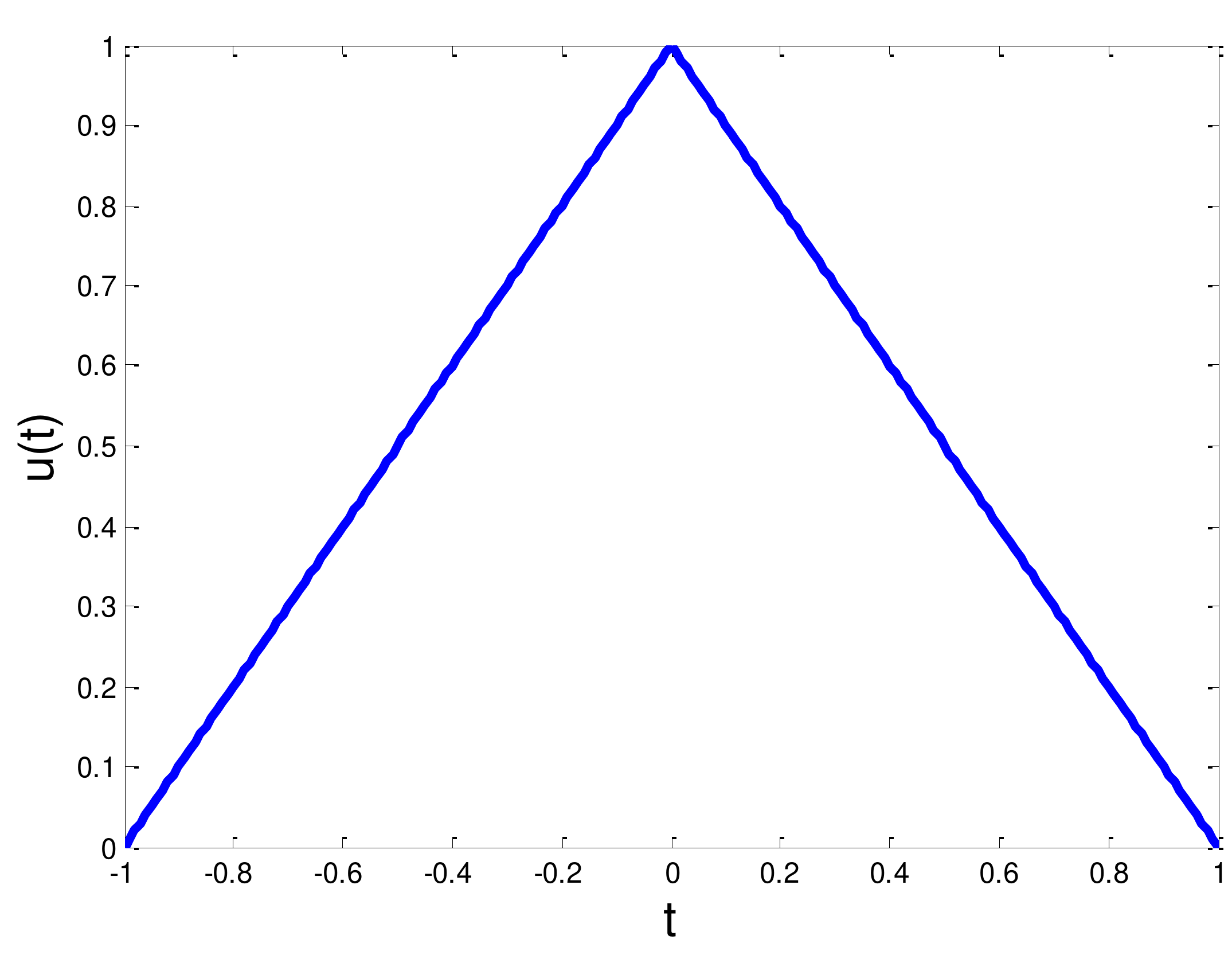}      }
}
\subfloat
 [$w$]
 { \label{eapc1w}
\scalebox{0.24}{   \includegraphics{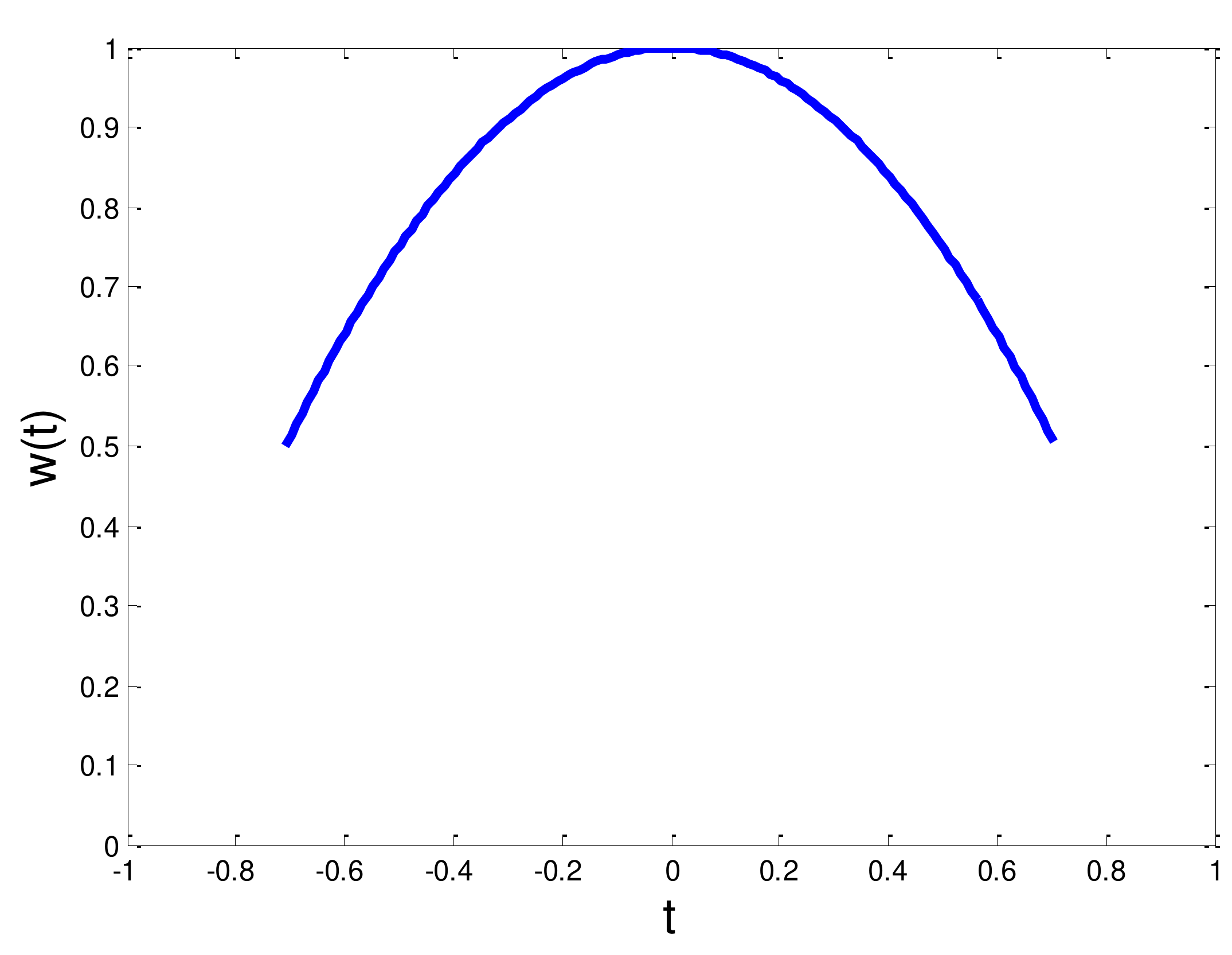}      }
}
 \subfloat
 [$u\nabla w$]
 { \label{eapc1uw}
\scalebox{0.24}{   \includegraphics{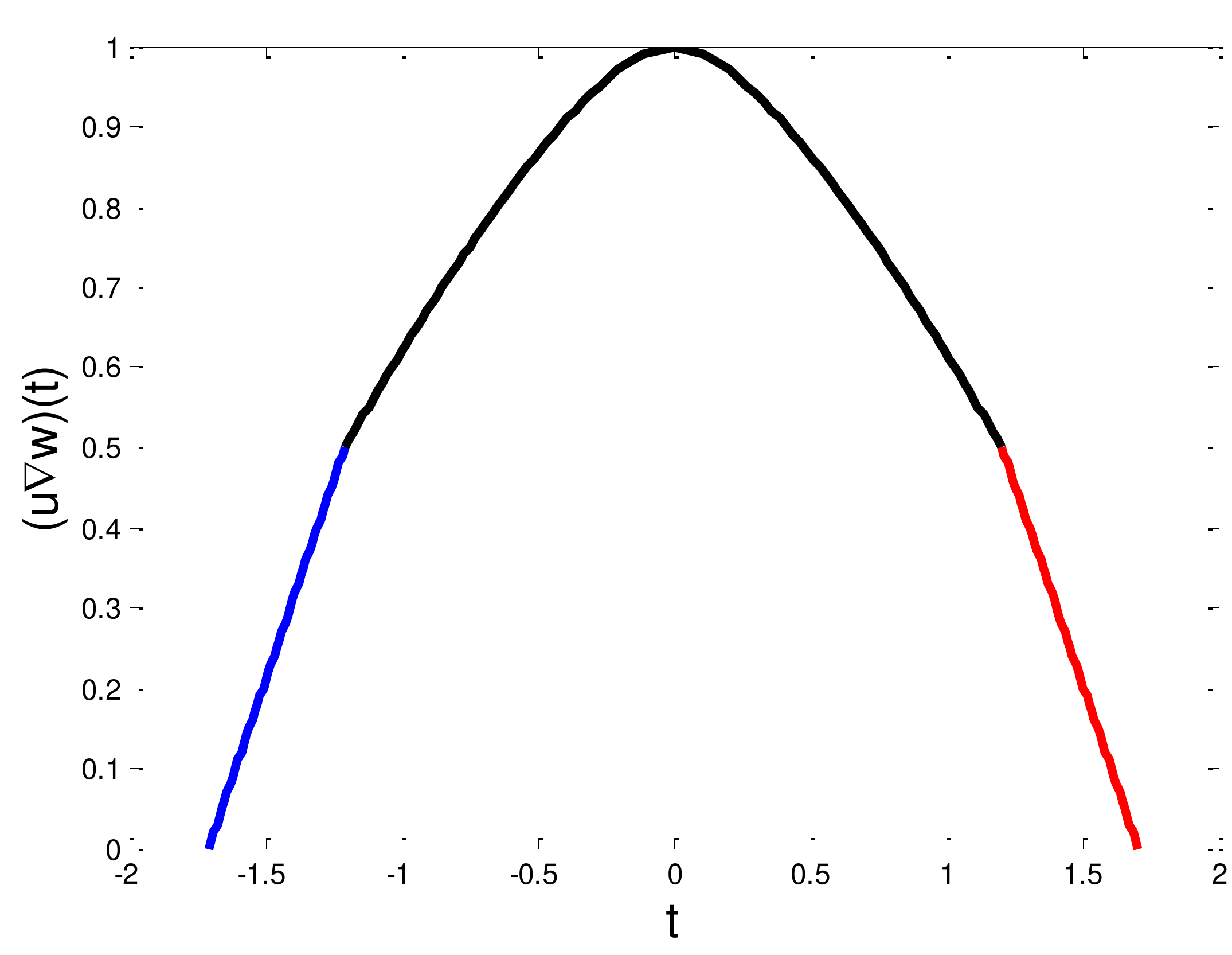}      }
}

 \caption{Fuzzy numbers in Example  \ref{eapc1m} }
 \label{eapc1f}
 \end{figure}

Suppose
\[
u(t) =
\left\{
  \begin{array}{ll}
   t+1 , &  t\in [-1,0], \\
    1-t, & t\in [0,1], \\
    0, &  t\notin [-1,1].
  \end{array}
\right.
\]
Observe that $u$ is differentiable on $(-1,1) \setminus \{0\}$.
We can see that $u\in \mathcal{F}_\mathrm{N}(\mathbb{R})\cap \mathcal{F}_\mathrm{C}(\mathbb{R})$
and
\[
[u]_\al =
   [-1+\alpha, \   1-\alpha]
\]
for all  $\alpha\in [0,1]$. See Fig. \ref{eapc1u} for the figure of $u$.
Let
\[
w(t) =
\left\{
  \begin{array}{ll}
   1-t^2 , &  t\in [-\sqrt{0.5},    \sqrt{0.5}], \\
       0, &  t\notin [-\sqrt{0.5},    \sqrt{0.5}],
  \end{array}
\right.
\]
then
\[
[w]_\al =
\left\{
  \begin{array}{ll}
   [-\sqrt{1-\alpha},         \sqrt{1-\alpha} ],    &  \alpha\in [0.5,    1], \\
   \mbox{}  [-\sqrt{0.5},      \sqrt{0.5}],         &  \alpha\notin [0,    0.5].
  \end{array}
\right.
\]
See Fig. \ref{eapc1w} for the figure of $w$.
It can be checked that
$w$ is differentiable on $(-0.5, 0.5)$,
i.e.
$w \in \mathcal{F}_\mathrm{D}(\mathbb{R})$.
Notice that
\begin{gather*}
  w'(w^-(1))  = w'(w^+(1)) =w'(0)=0,
 \\
  w(w^-(0))=0.5 \not= 0= u(u^-(0),
\\
 w(w^+(0))=0.5 \not= 0= u(u^+(0)),
\end{gather*}
thus we know that
$u$ and $w$ satisfy all but condition (\romannumeral1)
in Theorem \ref{rap}.
We will see that $w$ is not a smoother of $u$.
In
fact,
it can be
computed that
\[
[u\nabla w]_\al =[u]_\al  +   [w]_\al =
\left\{
  \begin{array}{ll}
   [-1+\alpha-\sqrt{1-\alpha},  \     1-\alpha+ \sqrt{1-\alpha} ],    &  \alpha\in [0.5,    1], \\
   \mbox{}  [-1+\alpha-\sqrt{0.5},  \    1-\alpha+  \sqrt{0.5}],     &  \alpha\notin [0,    0.5].
  \end{array}
\right.
\]
Now we can plot the figure of  $u\nabla w$ in Fig \ref{eapc1uw}, and obtain
\[
(u\nabla w) (t) =
\left\{
  \begin{array}{ll}
  1+t+\sqrt{0.5},    &  t \in [-1-\sqrt{0.5},   -0.5-\sqrt{0.5}],
\\
  t + 0.5\sqrt{1 - 4t} + 0.5 ,    &  t \in [-0.5-\sqrt{0.5},   0],
\\
 -t + 0.5\sqrt{1 + 4t} + 0.5 ,    &  t \in [0,   0.5+\sqrt{0.5}],
\\
 1-t+\sqrt{0.5},    &  t \in [0.5+\sqrt{0.5},   1+\sqrt{0.5}].
  \end{array}
\right.
\]
 It
is easy to check that $u\nabla w$ is not differentiable at $t=-0.5-\sqrt{0.5}$ and $t=0.5+\sqrt{0.5}$,
and therefore $u\nabla w \notin \mathcal{F}_\mathrm{D}(\mathbb{R})$.
So $w$ is not a smoother of $u$.
This
means
that
condition (\romannumeral1) in Theorem \ref{rap}
can not be omitted.

}
\end{eap}

%
%
%
%
%

In Proposition \ref{1}, it
 pointed out that if $u\in \mathcal{F}_\mathrm{T}(\mathbb{R})$,
  then
$u\nabla w_p \in \mathcal{F}_\mathrm{D}(\mathbb{R})$.
This means
that
$w_p$, $p>0$, can serve as smoothers for all fuzzy numbers in
$\mathcal{F}_\mathrm{T}(\mathbb{R})$.
In Theorem \ref{rap}, it finds
 that
$w_p$, $p>0$, can also work as smoothers
even for
fuzzy numbers not in $\mathcal{F}_\mathrm{T}(\mathbb{R})$.
The
following
  example is given to show this fact.

\begin{eap}  \label{unsm} {\rm

Suppose
\[ u(t)=\left\{
\begin{array}{ll}
-0.5(t^{2}+2t),
  &   \ t\in [-2,-1),\\
 0.5,   &  \ t\in[-1,   -0.5],\\
  2t^{2}+2t+1,    & \ t\in (-0.5,0],\\
  2t^{2}-2t+1,   & \ t\in (0,0.5),\\
0.5,    &  \ t\in [0.5,1],\\
-0.5(t^{2}-2t),    & \ t\in(1,2],\\
0,                 &  \ t\not\in [-2,2].
\end{array}
\right.
\]
See
 Fig. \ref{eapfug}
 for the figure of $u$.
\begin{figure}

 \subfloat
 [$u$]
 { \label{eapfug}
\scalebox{0.24}{   \includegraphics{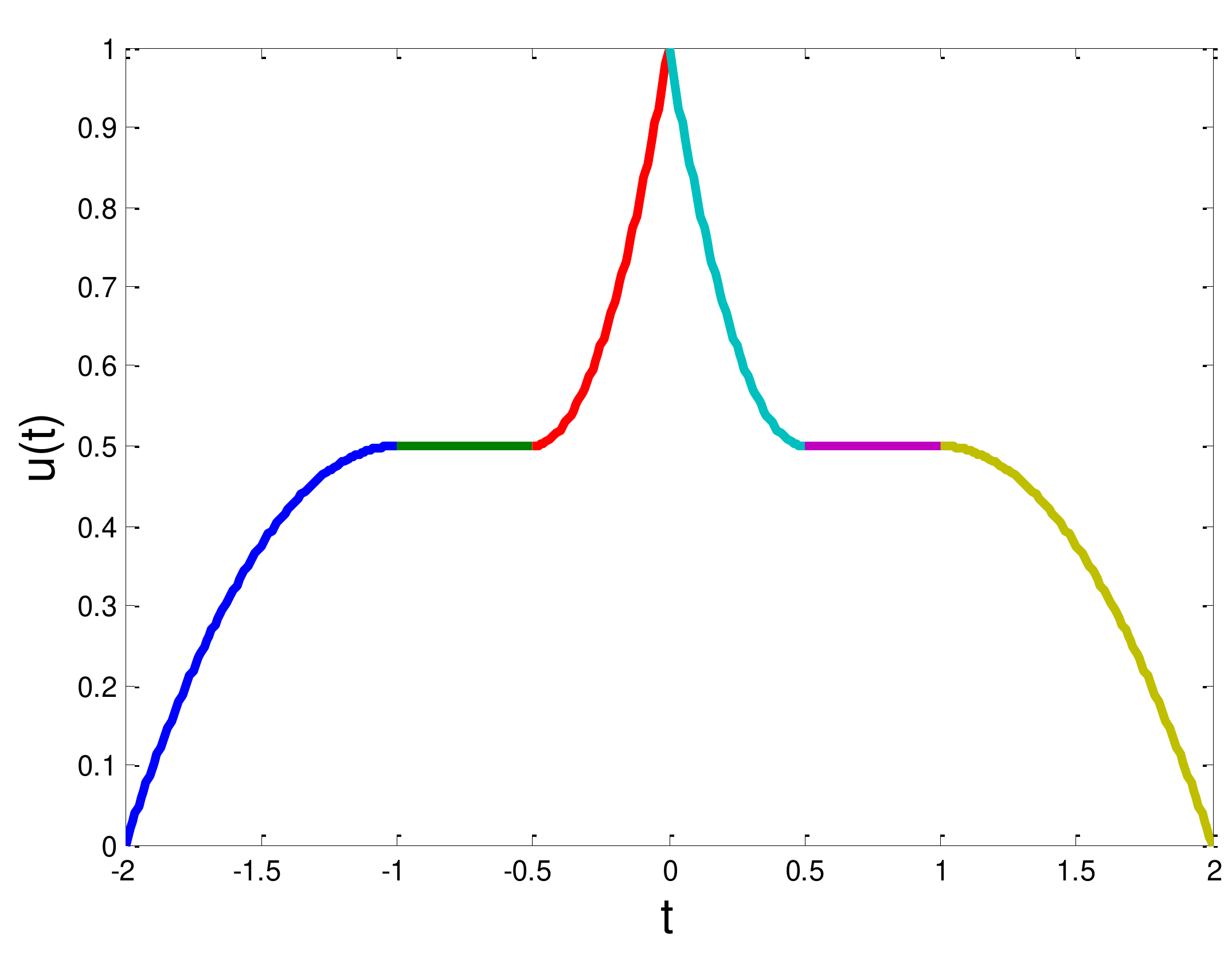}      }
}
\subfloat
 [$w_1$]
 { \label{eapfwg}
\scalebox{0.24}{   \includegraphics{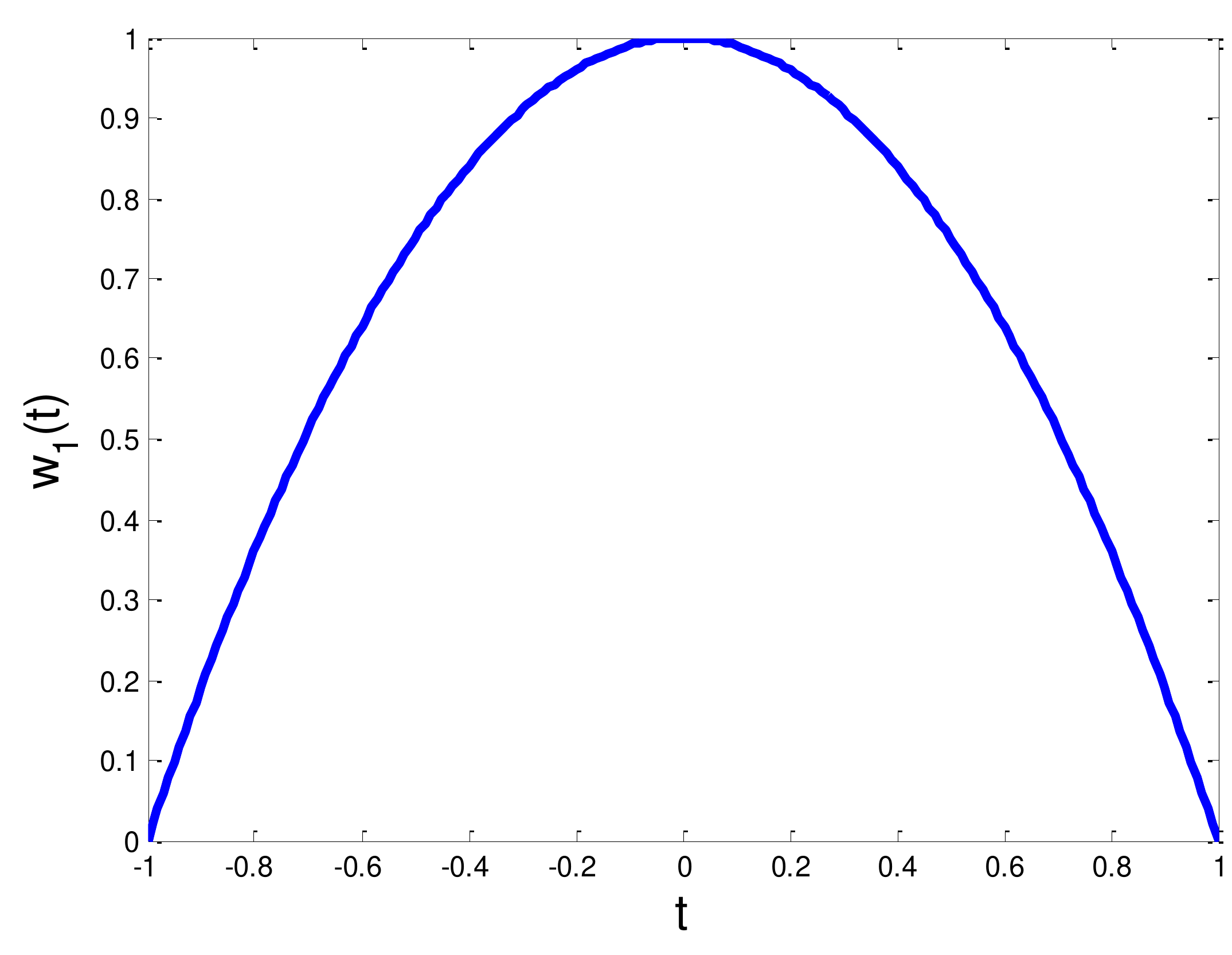}      }
}
 \subfloat
 [$u\nabla w_1$]
 { \label{eapfuwg}
\scalebox{0.24}{   \includegraphics{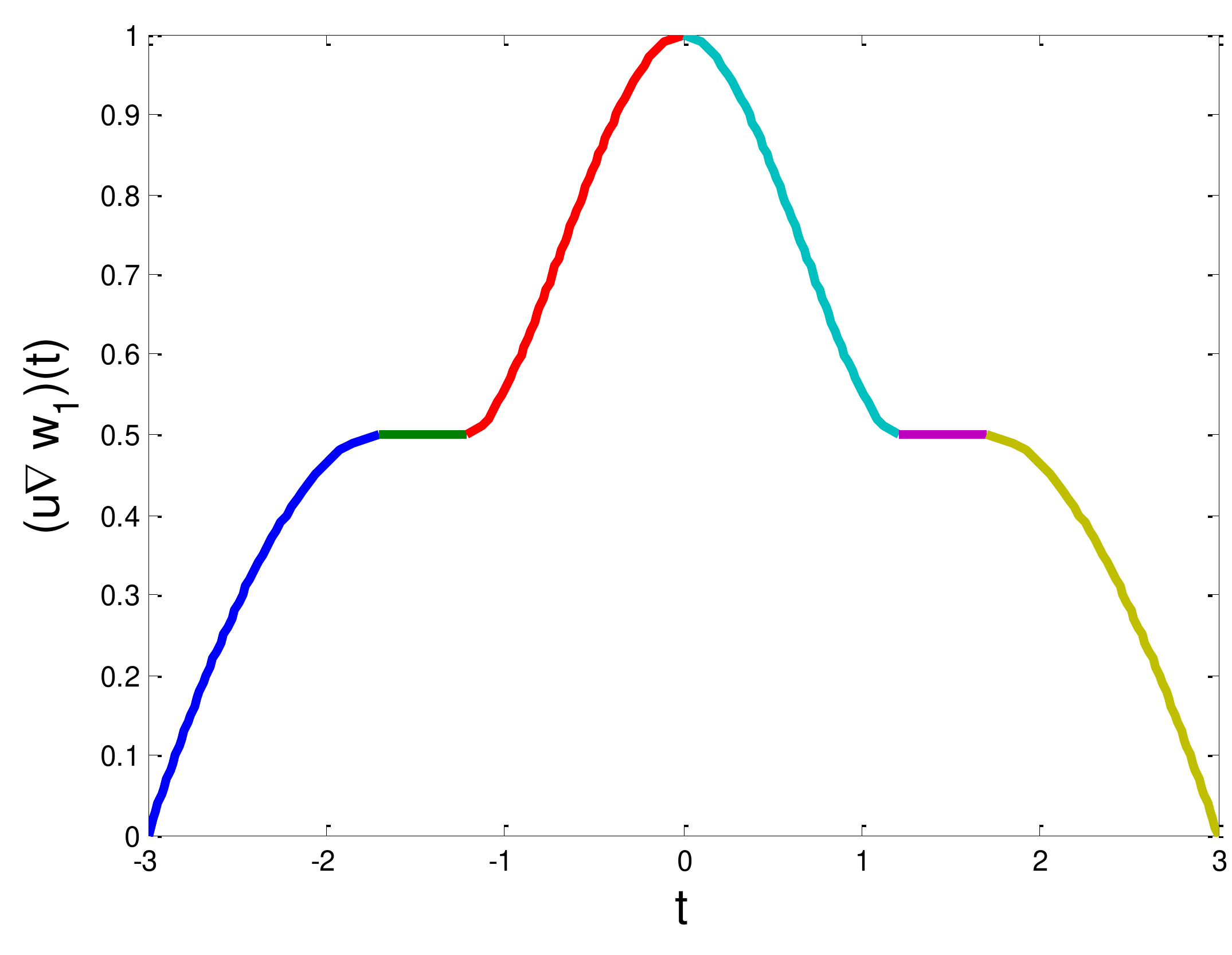}      }
}

 \caption{Fuzzy numbers $u$,  $w_1$ and $u\nabla w_1$ in Example \ref{unsm} }

 \label{eapf}

 \end{figure}
 Clearly
$u$ is not
strictly increasing on
$[-1,-0.5]$ and is not strictly decreasing on
$[0.5,1]$. Therefore $u$ is in
$\mathcal{F}_\mathrm{N}(\mathbb{R}) \cap \mathcal{F}_\mathrm{C}(\mathbb{R}) $ but     not in
$\mathcal{F}_\mathrm{T}(\mathbb{R})$.
Take
$w_1$ defined in \eqref{wpc}, i.e.
\[
w_1(t)=\left\{
\begin{array}{ll}
1-t^2,     &\ t\in [-1,1],\\
0,        & \ t\notin [-1,1].
\end{array}
\right.
\]
The figure of $w_1$ is in Fig \ref{eapfwg}.
It is easy to check that $u$ and $w_1$ satisfy all the conditions in
Theorem \ref{rap}.
Thus $w_1$ is a smoother for $u$.
Now, we validate this assertion by computing $u\nabla w_1$.
Note that
\begin{align*}
  [u & \nabla w_{1}]_\al=[u]_\al +[ w_{1}]_\al  \\
   &= \left\{
\begin{array}{ll}
[-\sqrt{1-2\alpha }-1-\sqrt{1-\al}, \ 1+ \sqrt{1-2\alpha} +\sqrt{1-\al}\, ],
& \  \alpha\in [0,0.5],
\\
\mbox{} [ 0.5( -1 +  \sqrt{ 2\alpha-1 }) -\sqrt{1-\al} , \ 0.5( 1  - \sqrt{ 2\alpha-1 } )   +   \sqrt{1-\al}   ],
& \    \alpha\in (0.5,1],
\end{array}
\right.
\end{align*}
and
 thus we can obtain that
\begin{align*}
  (u\nabla & w_{1})(t)   \\
  = & \left\{
\begin{array}{ll}
-(t+1)(    3t  +   2(2t^2+4t+3)^{0.5}   +   3    ),
&
 \hspace{-3cm}  t\in [-3,-1-\sqrt{0.5}],
\\
0.5,    & \hspace{-3cm}    t\in (-1-\sqrt{0.5},-0.5-\sqrt{0.5}],
\\
(4t(1 - 2t - 2t^2)^{0.5})/9 - (2t)/9
 + (2(1 - 2t - 2t^2)^{0.5})/9 - (2t^2)/9 + 7/9,
 \\
&  \hspace{-3cm}
t\in(-0.5-\sqrt{0.5},0],
\\
(-4t(1 + 2t - 2t^2)^{0.5})/9 + (2t)/9
 + (2(1 + 2t - 2t^2)^{0.5})/9 - (2t^2)/9 + 7/9,
 \\
&
\hspace{-3cm}
t\in[0,0.5+\sqrt{0.5}),
 \\
0.5,
&
\hspace{-3cm}  t\in [0.5+\sqrt{0.5},1+\sqrt{0.5}),
\\
(t-1)(   -3t + 2(2t^2-4t+3)^{0.5} + 3   ),
&
\hspace{-3cm}  t\in [1+\sqrt{0.5},3],
\\
0,
&
\hspace{-3cm} t\not\in [-3,3].
\end{array}
\right.
\end{align*}
The
figure of $u\nabla w_1$ is in Fig \ref{eapfuwg}.
It
 can be computed that
  \begin{gather*}
 (u\nabla w_{1})'(-1-\sqrt{0.5})= (u\nabla w_{1})'(-0.5-\sqrt{0.5})=0,
  \\
(u\nabla w_{1})'(0)=0,
\\
(u\nabla w_{1})'(1+\sqrt{0.5})=(u\nabla w_{1})'(0.5+\sqrt{0.5})=0.
\end{gather*}
Combined
with the expression of $u\nabla w_1$,
it now follows that
$u\nabla w_{1}$ is differentiable on $(-3,3)$, i.e., $u\nabla
w_{1}\in \mathcal{F}_\mathrm{D}(\mathbb{R})$.
This means that $w_1$ can work as a smoother for the fuzzy number
$u$ which is in
$(      \mathcal{F}_\mathrm{N}(\mathbb{R}) \cap \mathcal{F}_\mathrm{C}(\mathbb{R})     ) \backslash
\mathcal{F}_\mathrm{T}(\mathbb{R})$.

}
\end{eap}

However $w_p$, $p>0$, may not be smoothers for a fuzzy number $u$
which is not in $\mathcal{F}_\mathrm{C}(\mathbb{R})$.
The following example is given to show this fact.

\begin{eap} \label{eapcn}
  {\rm

\begin{figure}

 \subfloat
 [$u$]
 { \label{eapcnu}
\scalebox{0.31}{   \includegraphics{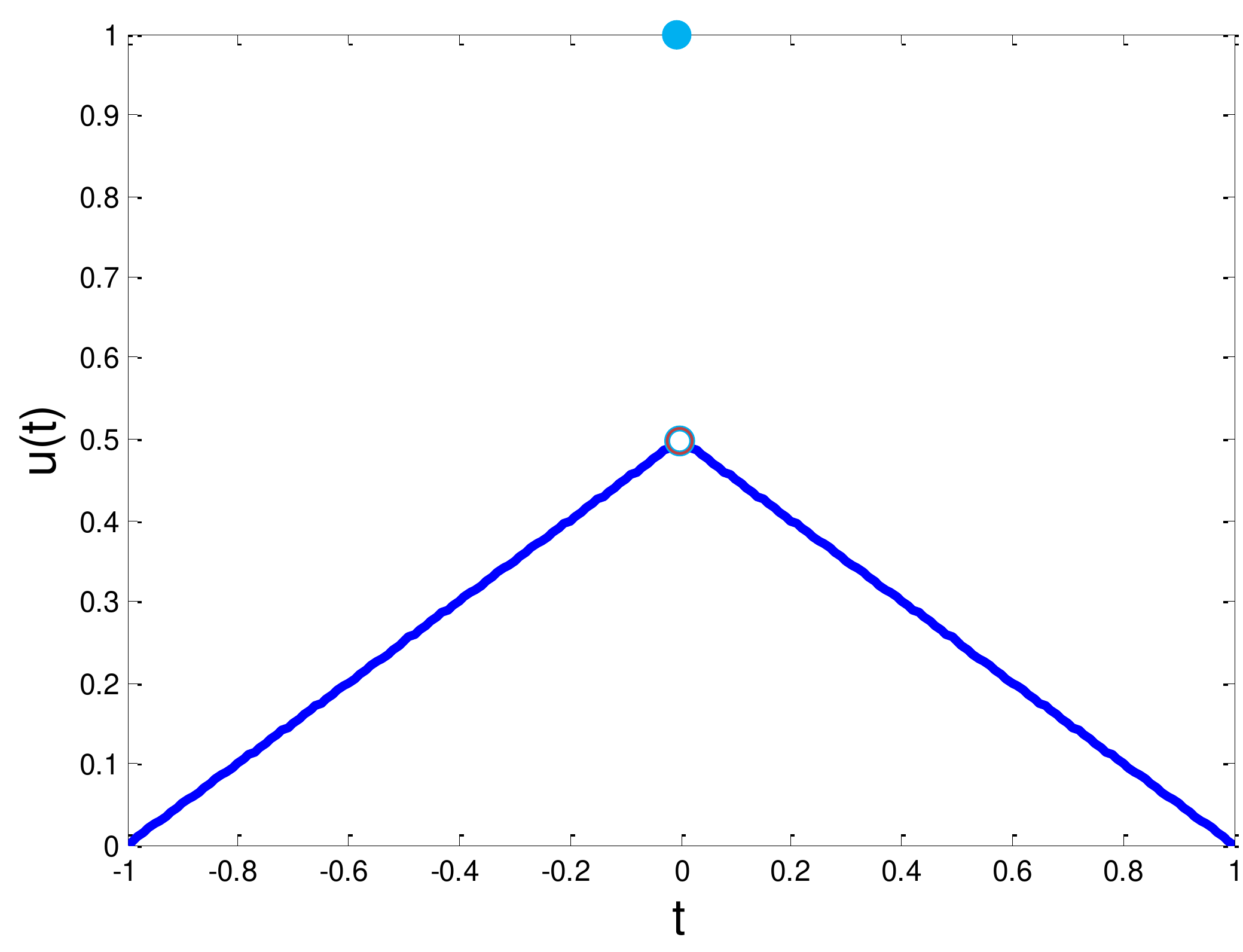}      }
}
\subfloat
 [$u\nabla w_1$]
 { \label{eapcnuw}
\scalebox{0.3}{   \includegraphics{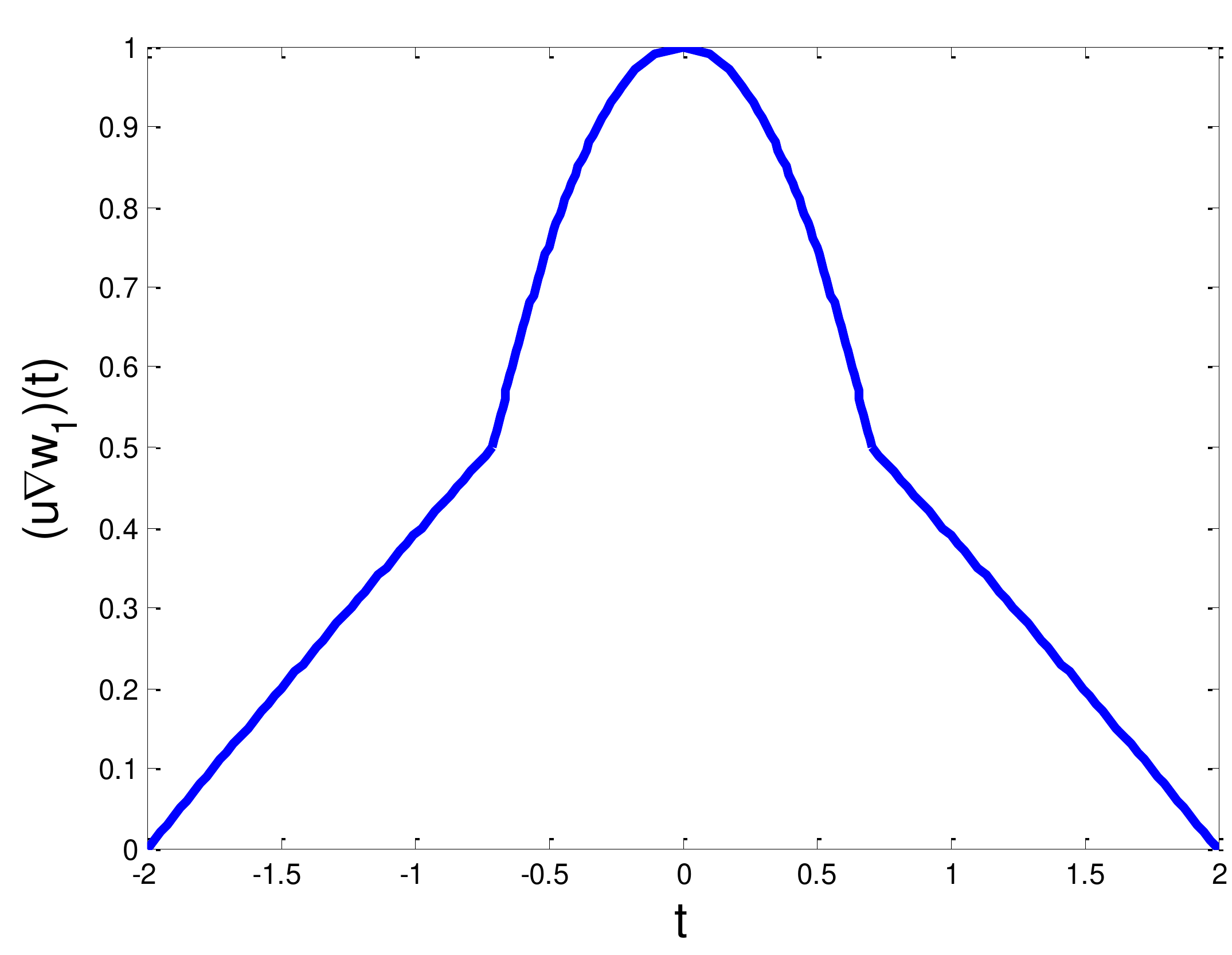}      }
}

 \caption{Fuzzy numbers in Example \ref{eapcn} }

 \label{eapcnf}

 \end{figure}

Suppose that
\[
u(t)=\left\{
\begin{array}{ll}
0.5t+0.5,    &\ t\in [-1,   0),
\\
-0.5t+0.5,   & \ t\in (0,  1],
\\
1,    &  \ t=0,
\\
0,    & \   t \notin [-1,   1].
\end{array}
\right.
\]
See Fig. \ref{eapcnu} for the figure of $u$.
It
 is easy to check that
$u\in ( \mathcal{F}_\mathrm{T}(\mathbb{R}) \cap \mathcal{F}_\mathrm{N}(\mathbb{R}) ) \backslash \mathcal{F}_\mathrm{C}(\mathbb{R})$.
But
\[
(u\nabla w_1)(t)
=
\left\{
\begin{array}{ll}
 0.125(-8t + 9)^{0.5} + 0.5t + 0.375, &   \   t\in [-2,  -\sqrt{0.5}],
\\
1-t^2,   & \ t\in [-\sqrt{0.5},  \sqrt{0.5}],
\\
 0.125(8t + 9)^{0.5} - 0.5t + 0.375, &   \   t\in [\sqrt{0.5},    2].
\end{array}
\right.
\]
See Fig. \ref{eapcnuw} for the figure of $u\nabla w_1$.
It is easy to check that
$u\nabla w_1$ is not differentiable at $t=\pm \sqrt{0.5}$.
Thus $u\nabla w_1 \notin \mathcal{F}_\mathrm{D}(\mathbb{R})$.
So $w_1$ is not a smoother for
$u$.
This means that $w_p$, $p>0$,
     may not work as smoothers for
fuzzy numbers in $\mathcal{F}_\mathrm{T}(\mathbb{R}) \cap \mathcal{F}_\mathrm{N}(\mathbb{R})  $
but
 not in
$\mathcal{F}_\mathrm{C}(\mathbb{R})$.

}
\end{eap}

\begin{re}\label{usc}
{\rm

Define
$$
\mathcal{F}^0 _\mathrm{NC}(\mathbb{R}) =
 \{
u\in \mathcal{F}_\mathrm{N}(\mathbb{R})\cap \mathcal{F}_\mathrm{C}(\mathbb{R}) : u( u^-(0) )   = u( u^+(0) ) =0 \}.
$$
Let $u$ be
a fuzzy number in $\mathcal{F}^0 _\mathrm{NC}(\mathbb{R})$.
Suppose that
$Z_p^f$ is a fuzzy number defined by \eqref{zpf} and satisfies the conditions in Proposition \ref{2}, i.e. $p>0$ and $f: [0,1]   \to  [0,1]$ is a differentiable and strictly decreasing function with $f(0)=1$, $f(1)=0$,
and $\lim_{\alpha \to 1-} f'(\al)   =  -\infty$.
Then
$u$ and $Z_p^f$
satisfy all conditions in
Theorem \ref{rap}.
So,
 by Theorem \ref{rap},
 we
 know
 that
 $Z_p^f$
presented in Proposition \ref{2}
are smoothers
for
fuzzy numbers in $  \mathcal{F}^0_\mathrm{NC}( \mathbb{R} )$.
Let
$f=\sqrt{1-t}$, then $Z_p^f$ is just $w_p$ defined in \eqref{wpc}.
So
$w_p$, $p>0$, are also smoothers for
fuzzy numbers in $ \mathcal{F}^0_\mathrm{NC}( \mathbb{R} )$.

}

\end{re}

Theorem \ref{rap} shows  that how to choose smoothers for fuzzy numbers in
$\mathcal{F}_\mathrm{N}(\mathbb{R})  \cap  \mathcal{F}_\mathrm{C}(\mathbb{R})$.
Now we discuss how to use this method to construct a sequence
of smooth fuzzy numbers to approximate the
original fuzzy number.

Suppose that $u_{l,r}$ is a fuzzy number in $\mathcal{F}_\mathrm{N}(\mathbb{R})  \cap  \mathcal{F}_\mathrm{C}(\mathbb{R})$
with
$u(u^-(0)) = l$ and $u(u^+(0)) = r$.
Define corresponding
$v_{l,r, p}$, $p>0$,    as follows.
\[
 v_{l,r, p}(t)=\left\{
\begin{array}{ll}
1-{\left(     t/p      \right)}^2,    &\ t\in [-p\sqrt{1-l},   p\sqrt{1-r}],
\\
0,    & \ t\notin [-p\sqrt{1-l}, \     p\sqrt{1-r}].
\end{array}
\right.
\]
It is easy to check that
$v_{l,r,p}(v_{l,r,p}^-(0))=l$,
$v_{l,r,p}(v_{l,r,p}^+(0))=r$,
and
 then we can see that,
for each $p>0$,
$u_{l,r}$ and $v_{l,r,p}$
satisfy all the conditions in Theorem \ref{rap}.
So $v_{l,r,p}$, $p>0$, are smoothers for $u_{l,r}$.
In
 the
 following, for simplicity, we use $u$ and $v_p$ to denote $u_{l,r}$
and
$v_{l,r,p}$, respectively.
It is easy to observe that if $l=r=0$, then $v_p$ is just $w_p$ defined in \eqref{wpc}.

Notice that
$d_\infty (  u\nabla v_p,\ u      ) \leq    \max \{   p\sqrt{1-l} ,\   p\sqrt{1-t}  \}    \leq     p$
and
    $u\nabla v_p   \in   \mathcal{F}_\mathrm{D}(\mathbb{R})$,
thus
we have the following
conclusion.
\begin{tm}
\label{nct}
Given $u$ in $\mathcal{F}_\mathrm{N}(\mathbb{R})  \cap  \mathcal{F}_\mathrm{C}(\mathbb{R})$,
then
$$d_\infty (   u\nabla v_{1/n},\ u    )     \leq    1/n.$$
So
the smooth fuzzy numbers sequence $\{     u\nabla v_{1/n}: \ n\in \mathbb{N}     \}  $
approximates
 $u$
according to the supremum metric $d_\infty$.
\end{tm}
The
above theorem shows that, for each fuzzy number $u \in    \mathcal{F}_\mathrm{N}(\mathbb{R})  \cap  \mathcal{F}_\mathrm{C}(\mathbb{R})$,
we can find
a sequence
of
smooth fuzzy numbers $ \{ u\nabla v_{1/n} \} $
which approximates $u$
in $d_\infty$ metric.
By
 Theorem \ref{rap}, we can construct other types of smoothers for $u_{l,r}$.
For example, define fuzzy numbers $\xi_{l,r, p}^{f,g}$ as follows:
\[
 \xi_{l,r, p}^{f,g}(t)=\left\{
\begin{array}{ll}
f(  \frac{t-pa}{pb-pa}  ),    &\ t\in [pa, pb],
\\
1,    & \ t\in [pb,  pc],
\\
g(  \frac{t-pc}{pd-pc}   ),   &  \  t\in [pc, pd],
\\
0,    & \ t\notin [pa,  pd],
\end{array}
\right.
\]
where $a<b<c<d$, $p>0$, $f:[0,1] \to [l,1]$ is an increasing and differentiable function which satisfies that $f(0)=l$, $f(1)=1$,
and $f_-'(1)=0$,
and
$g:[0,1] \to [r,1]$ is a decreasing and differentiable function which satisfies that $g(0)=1$, $g(1)=r$,
and $g_+'(0)=0$. Then it can be checked
that
$ \xi_{l,r, p}^{f,g}$, $p>0$, are smoothers of $u_{l,r}$
and that
$u_{l,r} \nabla \xi_{l,r, p}^{f,g}   $
converges to
$ u_{l,r} $
as $p\to 0$.

Next we consider how to use the convolution method to construct a smooth approximation
for an arbitrary
 fuzzy number $u$ in $\mathcal{F}_\mathrm{C}(\mathbb{R})$.
The following lemmas are needed.

\begin{lm}
\label{ase}
Suppose that
$w\in  \mathcal{F}_\mathrm{D}(\mathbb{R})$
and
that
$\al\in [0,1]$,
then
 \\
 (\romannumeral1) \ $w'(w^-(\al))=0$
is equivalent to
      $w'(w_s^-(\al))=0$;
\\
 (\romannumeral2) \ $w'(w^+(\al))=0$
is equivalent to
      $w'(w_s^+(\al))=0$.
\end{lm}

\pof \   (\romannumeral1) \  If
 $w^-(\al)=w_s^-(\al)$,
 then obviously $w'(w^-(\al))=0$ is equivalent to $w'(w_s^-(\al))=0$.
If
 $w^-(\al))<w_s^-(\al)$,
it then follows from
$w\in \mathcal{F}_\mathrm{D}(\mathbb{R})$
that
$w'(w^-(\al))=w'(w_s^-(\al))=0$.

(\romannumeral2) \
The desired conclusion can be proved similarly as (\romannumeral1).
\ep


\begin{lm}
  \label{cndm}
Let $u\in \mathcal{F} (\mathbb{R})$.
Suppose that $w \in \mathcal{F}_\mathrm{D}(\mathbb{R})$ satisfies condition (\romannumeral1)
and
the following condition (\romannumeral4):
\\
(\romannumeral4) \
Suppose that  $ x \in (u^-(0), u^+(0)) $ is a non-differentiable point of $u$,
 then
  \\
(\textbf{\romannumeral4-1}) \
if $  u^-(0)    < x < u^-(1)$, then $w_+'(w^-(\al))=0$, where $\alpha: = u(x)$;
\\
(\textbf{\romannumeral4-2}) \
if $u^+(1) < x < u^+(  0   ) $, then $w_-'(w^+(\beta))=0$, where $\beta: = u(x)$.
\\
Then, for each $z \in (   (u\nabla w)^-(0), (u\nabla w)^+(0) ) $ with $(u\nabla w)(z) < 1$,
\begin{description}
  \item[C1] $ (u\nabla w)' (z)=0$ when $z= (u\nabla w)^- (\alpha)$
and
$u^-(\al) \in (u^-(0), u^+(0))$ is a continuous and non-differentiable point of $u$.

 \item[C2] $ (u\nabla w)' (z)=0$ when $z= (u\nabla w)_s^- (\alpha)$
and
$u_s^-(\al) \in (u^-(0), u^+(0))$ is a continuous and non-differentiable point of $u$.

\item[C3] $ (u\nabla w)' (z)=0$ when $z= (u\nabla w)^+ (\alpha)$
 and
$u^+(\al)\in (u^-(0), u^+(0))$ is a continuous and non-differentiable point of $u$.

 \item[C4] $ (u\nabla w)' (z)=0$ when $z= (u\nabla w)_s^+ (\alpha)$
and
 $u_s^+(\al)\in (u^-(0), u^+(0))$ is a continuous and non-differentiable point of $u$.

\end{description}
\end{lm}

\pof \ We only prove statements \textbf{C1} and \textbf{C2}. Other statements can be proved similarly.
Set
$\al_0:= u(u^-(0)) = w(w^-(0))= (u\nabla w) (  (u\nabla w)^-(0)   ) $.

Suppose that $z=(u\nabla w)^- (\alpha)$
 and
that
 $u^-(\al)\in (u^-(0), u^+(0))$ is a continuous and non-differentiable point of $u$.
Then
$\al_0 < \al <1$, and hence $w^-(\al) \in  (w^-(0),   w^+(0))$.
By Proposition \ref{vec},
$u(u^-(\al)) =  \alpha$, and
therefore,
   by condition (\romannumeral4-1),
$w'(w^-(\al))=0$.
It
 thus follows from Theorem \ref{dzeron} (\romannumeral1), (\romannumeral2)
 that
$(u\nabla w)'   (z) =  (u\nabla w)'   (   (u\nabla w)^-(\al) )   = 0$. So statement \textbf{C1} is true.

Suppose that $z= (u\nabla w)_s^- (\alpha)$
and
$u_s^-(\al) \in (   u^-(0),  u^+(0)   )$ is a continuous and non-differentiable point of $u$.
 Then, by Proposition \ref{vec},
$u(  u_s^-(\al) ) =\al<1$,
and
hence,
  by condition (\romannumeral4-1),
$w_+'(w^-(\al))=0$.
If $\al>\al_0$, then both $w^-(\al)$ and $w_s^-(\al)$ are inner points.
Thus by Lemma \ref{ase}
$w'(  w_s^-(\al )  ) =0$.
So
it
   follows from Theorem \ref{dzeron} (\romannumeral3), (\romannumeral4)
 that
$ (u\nabla w)' (z) = (u\nabla w)'   (   (u\nabla w)_s^-(\al) )   = 0$.
If
 $\al=\al_0$.
Note that
$w^-(\al_0) \leq w_s^-(\al_0)$,
we
know
that
$w_+'(w_s^-(\al_0))=0$,
and
hence
$(u\nabla w)_+' (z) = (u\nabla w)_+'   (   (u\nabla w)_s^-(\al_0) )   = 0$.
Since
$z$ is an inner point of $[u \nabla w]_0$,
we
obtain
that
$(u\nabla w)' (z) = 0$.
Thus
statement \textbf{C2} is proved.
\ep

The following theorem presents a method
 to
find    smoothers
in
 $\mathcal{F}_\mathrm{C}(\mathbb{R})$.
From Theorem \ref{rap}, we have already given
a way
 to pick smoothers for
fuzzy numbers
in
 $\mathcal{F}_\mathrm{C}(\mathbb{R}) \cap \mathcal{F}_\mathrm{N}(\mathbb{R})$.
Now our considerations need include
 fuzzy numbers in
$\mathcal{F}_\mathrm{C}(\mathbb{R}) \backslash \mathcal{F}_\mathrm{N}(\mathbb{R})$
which
have
 one or more non-differentiable points in $[u]_0 \backslash [u]_1$.

\begin{tm}
\label{rndp}
Suppose that $u\in \mathcal{F}_\mathrm{C}(\mathbb{R})$
and
that $w\in      \mathcal{F}_\mathrm{D}(\mathbb{R})$,
then $w$ is a smoother of $u$, i.e.
$u\nabla w   \in \mathcal{F}_\mathrm{D}(\mathbb{R})$,
when
$w$ satisfies the conditions (\romannumeral1), (\romannumeral2),
and
 (\romannumeral4).

\end{tm}

\pof \ To  prove
 that $u\nabla w \in   \mathcal{F}_\mathrm{D}(\mathbb{R}) $,
we
   adopt the same procedure as in the proof of Theorem \ref{rap}.
The
 proof
is divided into
the same situations
as
the proof
 of
Theorem
 \ref{rap}.
Since $u\in \mathcal{F}_\mathrm{C}(\mathbb{R})$,
we know
that each inner point of
$[u]_0$ is a continuous point of $u$.
Hence if $w$
satisfies condition (\romannumeral4),
then
$w$ must satisfy condition (\romannumeral3).
Thus
 only
 case (A) need to be reconsidered.

Case (A)\ $x= (u\nabla w)^-(\al)$ ($(u\nabla w)_s^-(\al)$, $(u\nabla w)^+(\al)$
 $(u\nabla w)_s^+(\al)$)
with
 $(u\nabla w)(x) < 1$,
and
 both $u^-(\al)$ and $w^-(\al)$ ($u_s^-(\al)$ and $w_s^-(\al)$,  $u^+(\al)$ and $w^+(\al)$, $u_s^+(\al)$ and $w_s^+(\al)$) being inner points of $[u]_0$ and $[w]_0$, respectively.

It is easy to see that     $\al<1$.
Note that
 $u\in  \mathcal{F}_\mathrm{C}(\mathbb{R})$,
this implies
that
each inner point of
$[u]_0$ is a continuous point of $u$.
If
 $u^-(\al)$ ($u_s^-(\al)$,  $u^+(\al)$,   $u_s^+(\al)$) is a continuous but non-differentiable point of $u$,
then
by statements  \textbf{C1}--\textbf{C4},
we know that $(u\nabla w)'(x) = 0$.
If
 $u^-(\al)$ ($u_s^-(\al)$,  $u^+(\al)$,   $u_s^+(\al)$) is a differentiable point of $u$,
then
by statements  \textbf{A1}--\textbf{A8},
we can compute
$(u\nabla w)'(x)$.
\ep

%
%
%
%
%
%

The following is a concrete  example by which we illustrate how to
use the results in Theorem \ref{rndp} to    construct
smoothers
for fuzzy numbers in $ \mathcal{F}_\mathrm{C}(\mathbb{R})  $.

\begin{eap}\label{eapnd}{\rm
Suppose
\[ u(t)=\left\{
\begin{array}{ll}

 0.5+t, & \ t\in [-0.5,0),
 \\
0.5+0.5t,  & \ t\in [0,1],
\\
2-t,   & \ t\in (1,2],
\\
0,           & \ t\not\in [-0.5,2].
\end{array}
\right.
\]
\begin{figure}

 \subfloat
 [$u$]
 { \label{eapndfug}
\scalebox{0.26}{   \includegraphics{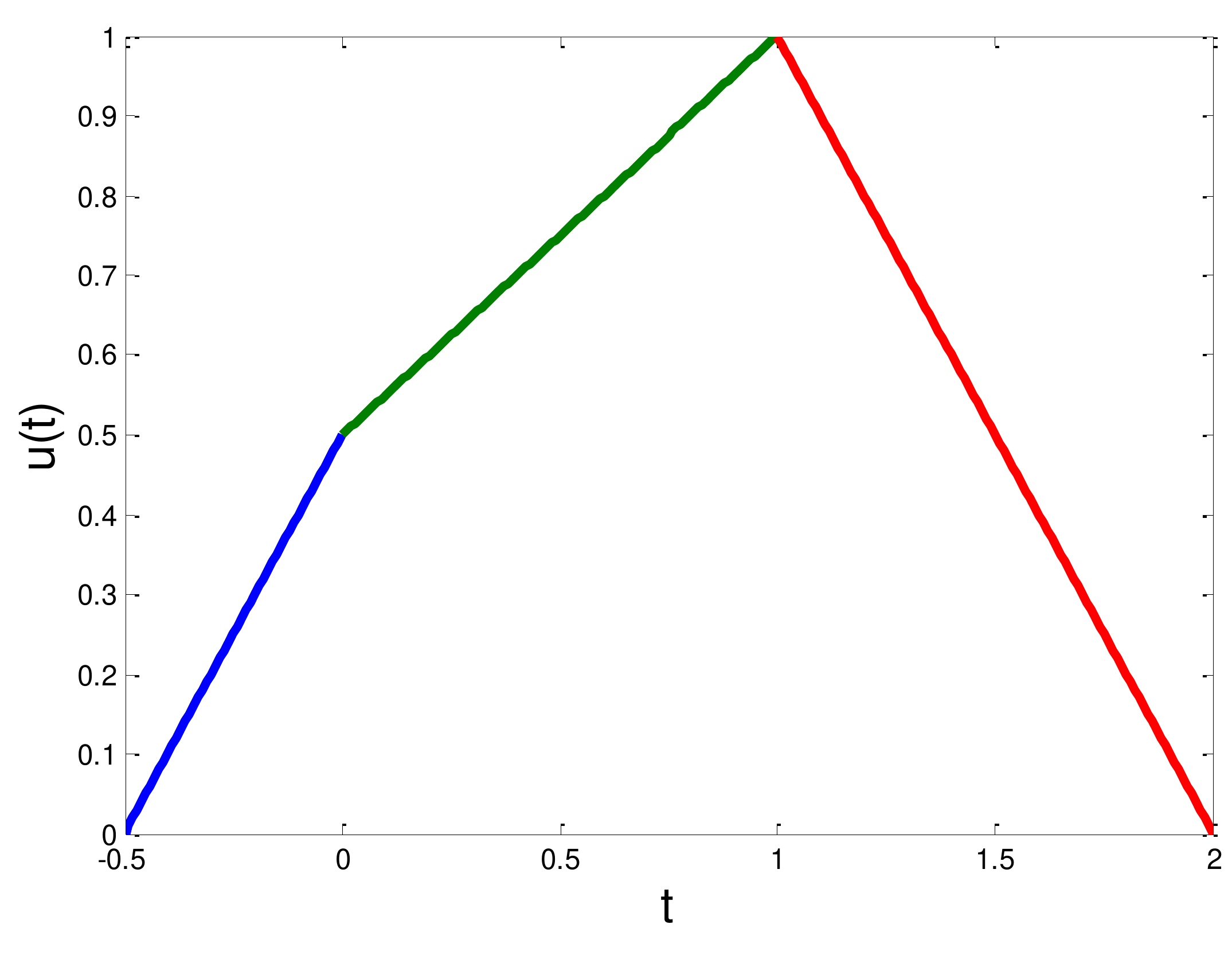}      }
}
 \subfloat
 [$u\nabla w_1$]
 { \label{eapndfuwg}
\scalebox{0.27}{   \includegraphics{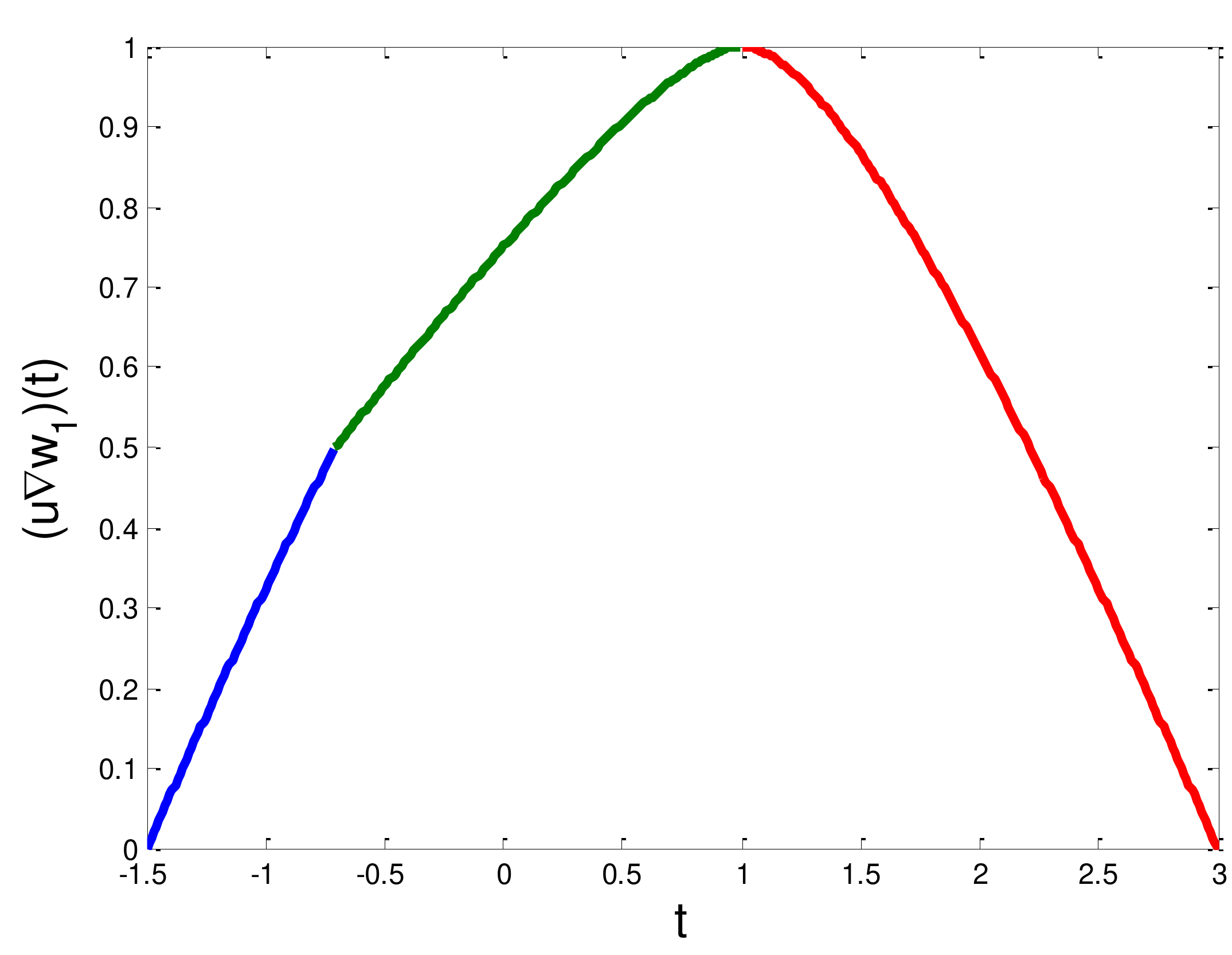}      }
}

\subfloat
 [$p$]
 { \label{eapsfpg}
\scalebox{0.27}{   \includegraphics{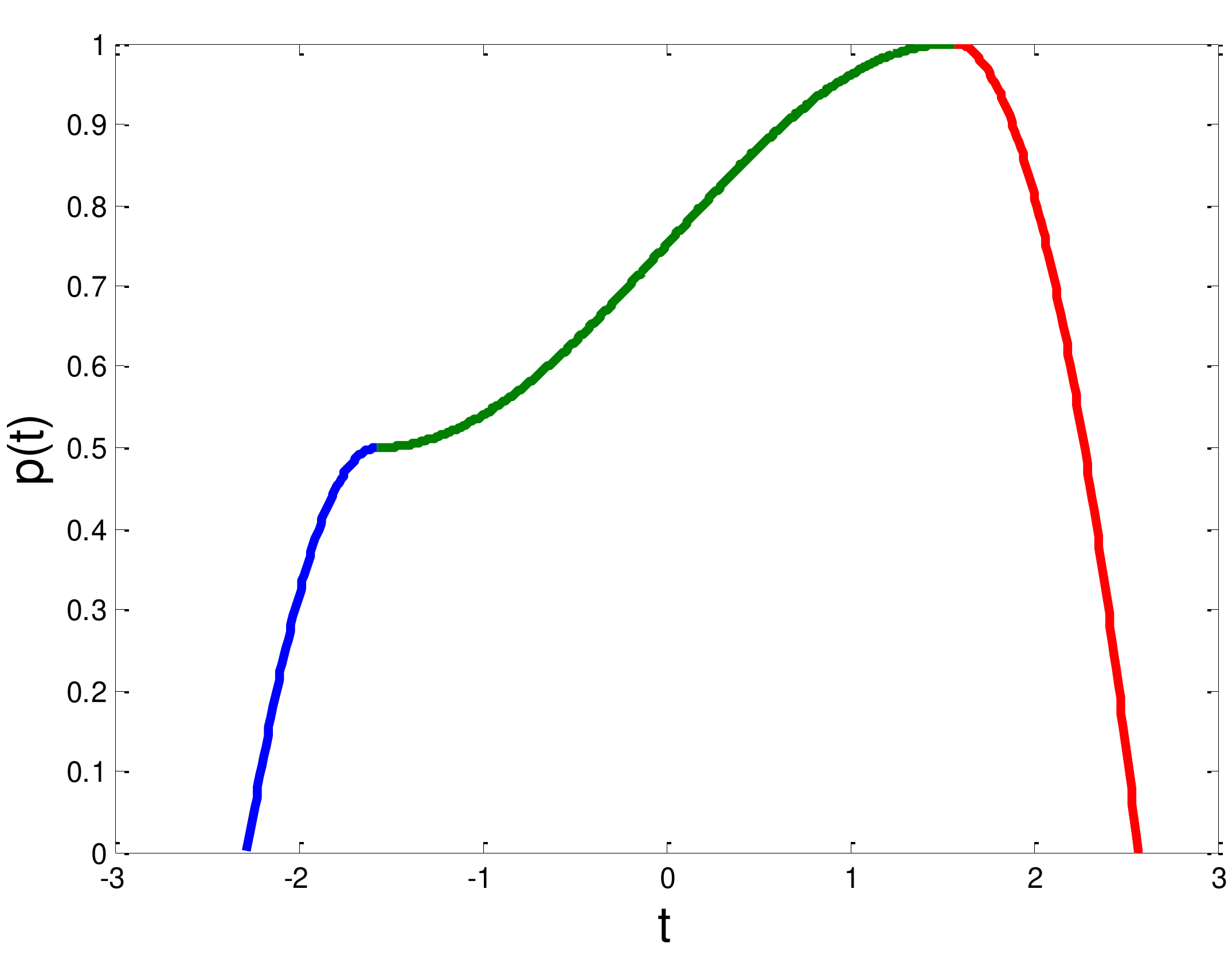}      }
}
\subfloat
 [$u\nabla p$]
 { \label{eapsfupg}
\scalebox{0.27}{   \includegraphics{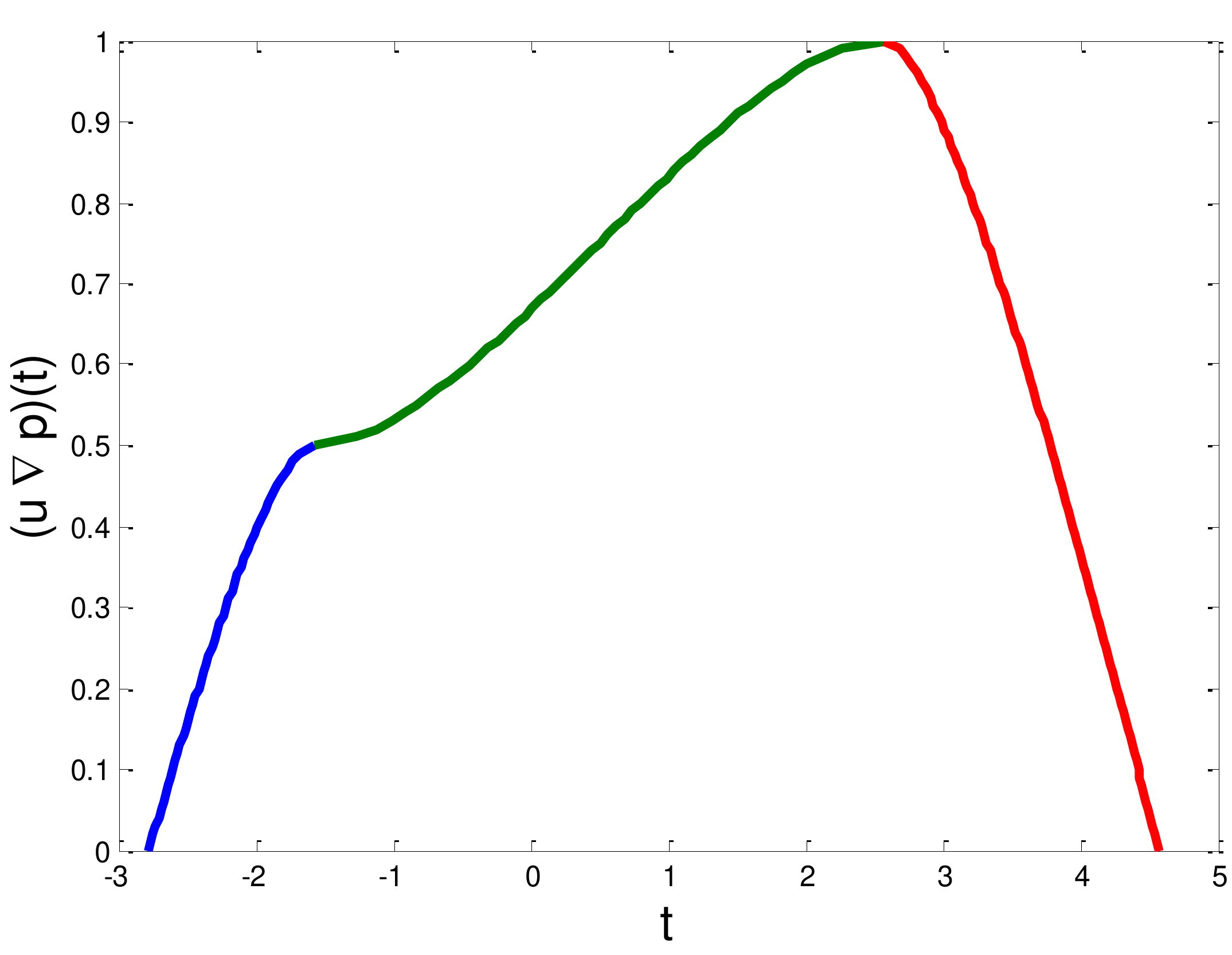}      }
}

 \caption{Fuzzy numbers in Example \ref{eapnd} }
 \label{eapndf}
 \end{figure}
Obviously, $u$ is in $\mathcal{F}_\mathrm{C}(\mathbb{R})$
and has a non-differentiable point $u^-(0.5)= 0$.
See Fig. \ref{eapndfug} for the figure of $u $.

First let's see    whether $w_1$ defined in \eqref{wpc} can work as a smoother for $u$.
Note that ${w_1}'(   { w_1}^-(0.5)   ) \not=0$, thus $u$ and $w_1$
do not satisfy the conditions in Theorem \ref{rndp}.
It can be deduced
that
\[
( u\nabla w_{1} ) (t)
=\left\{
\begin{array}{ll}
\sqrt{\frac{3}{4}-t}+t,  & \    \ t\in [-\frac{3}{2},-\sqrt{0.5}],
\\
\frac{\sqrt{9-8t}+3+4t}{8},   & \    \ t\in (-\sqrt{0.5},1],
\\
\frac{\sqrt{4t-3}-2t+3}{2},    & \     \ t\in(1,3],
\\
0,         & \    \ t\notin [-\frac{3}{2},   3].
\end{array}
\right.
\]
See Fig. \ref{eapndfuwg} for the figure of $u\nabla w_1 $.
Notice that
$
(u\nabla w_{1})'_-(-\sqrt{0.5})=2-\sqrt{2}$,
and
$(u\nabla w_{1})'_+(-\sqrt{0.5})=\frac{3-2\sqrt{2}}{14},
$
 it then follows that $u\nabla w_{1}$ is
not differentiable at $(u\nabla w_1)^-(0.5)   =   -\sqrt{0.5}$, so $u\nabla w_{1}\notin
\mathcal{F}_\mathrm{D}(\mathbb{R})$.
This means that $w_1$ is not a smoother of $u$.

Now we use the conclusions in Theorem \ref{rndp} to construct a smoother for $u$.
Consider
a fuzzy number $p$ defined by
\[p(t)
=\left\{
\begin{array}{ll}
0.5-(t+\pi/2)^2,          &   \   t\in [-\pi/2-\sqrt{0.5},  -\pi/2],
\\
0.25 \sin t+ 0.75,      & \ t \in [-\pi/2, \pi/2],
\\
1-(t-\pi/2)^2,     &   \   t\in [\pi/2, 1+ \pi/2].
\end{array}
\right.
\]
 Then $p'(p^-(0.5))=0$.
See Fig. \ref{eapsfpg} for the figure of $p$.
It's
easy to check that $u$ and $p$
satisfies
all
the conditions of Theorem \ref{rndp}.
Thus by Theorem \ref{rndp}, $p$ is a smoother of $u$.
To validate this assertion,
we
compute
that
\[
[u\nabla p]_\al
=\left\{
\begin{array}{ll}
[-0.5\pi-(0.5-\alpha)^{0.5}+\alpha- 0.5, \  0.5\pi +(1-\alpha)^{0.5} +2-\alpha    ],  &   \   \al\in [0,  0.5],
\\
\mbox{}
 [\arcsin (4\alpha-3)+ 2\alpha-1,    \  0.5\pi +(1-\alpha)^{0.5}+ 2-\alpha ],      & \ \alpha \in [0.5,1].
\end{array}
\right.
\]
Now we can plot $u\nabla p$ in Fig. \ref{eapsfupg},
and obtain
\begin{equation}
\label{ups}
(u\nabla p)(t)
=\left\{
\begin{array}{ll}
 \pi/2 + t + (1 - 4t - 2\pi)^{0.5}/2,  & \ t\in [-\pi/2-\sqrt{0.5}-0.5,  -\pi/2],
\\
f(t),  & \ t\in [-0.5\pi, 1+0.5\pi],
\\
\pi/2 - t + (4t - 2\pi - 3)^{0.5}/2 + 3/2, & \ t\in [ 1+0.5\pi, 3+0.5\pi],
\end{array}
\right.
\end{equation}
  where $f(t)$ is the inverse function of $g(\al)=\arcsin (4\alpha-3)+ 2\alpha-1$ when $\al\in [0.5,1]$.

  It can be checked that $u\nabla p \in \mathcal{F}_\mathrm{D}(\mathbb{R})$.
  In
  fact, from the expression of $u\nabla p$,
we know that
  it only need to show that $u\nabla p$ is differentiable at $-0.5\pi$ and $1+0.5\pi$.
   Use
implicit differentiation
 to take $f'$ as follows:
$$f'(t)=(2+\frac{4}{\sqrt{1-(4\alpha(t)-3)^2}})^{-1}.$$
Then
we have that
 $
  (u\nabla p)_-'(-0.5\pi)=\frac{d}{dt}(\pi/2 + t + (1 - 4t - 2\pi)^{0.5}/2)|_{t=-0.5\pi}=0
  $
and
$
   (u\nabla p)_+'(-0.5\pi)=f'(-0.5\pi)=0,
$
   hence
 \begin{equation}\label{upt} (u\nabla p)'(-0.5\pi)=0.\end{equation}
 Similarly,
it can   be computed
that
\begin{equation}\label{upg}(u\nabla p)'(1+0.5\pi)=0.\end{equation}
It
now
follows from \eqref{ups},\eqref{upt} and \eqref{upg} that
 $u\nabla p$
is differentiable on $(-\pi/2-\sqrt{0.5}-0.5, \ 3+0.5\pi)$, i.e,
$u\nabla p \in \mathcal{F}_\mathrm{D}(\mathbb{R})$.
So $p$ is a smoother of $u$.
From this example we can see how condition (\romannumeral4) takes effect.

 }
\end{eap}

Look at the construction of $p$. To make $u$ and $p$ satisfy the condition (\romannumeral4),
i.e. to assure $p'(p^-(0.5))   =   p'(-\pi/2)=0$,
we use
a polynomial function and a sine function to construct $p$.
Along this line,
given a fuzzy number $u\in \mathcal{F}_\mathrm{C}(\mathbb{R})$
with finite non-differentiable points,
we can   use
 polynomial functions,     sine functions and   cosine functions
to construct a     fuzzy number $\eta \in \mathcal{F}_\mathrm{D}(\mathbb{R})$ such that $u$ and $\eta$
satisfy all conditions in
Theorem \ref{rndp}.
To
construct a smooth fuzzy numbers sequence to approximate $u$,
put $ \eta_p
:=p\cdot \eta$, i.e.
\begin{equation}\label{vps}
  \eta_p(t)
:=(p\cdot \eta) (t)=\eta(t/p)
=\left\{
\begin{array}{ll}
\eta(t/p),  &   \   t\in    p[\eta]_0,
\\
0,      & \ t\notin     p[\eta]_0,
\end{array}
\right.
\end{equation}
where $p>0$ is a real number.
Clearly $\eta=\eta_1$.
It can be checked that
\begin{equation}\label{lrvu}
  \eta(\eta^-(0)) =  \eta_p( {\eta_p}^-(0)  ),\     \eta(\eta^+(0)) =  \eta_p( {\eta_p}^+(0)  ),
\end{equation}
and
\begin{gather*}
  {\eta_p}_-'(t) =  \eta_-'(t/p)/p,
\\
{\eta_p}_+'(t) =  \eta_+'(t/p)/p
\end{gather*}
for all $t\in \mathbb{R}$.
Thus we know
$\eta_p \in \mathcal{F}_\mathrm{D}(\mathbb{R})$,
\begin{gather}
 {\eta_p}_-'( {\eta_p}^-(1) )=\eta'_-(\eta^-(1))/p,   \label{l1dv}
\\
{\eta_p}_+'( {\eta_p}^+(1) )=\eta'_+(\eta^+(1))/p,     \label{r1dv}
\end{gather}
and
\begin{gather}
   \eta_p'( {\eta_p}^-(\al) )=\eta_p'(p\cdot \eta^-(\al)) = \eta'(\eta^-(\al))/p,    \label{ladv}
\\
 \eta_p'( {\eta_p}^+(\al) )=\eta_p'(p\cdot \eta^+(\al)) = \eta'(\eta^+(\al))/p       \label{ladvbe}
 \end{gather}
for each $\al\in (0,1)$.

From eqs. \eqref{lrvu}--\eqref{ladvbe},
we know that
$u$ and $\eta$ satisfy all the conditions in Theorem \ref{rndp}
is equivalent to
that $u$ and $\eta_p$, $p>0$, satisfy all the conditions in Theorem \ref{rndp}.
So
we
have the following conclusion.

\begin{tm}
  Given $u$ in $ \mathcal{F}_\mathrm{C}(\mathbb{R})$.
If the number of   non-differentiable points of   $u$ is finite,
then
there exists smoothers for $u$.
Moreover, if $\eta$ is a smoother of $u$,
then
$p\cdot \eta$, $p>0$, are also smoothers of $u$,
and
 $$d_\infty ( u,   u\nabla (p\cdot \eta)  ) \leq  p\cdot \max \{|\eta^+(0)|, |\eta^-(0)| \}.$$
So
  $  \{ u\nabla (p\cdot \eta)\}$
 converges to $u$
according to the supremum metric $d_\infty$
as $p\to 0$.
\end{tm}

Finally, we consider how to smooth an arbitrarily given fuzzy number $u$, which may have non-continuous points
in $(u^-(0),  u^+(0))$.


\begin{lm}  \label{s3e}
Suppose that $u\in \mathcal{F} (\mathbb{R})$
and
that
 $w \in \mathcal{F}_\mathrm{D}(\mathbb{R})$ satisfies condition (\romannumeral1) and condition (\romannumeral5)
listed below:
\begin{description}

  \item[(\romannumeral5-1)]    if $u^-(0)< x \leq u^-(1)$ and $u$ is not left-continuous at $x$, then $w_+'(w^-(\beta))=0$, where $\beta=\lim_{y\to x-}   u(y)$;

  \item[(\romannumeral5-2)]    if $u^+(1) \leq x < u^+(0)$ and $u$ is not right-continuous at $x$, then $w_-'(w^+(\gamma))=0$, where $\gamma=\lim_{y\to x+}   u(y)$.

\end{description}
  Then, given $x \in (  (u\nabla w)^-(0), (u\nabla w)^+(0)   ) $ with $(u\nabla w)(x)<1$,
the following statements hold.
\begin{description}
\item[B5]   $ (u\nabla w)' (x)=0$ and $x=(u\nabla w)_s^-(0)$
 when $x= (u\nabla w)_s^- (\alpha)$, $w_s^-(\al)= w^-(0) $ and $u^-(0)<u_s^-(\al) = u^+(0) $.

\item[B6]   $ (u\nabla w)' (x)=0$ and $x=(u\nabla w)_s^+(0)$ when $x= (u\nabla w)_s^+ (\alpha)$, $w_s^+(\al)= w^+(0) $ and $u^-(0)=u_s^+(\al) < u^+(0) $.
\end{description}

\end{lm}

\pof \   We only prove statement \textbf{B5}. Statement \textbf{B6} can be proved similarly.
 Set
 $w(w^-(0))=u(u^-(0))=(u\nabla w) (  (u\nabla w)^-(0)    ):=\alpha_0$.
From
 $w_s^-(\al)= w^-(0) $
and
 $u^-(0)<u_s^-(\al_0) = u^+(0) $, we know $\al=\al_0<1$ and
$u^-(0)=u^-(\al_0)<u_s^-(\al_0) =  u^-(1) = u^+(0) $.
Hence
$u(  u_s^-(\al_0)     ) =1$ and $\al_0=\lim_{y \to u_s^-(\al_0)- }  u(y) <1$.
Therefore
by
condition (\romannumeral5-1), we know
$  w'_+(  w_s^-(\al_0)   ) = w'_+(  w^-(0)  ) =   w'_+(  w^-(\al_0)   ) = 0$,
and so, from Theorem \ref{dzeron}(\romannumeral4),
$
 (u\nabla w)_+' (x) =(u\nabla w)_+' (   (u\nabla w)_s^-(\al_0) ) = 0$.
   Notice that $(u\nabla w)_s^-(\al_0)$ is an inner point
of
$[u \nabla w]_0$,
we thus obtain
$$ (u\nabla w)' (x) =(u\nabla w)' (   (u\nabla w)_s^-(\al_0) ) = 0. \quad \Box$$


\begin{lm}
  Suppose that $u\in \mathcal{F} (\mathbb{R})$
and
that
 $w \in \mathcal{F}_\mathrm{D}(\mathbb{R})$ satisfies conditions (\romannumeral1), (\romannumeral4)
and
  (\romannumeral5).
Then, given $x \in (  (u\nabla w)^-(0), (u\nabla w)^+(0)   ) $ with $(u\nabla w)(x)<1$,
the following statements hold.
\begin{description}
  \item[D1]  $ (u\nabla w)' (x)=0$
 when $x= (u\nabla w)^- (\alpha)$, $w^-(0)<w^-(\al) < w^+(0) $, $u^-(0)<u^-(\al) < u^+(0) $,
$u$ is not continuous at $u^-(\al)$
and
$u(  u^-(\al)  ) = \al$.

\item[D2]  $ (u\nabla w)' (  x   )= w'(w^-(\al))$
 when $x= (u\nabla w)^- (\alpha)$, $w^-(0)<w^-(\al) < w^+(0) $, $u^-(0)<u^-(\al) \leq u^+(0) $,
$u$ is not continuous at $u^-(\al)$,
$
u(   u^-(\al)   ) >  \al$ and $\lim_{ y\to  u^-(\al)- } u(y) < \alpha$.

 \item[D3]  $ (u\nabla w)' (x)=0$
 when $x= (u\nabla w)^- (\alpha)$, $w^-(0)<w^-(\al) < w^+(0) $, $u^-(0)<u^-(\al) \leq u^+(0) $,
$u$ is not continuous at $u^-(\al)$,
$
u(u^-(\al)) >\al$, and $\lim_{   y\to u^-(\al)-   }   u(y) = \alpha$.

 \item[D4] $ (u\nabla w)' (x)=0$
 when $x= (u\nabla w)_s^- (\alpha)$, $w^-(0)<w_s^-(\al) < w^+(0) $, $u^-(0)<u_s^-(\al) < u^+(0) $,
$u$ is not continuous at $u_s^-(\al)$
 and
$u(u_s^-(\al))=\al$.

 \item[D5] $ (u\nabla w)' (x)=    w'( w_s^-(\al)  $
 when $x= (u\nabla w)_s^- (\alpha)$, $w^-(0)< w_s^-(\al) < w^+(0) $, $u^-(0)<u_s^-(\al) \leq u^+(0) $,
$u$ is not continuous at $u_s^-(\al)$,
$u(   u_s^-(\al)   ) >  \al$ and $\lim_{ y\to  u_s^-(\al)- } u(y) < \alpha$.

  \item[D6] $ (u\nabla w)' (x)=0$
 when $x= (u\nabla w)_s^- (\alpha)$, $w^-(0)< w_s^-(\al) < w^+(0) $, $u^-(0)<u_s^-(\al) \leq u^+(0) $,
$u$ is not continuous at $u_s^-(\al)$,
$
u(u_s^-(\al)) > \al$     and        $\lim_{   y\to u_s^-(\al)-   }   u(y) = \alpha$.

 \item[D7]  $ (u\nabla w)' (x)=0$
 when $x= (u\nabla w)^+ (\alpha)$, $w^-(0)<w^+(\al) < w^+(0) $, $u^-(0)<u^+(\al) < u^+(0) $,
$u$ is not continuous at $u^+(\al)$
and
$u(  u^+(\al)  ) = \al$.

\item[D8]  $ (u\nabla w)' (  x   )= w'(w^+(\al))$
 when $x= (u\nabla w)^+ (\alpha)$, $w^-(0)<w^+(\al) < w^+(0) $, $u^-(0) \leq   u^+(\al) < u^+(0) $,
$u$ is not continuous at $u^+(\al)$,
$
u(   u^+(\al)   ) >  \al$ and $\lim_{ y\to  u^+(\al)+ } u(y) < \alpha$.

\item[D9]  $ (u\nabla w)' (x)=0$
 when $x= (u\nabla w)^+ (\alpha)$, $w^-(0)<w^+(\al) < w^+(0) $, $u^-(0) \leq    u^+(\al) < u^+(0) $,
$u$ is not continuous at $u^+(\al)$,
$
u(u^+(\al)) >\al$, and $\lim_{   y\to u^+(\al)+   }   u(y) = \alpha$.

 \item[D10] $ (u\nabla w)' (x)=0$
 when $x= (u\nabla w)_s^+ (\alpha)$, $w^-(0)<w_s^+(\al) < w^+(0) $, $u^-(0)<u_s^+(\al) < u^+(0) $,
$u$ is not continuous at $u_s^+(\al)$
 and
$u(u_s^+(\al))=\al$.

 \item[D11] $ (u\nabla w)' (x)=   w'( w_s^+(\al) $
 when $x= (u\nabla w)_s^+ (\alpha)$, $w^-(0)< w_s^+(\al)  < w^+(0) $, $u^-(0) \leq u_s^+(\al) < u^+(0) $,
$u$ is not continuous at $u_s^+(\al)$,
$u(   u_s^+(\al)   ) >  \al$ and $\lim_{ y\to  u_s^+(\al)+ } u(y) < \alpha$.

\item[D12] $ (u\nabla w)' (x)=0$
 when $x= (u\nabla w)_s^+ (\alpha)$, $w^-(0)< w_s^+(\al) < w^+(0) $, $u^-(0) \leq  u_s^+(\al) < u^+(0) $,
$u$ is not continuous at $u_s^+(\al)$,
$
u(u_s^+(\al)) > \al$     and        $\lim_{   y\to u_s^+(\al)+   }   u(y) = \alpha$.

\end{description}

\end{lm}


\pof \  We only prove statements \textbf{D1}-- \textbf{D6}.
The
remainder
statements can be proved similarly.
Clearly, $\al < 1$.

To show statement \textbf{D1}, notice that $u^-(\al)$ is also a non-differentiable point,
so
by condition (\romannumeral4-1), we know that $w'(w^-(\al))=0$, and thus
$(u\nabla w)' (    (u\nabla w)^-(\al)    )=0$.
This is
statement \textbf{D1}.

To prove statement \textbf{D2},
note that
$w( w^-(\al) )=\alpha$.
Thus
 by Theorem \ref{dbzeron} (\romannumeral3), (\romannumeral4),
 we know that
$ (u\nabla w)' (   (u\nabla w)^-(\al)    )= w'(w^-(\al))$.
Hence
statement \textbf{D2}
holds.

To show statement \textbf{D3},
observe that $u^-(\al)$ is a non-continuous point,
hence,
by condition (\romannumeral5-1),
$w'(  w^-(\al)  )=  w'_+(  w^-(\al)  )= 0$,
and thus
by Theorem \ref{dzeron} (\romannumeral1), (\romannumeral2),
$(u\nabla w)' (    (u\nabla w)^-(\al)    )=0$.
So
statement \textbf{D3}
is true.

To demonstrate statement \textbf{D4},
 observe that $u_s^-(\al)$ is also a non-differentiable point,
then
by condition (\romannumeral4-1), we know that $w'(w^-(\al))=0$, and
hence, by Lemma \ref{ase}, $w'(w_s^-(\al))=0$.
Thus
$(u\nabla w)' (    (u\nabla w)_s^-(\al)    )=0$.
So statement \textbf{D4}
is
proved.

To show statement \textbf{D5},
from $u(   u_s^-(\al)   ) >  \al$ and $\lim_{ y\to  u_s^-(\al)- } u(y) < \alpha$,
we know
that
$ u_s^-(\al) = u^-(\al)$.
If
$ w_s^-(\al) = w^-(\al) $,
                      then
statement \textbf{D5} is just statement \textbf{D2}.
Hence
$(u\nabla w)'(x) = w'(  w^-(\al)  ) =  w'(  w_s^-(\al)  ) $.
If
 $ w^-(\al)  < w_s^-(\al) $, then $w'(  w_s^-(\al) )   )=0$, and thus
from
Theorem \ref{dzeron} (\romannumeral3), (\romannumeral4),
we know that
$(u\nabla w)' (   (u\nabla w)_s^-(\al)    ) = 0=  w'(  w_s^-(\al)    )   $.
So
statement \textbf{D5}
is
proved.

To prove statement \textbf{D6}.
If $ w_s^-(\al)  >    w^-(\al) $,
then, from $w \in  \mathcal{F}_\mathrm{D}(\mathbb{R})$,
we know
that
$w'(  w_s^-(\al)  )=0$.
If
 $ w_s^-(\al)  =    w^-(\al) $,
note
 that
 $u_s^-(\al)$ is a non-left-continuous point,
hence,
by condition (\romannumeral5-1),
$w_+'(  w^-(\al)  )=0$,
and
therefore
$w'(  w_s^-(\al)  )  =  w_+'(  w_s^-(\al)  )   =0$.
Thus
by Theorem \ref{dzeron} (\romannumeral3), (\romannumeral4),
$(u\nabla w)' (    (u\nabla w)_s^-(\al)    )=0$.
So
statement \textbf{D6} is proved.
\ep

\begin{re} \label{cmn}
  It can be checked
that
for each $u, w \in \mathcal{F}(\mathbb{R})$,
if $w$ satisfies conditions (\romannumeral4) and (\romannumeral5),
then
$w$ also
satisfies condition (\romannumeral3).

\end{re}


\begin{tm}
\label{rncp}
Suppose that $u\in \mathcal{F}(\mathbb{R})$
and
 $w\in      \mathcal{F}_\mathrm{D}(\mathbb{R})$,
then $w$ is a smoother of $u$, i.e.
$u\nabla w   \in \mathcal{F}_\mathrm{D}(\mathbb{R})$,
when
$w$ satisfies conditions (\romannumeral1), (\romannumeral2), (\romannumeral4)
 and
(\romannumeral5).

\end{tm}

\pof \ To  prove
 that $u\nabla w \in   \mathcal{F}_\mathrm{D}(\mathbb{R}) $,
we
   adopt the same procedure as in the proof of Theorem \ref{rap}.
The
 proof
is divided into
the same situations
as
the proof
 of
Theorem
 \ref{rap}.
We can see that only cases (A) and (B) need to be reconsidered.
It is clear that
$\al<1$ in these two cases.

Case (A)\ $x= (u\nabla w)^-(\al)$ ($(u\nabla w)_s^-(\al)$, $(u\nabla w)^+(\al)$,  $(u\nabla w)_s^+(\al)$), $(u\nabla w)(x) < 1$,
with
 $u^-(\al)$ and $w^-(\al)$ ($u_s^-(\al)$ and $w_s^-(\al)$,  $u^+(\al)$ and $w^+(\al)$, $u_s^+(\al)$ and $w_s^+(\al)$) being inner points of $[u]_0$ and $[w]_0$, respectively.

If
 $u^-(\al)$ ($u_s^-(\al)$,  $u^+(\al)$,   $u_s^+(\al)$) is a non-continuous point of $u$,
then
by statements  \textbf{D1}--\textbf{D12},
we
can compute
 $(u\nabla w)'(x) $.
From
 the proof of case (A) in Theorem \ref{rndp}, we know $ u\nabla w $
is
 differentiable at $x $
when
$u^-(\al)$ ($u_s^-(\al)$,  $u^+(\al)$,   $u_s^+(\al)$) is a continuous point of $u$.

Case (B)\ $x= (u\nabla w)^-(\al)$ ($(u\nabla w)_s^-(\al)$, $(u\nabla w)^+(\al)$,  $(u\nabla w)_s^+(\al)$), $(u\nabla w)(x) < 1$,
and $x$ is not in
Case (A).

From
Remark \ref{cmn}, $w$ satisfies condition (\romannumeral3).
So,
for $x$ in subcases \textbf{B\romannumeral1}-- \textbf{B\romannumeral4},
we
can prove
that
$u \nabla   w$
is
differentiable
at
$x$
by using statements \textbf{B1}--\textbf{B4}.

Note that $u$ may not in
$\mathcal{F}_\mathrm{C}(\mathbb{R})$,
 we also need to consider the following subcases.
\begin{description}
\item[B\romannumeral5] $x= (u\nabla w)_s^- (\alpha)$,
$w_s^-(\al)= w^-(0) $ and $u^-(0) < u_s^-(\al) = u^+(0) $.

\item[B\romannumeral6]  $x= (u\nabla w)_s^+ (\alpha)$,
$w_s^+(\al)= w^+(0) $ and $u^-(0)   =    u_s^+(\al) < u^+(0) $.

\item[B\romannumeral7] $x= (u\nabla w)_s^- (\alpha)$,
$ w^-(0)< w_s^-(\al) <  w^+(0)$ and $u^-(0) < u_s^-(\al) = u^+(0) $.

\item[B\romannumeral8]  $x= (u\nabla w)_s^+ (\alpha)$,
$ w^-(0) <   w_s^+(\al) < w^+(0) $ and $u^-(0)   =    u_s^+(\al) < u^+(0) $.

 \item[B\romannumeral9]    $x = (u\nabla w)^-(\al)$,   $w^-(0)< w^-(\al) <w^-(1)$
and
 $u^-(0)< u^-(\al) = u^+(0)$.

  \item[B\romannumeral10]   $x = (u\nabla w)^+(\al)$,
$w^+(1) < w^+(\al) <w^+(0)$
and
$u^-(0)= u^+(\al) < u^+(0)$.
\end{description}
By
statements \textbf{B5} and \textbf{B6} in Lemma \ref{s3e},
we
can
deduce that
$(u\nabla w)'(x)=0$
when
$x$ is in subcases \textbf{B\romannumeral5} and  \textbf{B\romannumeral6}.

From
$u^-(0) < u_s^-(\al) = u^+(0) $ , we can deduce that $u( u_s^-(\al)    ) =1 >\al$ and $\lim_{   y\to u_s^-(\al)-   }   u(y) \leq \alpha$.
So,
by statements \textbf{D5} and \textbf{D6},
we can
compute
$(u\nabla w)'(x)$ when $x$ is in
subcases \textbf{B\romannumeral7}.
Similarly,
from
statements \textbf{D11} and \textbf{D12},
we can
compute
$(u\nabla w)'(x)$ when $x$ is in
subcases \textbf{B\romannumeral8}.

Suppose that $u^-(0)< u^-(\al) = u^+(0)$,
then
$u^-(0)< u^-(\al) = u^-(1)= u^+(1)=u^+(0)$ and
$
u( u^-(\al)    ) =1 > \al$.
So, by using statements \textbf{D2}, \textbf{D3},
we can
compute
$(u\nabla w)'(x)$ when $x$ is in
subcases \textbf{B\romannumeral9}.
Similarly,
from
statements
\textbf{D8} and \textbf{D9},
we can
compute
$(u\nabla w)'(x)$ when $x$ is in
  \textbf{B\romannumeral10}.
\ep

\begin{re}
Suppose that $u$ is a fuzzy number
and that
$x\in  (u^-(0), u^+(0)) $     is a non-continuous point of $u$.
Set $u(x):=\alpha<1$.
Then it can be checked that
 \begin{center}  if $x< u^-(1)$, then $x=u^-(\alpha)$; if $x> u^+(1)$, then $x=u^+(\alpha)$. \end{center}
Actually,
   if $x< u^-(1)$, then $x \geq u^-(\al)$.
   Assume     that $x >  u^-(\al)$,
    then we know that
    $u(y)\equiv \alpha$ for all $y\in [ u^-(\al), x]$,
    thus $u$ is left-continuous at $x$. Since $u$ is right-continuous on $(u^-(0),  u^-(1))$,
we know that
 $u$ is continuous at $x$,
 which is a contradiction. So $x=u^-(\al)$.
 In
 a similar way, we can show that
 if $x> u^+(1)$, then $x=u^+(\alpha)$.
\end{re}

The following example
 shows that how to
use the results in Theorem \ref{rncp} to    construct
smoother
for a fuzzy number $u   \notin \mathcal{F}_\mathrm{C}(\mathbb{R})  $.

\begin{eap}\label{eapnc}
{\rm

Suppose that $u$ is a fuzzy number defined by
$$
u(t)
=\left\{
\begin{array}{ll}
 t-1,  & \ t\in [1,2],
\\
-t+3,  & \ t\in [2,   2.5],
\\
-t+2.8, & \ t\in (2.5,   2.8],
\\
0,   & \  t\not\in [1, 2.8].
\end{array}
\right.
$$
The
 figure of $u$ is in Fig \ref{eapncfug}.
We can see that $u$ is discontinuous at $u^+(0.5) = 2.5$.
So $u\notin \mathcal{F}_\mathrm{C}(\mathbb{R}) $.
\begin{figure}

 \subfloat
 [$u$]
 { \label{eapncfug}
\scalebox{0.247}{   \includegraphics{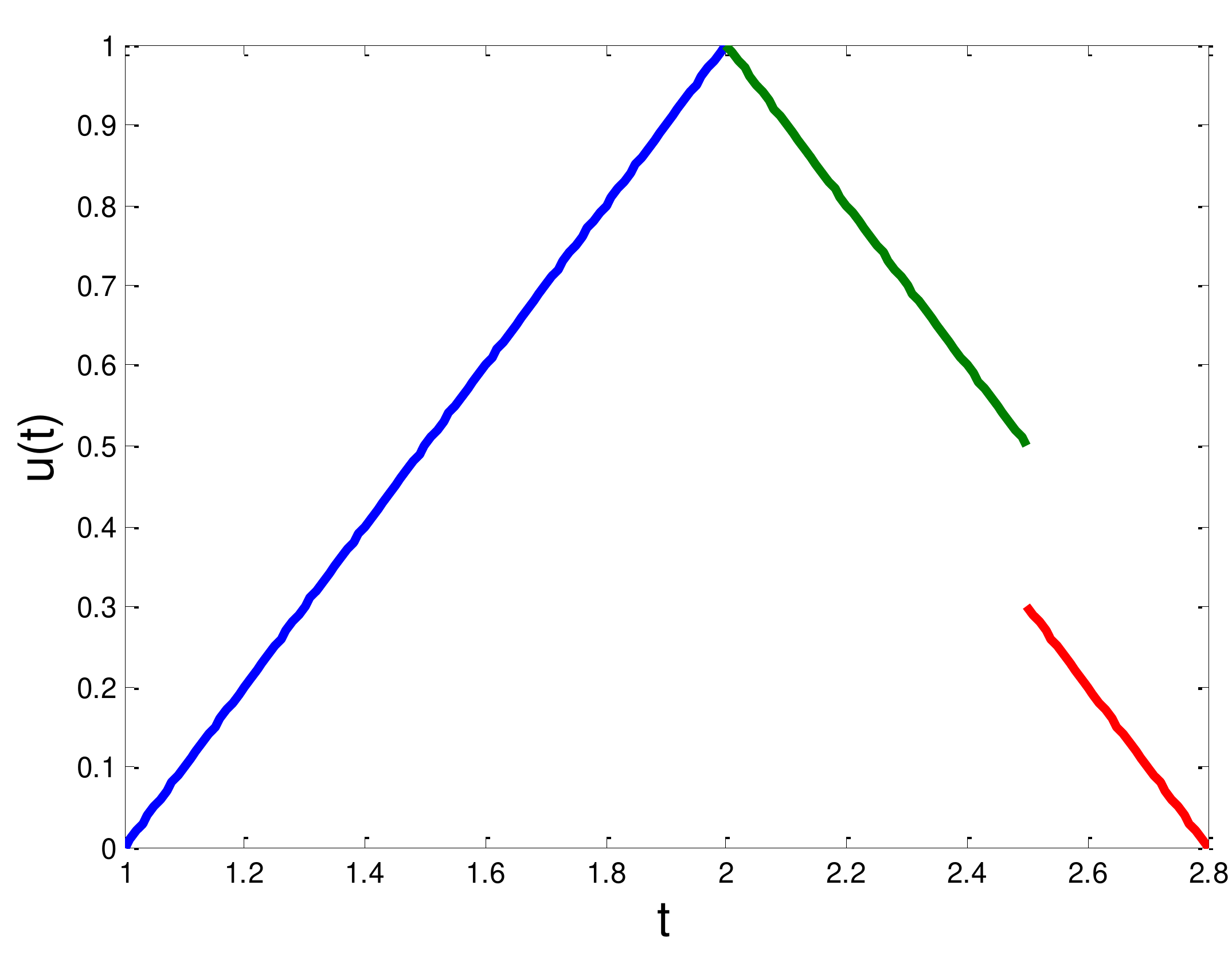}      }
}
\subfloat
 [$z$]
 { \label{eapncfzg}
\scalebox{0.247}{   \includegraphics{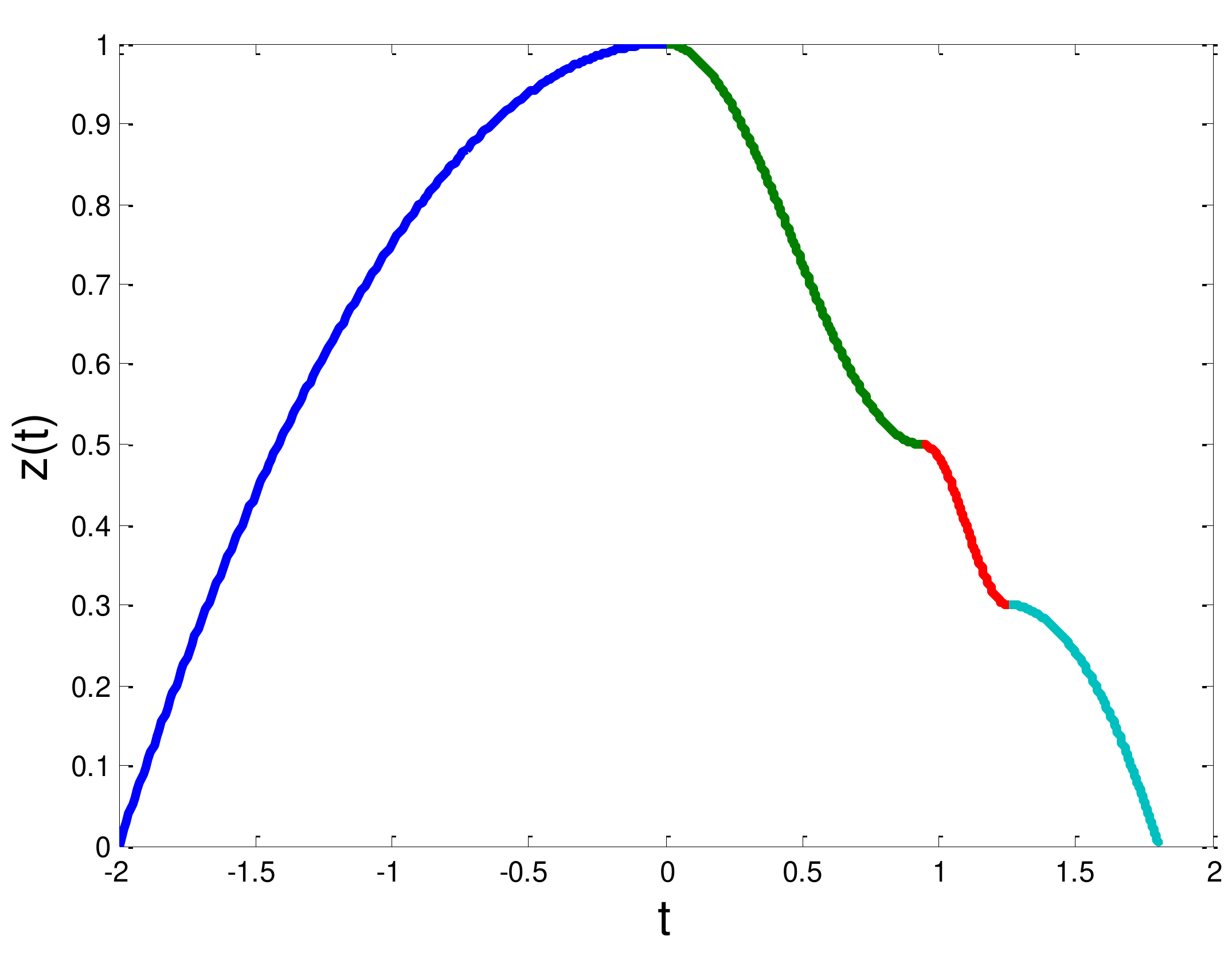}      }
}
\subfloat
 [$u\nabla z$]
 { \label{eapncfuzg}
\scalebox{0.247}{   \includegraphics{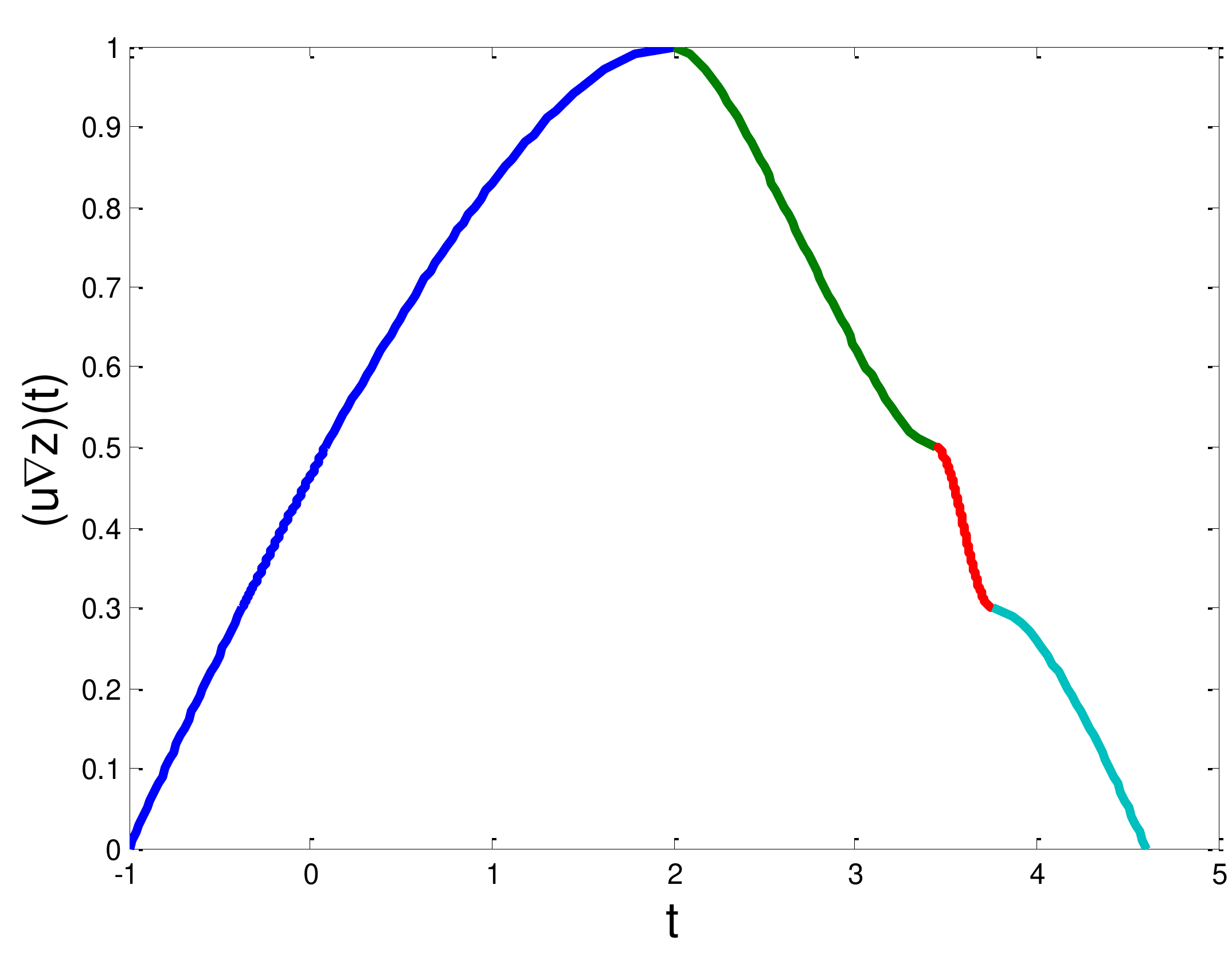}      }
}

\label{eapncf}
 \caption{Fuzzy numbers in Example \ref{eapnc} }
 \label{eapsf}
 \end{figure}
Now
 we use Theorem \ref{rncp} to construct a smoother $z$   for $u$.
Observe that the only non-differentiable
 point of $u$ in $(    u^-(0),   u^+(0)   )$ is $u^+(0.5) = 2.5$,
which is also the only non-continuous point of $u$.
 Since $u(2.5)=0.5$ and $\lim_{y\to 0.5+}   u(y)=0.3$,
 by conditions (\romannumeral4-2) and(\romannumeral5-2), it
must holds
that
$z'(  z^+(0.5)  )=0$
and
 $z'(   z^+(  0.3   )    )=0$.
 Consider
 \[
z(t)=\left\{
\begin{array}{ll}
1-t^2/4,& \ t\in [-2,0],
\\
0.25 \cos (10t/3) + 0.75, & \ t\in [0,  0.3\pi],
\\
0.4+0.1 \cos (   10( t-0.3\pi)  ), & \ t\in [0.3\pi,  0.4\pi],
\\
0.3-(t-0.4\pi)^2,  & \ t\in [0.4\pi, 0.4\pi+\sqrt{0.3}],
\\
0,    & \  t\notin [-2, 0.4\pi+\sqrt{0.3}].
\end{array}
\right.
\]
See Fig. \ref{eapncfzg} for the figure of $z$.
It can be checked that $z\in \mathcal{F}_\mathrm{D}(\mathbb{R})$
and that
\begin{gather*}
  z'(  z^+(0.5)  )=z'(0.3\pi)=0,
\\
  z'(   z^+(  0.3   )    )=z'(0.4\pi)=0.
\end{gather*}
We can see
 that $u$ and $z$ also satisfy all other conditions in Theorem \ref{rncp}.
Thus, by Theorem \ref{rncp}, $z$ is a smoother of $u$. To validated this assertion,
we compute that
\[
(u\nabla z)(t)
=\left\{
\begin{array}{ll}
t + 2(3 - t)^{0.5} - 3,   & \ t\in [-1,2],
\\
 f(t),     & \ t\in [2, 0.3\pi + 2.5],
\\
2/5 - \cos(10t - 25)/10,    &  \   t\in [0.3\pi + 2.5, \  0.4\pi + 2.5 ],
\\
 0.4\pi - t + 0.5(4t - 1.6\pi - 9)^{0.5} + 2.3 ,    &  \   t\in [  0.4\pi + 2.5, \   0.4\pi + 2.8 + 0.3^{0.5} ],
\\
0, &  t\notin [-1,  0.4\pi + 2.8 + 0.3^{0.5}]
\end{array}
\right.
\]
where $f$ is the inverse function of  $0.3\arccos(4\alpha - 3) +3-\al $ when $\alpha\in [0.5,1]$.
See Fig. \ref{eapncfuzg} for the figure of $u\nabla z$.
 It
 can be verified that $u\nabla z \in \mathcal{F}_\mathrm{D}(\mathbb{R})$.
  In
  fact,
  it only need to show that $u\nabla z$ is differentiable at points $2$, $0.3\pi + 2.5$, and    $0.4\pi + 2.5$.
   Use
implicit differentiation
 to take $f'$ as follows:
$$f'(t)=\left(      \frac{-1.2}{\sqrt{1-(4\alpha(t)-3)^2}}   -   1      \right)^{-1}.$$
Then
 $
  (u\nabla z)_-'(2)=\frac{d}{dt}(  1 + 0.5(3-t)^{-0.5}\cdot(-1)\cdot 2   )|_{t=2}=0$,
  and
   $(u\nabla z)_+'(2)=f'(2)=0$.
 Hence
 \begin{equation}\label{uzt} (u\nabla z)'(2)=0.\end{equation}

 Similarly, it can   be computed that
\begin{gather}
(u\nabla z)'(0.3\pi + 2.5)=0, \label{uzpi}
\\
(u\nabla z)'(0.4\pi + 2.5)=0. \label{uz4pi}
\end{gather}
It
now
follows from \eqref{uzt}, \eqref{uzpi}, and \eqref{uz4pi} that
 $u\nabla z$
is differentiable on $(   -1,    0.4\pi + 2.8 + 0.3^{0.5}   )$, i.e,
$u\nabla z \in \mathcal{F}_\mathrm{D}(\mathbb{R})$.
This means that
$z$ is a smoother of $u$.
From this example we can see how conditions (\romannumeral4) and (\romannumeral5)  take effect.

}
\end{eap}

By using a similar procedure as described in Example \ref{eapnc},
we can
 construct a smoother $\zeta$ for $u\in \mathcal{F}(\mathbb{R})$
which has finite number of non-differentiable points.
Note that a non-continuous point is also a non-differentiable point,
this
 means
 that
 the number of non-continuous points
and the number of continuous but non-differentiable points are both finite.
First, it picks out all the  non-continuous points and continuous but non-differentiable points
in $(u^-(0),  u^+(0))$.
Based on this,
 by using conditions (\romannumeral1), (\romannumeral2), (\romannumeral4) and (\romannumeral5),
we then give some requirements on $\zeta$    which ensure $\zeta$
to be a smoother for $u$.
Finally,
we
 use
 polynomial functions,     sine functions and   cosine functions
to construct a concrete fuzzy number $\zeta\in \mathcal{F}_\mathrm{D}(\mathbb{R})$
which meets these requirements.
From
Theorem \ref{rncp},
we know that
$\zeta$ is a smoother of $u$.
Since the number of non-differentiable points of $u$ is finite, this procedure can be completed in finite steps.

Now we discuss how
to construct a smooth fuzzy numbers sequence to approximate a fuzzy number in $\mathcal{F}(\mathbb{R})$.
In fact, it can   proceed   as in
the construction of    smooth fuzzy numbers sequences for   fuzzy numbers
in $\mathcal{F}_\mathrm{C}(\mathbb{R})$.

Suppose that $u, \zeta \in \mathcal{F} (\mathbb{R})$.
From
eqs. \eqref{lrvu}--\eqref{ladvbe} and Theorem \ref{rncp},
we know
 that
$\zeta$   is a smoother of $u$
is equivalent to
    $p\cdot \zeta$, $p>0$, are smoothers of $u$.
So
 we have the following statement which shows that,
by using the convolution method, it can produce a smooth approximation for
an arbitrary fuzzy number with finite non-differentiable points.
\begin{tm} \label{san}
 Suppose that $u$ is a fuzzy number.
If the number of   non-differentiable points of   $u$ is finite,
then there exist smoothers for $u$.
Moreover, if $\zeta$ is a smoother of $u$,
then $p\cdot \zeta$, $p>0$, are also smoothers of $u$,
 and
\begin{equation*}\label{cveu}
  d_\infty ( u,   u\nabla (p\cdot\zeta)  ) \leq  p\cdot \max \{|\zeta^+(0)|, |\zeta^-(0)| \}.
\end{equation*}
So
 $  \{ u\nabla (p\cdot \zeta) \}$
 converges to $u$
according to the supremum metric $d_\infty$ as $p\to 0$.
\end{tm}

\section{Properties of the approximations generated by the convolution method}

In this section we discuss the properties of the approximations generated by the convolution method.

We affirm that
 the convolution method
    can produce smooth approximations which preserve the core,
where an   approximation $\{v_n\}$ of a fuzzy number $u$ preserves the core means
that
Core($v_n$)=Core($u$)
for all $n$.
In
fact, if a smoother $v$ of $u$ satisfies the condition $[v]_1=\{0\}$,
then $[p\cdot v]_1=\{0\}$, and thus the corresponding approximation
$\{u\nabla (\frac{1}{n} \cdot v), n=1,2,\ldots\}$ for the fuzzy number $u$
preserves the core.
It
is easy to select a smoother $v$ which satisfies the condition $[v]_1=\{0\}$.
For example,
 if $w$ is a smoother of a fuzzy number $u$,
then we can     define $v$ by
$$
[v]_\al=[w^-(\al)-w^-(1),  w^+(\al)-w^+(1)]
$$
for all $\al\in [0,1]$.
It can be check that $[v]_1=\{0\}$,
and that,
 from Theorems \ref{rap}, \ref{rndp}, \ref{rncp},
$v$   is also   a smoother for $u$.

We also find that
 the convolution method can generate Lipschitz and smooth approximation, where a
Lipschitz   approximation is a approximation which is constructed by Lipschitz fuzzy numbers.
To
 show this assertion, we need the following lemma.

\begin{lm} \label{glec} Suppose that
$u$
is
a fuzzy number. Then the following statements are equivalent.
\\
(\romannumeral1) \ $u$ is Lipschitz with Lipschitz constant $K$.
\\
(\romannumeral2) \
$|\alpha - \beta|   \leq  K  | u^-(\al)  - u^-(\beta) |$
for all
$\alpha, \beta \in [u(u^-(0)),1]$,
and
$|\gamma - \delta|   \leq  K  | u^+(\gamma)  - u^+(\delta) |$
for all $\gamma$, $\delta$ in $[u( u^+(0)),1]$.
\\
(\romannumeral3) \
$|\alpha - \beta|   \leq  K  | u_s^-(\al)  - u^-(\beta) |$
for all
$\alpha, \beta \in [u(u^-(0)),1]$,
and
$|\gamma - \delta|   \leq  K  | u_s^+(\gamma)  - u^+(\delta) |$
for all $\gamma$, $\delta$ in $[u( u^+(0)),1]$.

\end{lm}

\pof \ (\romannumeral1)$\Rightarrow$(\romannumeral2).
If
 $u$ is Lipschitz with Lipschitz constant $K$,
then
$u\in \mathcal{F}_\mathrm{C}(\mathbb{R})$.
Given
  $\alpha, \beta \in [u(u^-(0)),1]$,
 by Corollary \ref{ncp},
             $u(u^-(\al))=\al$ and $u(u^-(\beta))=\beta$.
Thus it
  follows from $u$ is Lipschitz
that
$$|\alpha - \beta|   \leq  K  | u^-(\al)  - u^-(\beta) |.$$
Similarly,
we can prove that
$|\gamma - \delta|   \leq  K | u^+(\gamma)  - u^+(\delta) |$
for all $\gamma$, $\delta$ in $[u( u^+(0)),1]$.

(\romannumeral2)$\Rightarrow$(\romannumeral3).
Given
$\alpha, \beta \in [u(u^-(0)),1]$,
if
$\beta\leq \alpha$,
then
$$|\alpha - \beta|   \leq K  | u^-(\al)  - u^-(\beta) | \leq  K  | u_s^-(\al)  - u^-(\beta) |.$$
If
$\beta > \alpha$,
note that
$u_s^-(\al)=\lim_{\lambda \to \al+}  u^-(\lambda)   $,
since
$$| \beta -\lambda |   \leq K  | u^-(\beta)  - u^-(\lambda) | ,$$
let $\lambda \to \alpha+$,
then
we
obtain
that
$$| \beta -\alpha |   \leq K  | u^-(\beta)  - u_s^-(\alpha) |.$$
Similarly,
we can prove that
$|\gamma - \delta|   \leq  K | u_s^+(\gamma)  - u^+(\delta) |$
for all $\gamma$, $\delta$ in $[u( u^+(0)),1]$.

(\romannumeral3)$\Rightarrow$(\romannumeral1).
Given $x,y \in [u]_0$, put $u(x)=\alpha$ and $u(y)=\beta$.
Assume
that
$\alpha<\beta$ with no loss of generality.
Then we have
$$ |x-y|  \geq  \min \{   |u_s^-(\al) - u^-(\beta)|, \     |u_s^+(\al) - u^+(\beta)|   \}.    $$
From statement (\romannumeral3),
\begin{gather*}
  |\alpha-\beta| \leq  K |u_s^-(\al) - u^-(\beta)|,
\\
   |\alpha-\beta| \leq  K |u_s^+(\al) - u^+(\beta)|,
\end{gather*}
and thus
$$|u(x)  -u(y)| = |\alpha-\beta| \leq  K |x-y|. $$
This means that $u$ is Lipschitz with Lipschitz constant $K$.
\ep

The following theorem shows that
the convolution transform can retain the
Lipschitz property under an assumption which is general for smoothers.

\begin{tm} \label{lsf}
  Suppose that $u$ is a fuzzy number, and that $v$ is a Lipschitz fuzzy number with Lipschitz constant $K$.
If
  $v(v^-(0))   \leq  u(u^-(0)) $
and
   $v(v^+(0))   \leq  u(u^+(0))$,
then
$u\nabla v $ is also a Lipschitz fuzzy number.
\end{tm}

\pof \ Note that, by Lemma \ref{evf},
$$
(u\nabla v) ((u\nabla v)^-(0))=u(u^-(0))\wedge v(v^-(0))=v(v^-(0)).
$$
Given
$\alpha$, $\beta$
in
$[(u\nabla v) ( (u\nabla v)^-(0)), 1]$,
then
$\alpha$, $\beta$
is also in
$[v( v^-(0)), 1]$.
Since $v$ is Lipschitz with Lipschitz constant $K$,
 it
 follows from    Lemma \ref{glec} (\romannumeral2)
that
$$|\alpha-\beta| \leq K | v^-(\alpha)   - v^-(\beta) | \leq K |(u\nabla v)^-(\al) -  (u\nabla v)^-(\beta)    |.$$
Similarly, we can obtain that
$$
|\gamma-\delta| \leq K | v^+(\gamma)   - v^+(\delta) | \leq K |(u\nabla v)^+(\gamma) -  (u\nabla v)^+(\delta) |
$$
for all $\gamma, \delta \in [(u\nabla v) ( (u\nabla v)^+(0)), 1]$.
So
$u\nabla v$ is also Lipschitz with Lipschitz constant $K$.
\ep

Since
we
 use
 polynomial functions,     sine functions and   cosine functions
to construct smoothers,
it is easy to
make the smoothers
 to be a Lipschitz fuzzy number.
Note that in the construction process of a smoother,
it requires that
$v(v^-(0)) = u(u^-(0))$ and $v(v^+(0)) = u(u^+(0))$,
where $u$ is the original fuzzy number and $v$ is its smoother (see condition (\romannumeral1)).
Thus, by Theorems \ref{san} and \ref{lsf}, we can
produce a
Lipschitz and smooth approximation for
a fuzzy number with finite non-differentiable points.
From the above discussions,
we
can
further
ensure
this
Lipschitz and smooth approximation
preserves the core at the same time.

\begin{re}{ \rm
  From Theorems \ref{dzeron}, \ref{dbzeron} and \ref{rdbzeron}, we know that
if $v$ is a smoother for $u \in \mathcal{F}(\mathbb{R})$,
then
\begin{gather*}
  v'(v^-(\al)) \geq  (u\nabla v)' (  (u\nabla v)^-(\al)   ) ,
 \\
  v'(v^+(\al)) \geq  (u\nabla v)' (  (u\nabla v)^+(\al)   ),
\end{gather*}
for
 all
$\al\in (0,1]$.
It follows immediately that
if
  a smoother $v$ of $u$ is Lipschitz,
then
 $u\nabla v$
is also Lipschitz.

}
\end{re}

\section{Conclusions}

This paper discusses how to smooth fuzzy numbers and then construct smooth approximations for
 fuzzy numbers by using
the convolution method.
The main contents
             are illustrated in the
                        following.
\begin{description}

\item[1] It shows that how to use the convolution method to produce smooth approximations for
 fuzzy numbers which have  finite non-differentiable points. This type of fuzzy numbers are quite general in real world applications.

\vspace{3mm}

  \item[2] It further points out that the convolution method can generate smooth and Lipschitz approximations which preserve the core at the same time.

\vspace{3mm}

\item[3] The constructing of smoothers is the key step in the construction processes of approximations
in
the above results.
          Theorems \ref{rap}, \ref{rndp} and \ref{rncp}
provide principles for constructing smoothers,
therein
conditions are given to ensure that the constructed fuzzy numbers are smoothers for a given type of fuzzy numbers.
These conditions are general.
In fact, by the conditions in Theorem \ref{rap}, we can judge that
the classes of fuzzy numbers $\{w_p\} $ and  $\{Z_p^f\}$
introduced in \cite{yo,yo2} are
smoothers
for
fuzzy numbers in $  \mathcal{F}^0_\mathrm{NC}( \mathbb{R} )$. See Remark \ref{usc} for details.

\end{description}

\appendix
  \renewcommand{\appendixname}{Appendix~\Alph{section}}

 \section{Proof of Theorem \ref{dzeron}   }

We
only
prove
statements (\romannumeral1) and (\romannumeral2).
The
 remainder
statements (\romannumeral3) -(\romannumeral8)
can be proved in the same way.

Set $v(v^-(0))=\al_0$ and $v(v^+(0))=\beta_0$, then, by Lemma \ref{evf},
\begin{gather*}
(u\nabla v)((u\nabla v)^-(0))\leq \alpha_0, \
(u\nabla v)((u\nabla v)^+(0))\leq \beta_0,
\\
(u\nabla v)((u\nabla v)^-(\alpha_0))= \alpha_0, \
(u\nabla v)((u\nabla v)^+(\beta_0))= \beta_0.
\end{gather*}

(\romannumeral1) \ If $\alpha=\alpha_0$,
then $ v^-(\al)=v^-(\al_0)=v^-(0)   $.
From
  $v'_-(   v^-(0)   )=0$,
we know
that
$v$ is left-continuous at $v^-(0)$.
Note that
$v(x)=0$ for all $x< v^-(0)$,
 thus $v( v^-(0)) =0=\alpha_0$.
Hence
 $(u\nabla v)((u\nabla v)^-(0))=0$,
 and therefore
 $(u\nabla v)'_-((u\nabla v)^-(0))=0$.

If $\al>\al_0$, then
by
 $v'_-(v^-(\al))=0$,
 we know that
 $v$ is left-continuous at $v^-(\al)$. So, by Proposition \ref{bse}, $v(v^-(\al))=\al$
 and hence, by Lemma \ref{evf},  $(u\nabla v)((u\nabla v)^-(\al))=\al$.

Now
we prove that
$(u\nabla v)'_-((u\nabla v)^-(\al))=0$.
Note that
$$
(u\nabla v)'_-((u\nabla v)^-(\al))
=
\lim_{z\to (u\nabla v)^-(\al)-}\frac{(u\nabla v)((u\nabla v)^-(\al))-(u\nabla v)(z)}{(u\nabla v)^-(\al)-z}.
$$
Given
$z\in ((u\nabla v)^-(\al_0),  (u\nabla v)^-(\al) )$, obviously
\begin{equation}
\frac{(u\nabla v)((u\nabla v)^-(\al))-(u\nabla v)(z)}     {(u\nabla v)^-(\al)-z}    \geq    0.
\label{deg}
\end{equation}
 On the other hand,
 set $(u\nabla v)(z)=\al-\delta$, where $\delta>0$, then by Proposition \ref{bse},
$$
z\leq (u\nabla v)_s^-(\alpha -   \delta),
$$
and thus
\begin{align}
       & \frac{(u\nabla v)((u\nabla v)^-(\al))-(u\nabla v)(z)}{(u\nabla v)^-(\al)-z}      \nonumber
\\
& \leq
\frac{\alpha-(\alpha-\delta)}
 { (u\nabla
v)^-(\al)-(u\nabla
v)_s^-(\alpha-  \delta)  }      \nonumber
\\
&= \frac{\delta }{
u^-(\al)-u_s^-(\alpha-\delta)
+
v^-(\al)-v_s^-(\al-\delta)
}                                    \nonumber
\\
&\leq \frac{\delta}{v^-(\al)-v_s^-(\alpha-\delta)}. \label{lde}
        \end{align}
By Proposition \ref{bse},
it holds that
$$\lim_{x \to v_s^-(\alpha-\delta)- }   v(x)  \leq  \alpha-\delta,$$
and therefore
\begin{equation}
\lim_{   x \to v_s^-(\alpha-\delta)-  }     \frac{  v(v^-(\al)) - v(x) }  { v^-(\al)-x }
\geq
\frac{\alpha-(\alpha-\delta)}{ v^-(\al)-  v_s^-(\alpha-\delta)}
=
\frac{\delta}   {   v^-(\al)-  v_s^-(\alpha-\delta)  }.
\label{leu}
\end{equation}
 Given $\varepsilon>0$,
 since $v'_-(v^-(\al))=0$,
 there is a $\xi (\varepsilon)>0$ such that
for all $y\in (v^-(\al)-\xi, v^-(\al))$, it holds that
\begin{equation}
\frac{v(v^-(\al))-v(y)}{v^-(\al)-y}  =    \frac{\alpha-v(y)}{v^-(\al)-y}  \leq   \varepsilon.
\label{dev}
\end{equation}
Hence,
 by \eqref{dev},
 if $v^-(\al)-v_s^-(\alpha-\delta)<\xi$,
 then
\begin{equation}
\lim_{   x \to v_s^-(\alpha-\delta)-  }
 \frac{  v(v^-(\al)) - v(x) }  { v^-(\al)-x }
 \leq
 \varepsilon,
 \label{les}
\end{equation}
and thus,
combined with \eqref{leu} and \eqref{les},
we have
\begin{equation}
\frac{\delta}   {   v^-(\al)-  v_s^-(\alpha-\delta) }     \leq       \varepsilon.
\label{esn}
\end{equation}
So if  $(u\nabla v)^-(\al)-z<\xi$, then
$v^-(\al)-v_s^-(\alpha-\delta)< (u\nabla v)^-(\al)-(u\nabla v)_s^-(\al-\delta)<(u\nabla v)^-(\al)-z <\xi$,
hence, by
\eqref{lde} and \eqref{esn},
we get that
$$\frac{(u\nabla v)((u\nabla v)^-(\al))-(u\nabla v)(z)}{(u\nabla v)^-(\al)-z}<\varepsilon,$$
from the arbitrariness of $\varepsilon>0$ and \eqref{deg},
we have
$$(u\nabla v)_-'((u\nabla v)^-(\al))=0.$$
So statement (\romannumeral1) holds.

(\romannumeral2) \
Since  $v(v^-(\al))=\beta$, we have $\beta \geq \alpha$,
and
$v^-(\al)=v^-(\beta)$. This implies that
$(u\nabla v)((u\nabla v)^-(\beta))=\beta$.

Now we show that $(u\nabla v)_+'((u\nabla v)^-(\beta))=0$
when
$v'_+(v^-(\al))=0$.
If
 $(u\nabla v)^-(\beta)   <          (u\nabla v)_s^-(\beta)$,
then it follows immediately
that
 $(u\nabla v)_+'((u\nabla v)^-(\beta))=0$.

 In the following, we
suppose that $(u\nabla v)^-(\beta)   =          (u\nabla v)_s^-(\beta)$,
then
$u^-(\beta)   =          u_s^-(\beta)$
and
$v^-(\beta)   =          v_s^-(\beta)$.
Note
 that
$$
(u\nabla v)'_+((u\nabla v)^-(\beta))
=
\lim_{z\to (u\nabla v)^-(\beta)+}
\frac
{(u\nabla v)(z)   -   (u\nabla v) ((u\nabla v)^-(\beta))}
{z   -   (u\nabla v)^-(\beta)}.
$$
The proof is divided into two cases.

Case (A) \
  $\beta<1$.

In this case, obviously,
\begin{equation}
\liminf_{z\to (u\nabla v)^-(\beta)+}
\frac
{(u\nabla v)(z)   -   (u\nabla v) ((u\nabla v)^-(\beta))}
{z   -   (u\nabla v)^-(\beta)}
\geq  0.
\label{lfu}
\end{equation}
Given
$z\in (  (u\nabla v)^-(\beta),   (u\nabla v)^-(1)   )$,
 set
 $(u\nabla v)(z)=\beta+\delta$, where $\delta>0$, then by Proposition \ref{bse},
$$
z\geq  (u\nabla v)^-(\beta  +   \delta),
$$
and thus
\begin{align}
& \frac
{(u\nabla v)(z)   -   (u\nabla v) ((u\nabla v)^-(\beta))}
{z   -   (u\nabla v)^-(\beta)}    \nonumber
\\
 \leq &
  \frac
{(u\nabla v)(z)   -   \beta}
{    (u\nabla v)^-(\beta+\delta)-   (u\nabla v)^-(\beta)      }    \nonumber
\\
=& \frac{\delta }{
u^-(\beta + \delta)-u^-(\beta)
+
v^-(\beta+\delta)-v^-(\beta)
}                                    \nonumber
\\
\leq&
 \frac{\delta}           { v^-(\beta+\delta)-v^-(\beta) }           \nonumber
 \\
 \leq&
 \frac{v(v^-(\beta+\delta))  -   v(v^-(\beta))}           { v^-(\beta+\delta)-v^-(\beta) }.
 \label{dre}
\end{align}
Let $z \to (u\nabla v)^-(\beta)+$, then $\delta \to 0+$. It thus follows from
\eqref{dre}
that
\begin{align*}
&\limsup_{z \to (u\nabla v)^-(\beta)+}
\frac
{(u\nabla v)(z)   -   (u\nabla v) ((u\nabla v)^-(\beta))}
{z   -   (u\nabla v)^-(\beta)}
\\
&\leq
\lim_{\delta \to 0+}
\frac{v(v^-(\beta+\delta))  -   v(v^-(\beta))}           { v^-(\beta+\delta)-v^-(\beta) }
\\
&=v'_+(v^-(\beta))=0,
\end{align*}
and then, combined with \eqref{lfu}, we know
$$ (u\nabla v)'_+((u\nabla v)^-(\beta))=0.   $$

Case (B) \
 $\beta=1$.

 In this case,
 if $(u\nabla v)^-(1)< (u\nabla v)^+(1) $, then
$ (u\nabla v)'_+   ((u\nabla v)^-(1))=0$.
If
 $(u\nabla v)^-(1)= (u\nabla v)^+(1) $, then $u^-(1)=u^+(1)$ and $v^-(1)=v^+(1)$.
So,
$$(u\nabla v)'_+((u\nabla v)^-(1))
=
(u\nabla v)'_+((u\nabla v)^+(1))
=
\lim_{z\to (u\nabla v)^+(1)+}\frac{(u\nabla v)((u\nabla v)^+(1))-(u\nabla v)(z)}{(u\nabla v)^+(1)-z}.$$
It is easy to see that
\begin{equation}
\limsup_{z\to (u\nabla v)^+(1)+}\frac{(u\nabla v)((u\nabla v)^+(1))-(u\nabla v)(z)}{(u\nabla v)^+(1)-z}
=
\limsup_{z\to (u\nabla v)^+(1)+}
\frac{1-(u\nabla v)(z)}     {(u\nabla v)^+(1)-z}
\leq
0.
\label{suprde}
\end{equation}
%
%
%
%
%
%
On the other hand, given $\varepsilon>0$, note that $v'_+(v^+(1))=v'_+(v^-(1))=0$, so
there is a $\xi>0$ such that
for each $y\in (v^+(1), v^+(1)+ \xi)$,
$$\frac{1-v(y) }{v^+(1)-y}      \geq      -\varepsilon,$$
proceed as in statement (\romannumeral1), we can get that
if $z-(u\nabla v)^+(1)< \xi$,
then
$$\frac{1-(u\nabla v)(z)}     {(u\nabla v)^+(1)-z}\geq -\varepsilon,$$
and hence from the arbitrariness of $\varepsilon$, we know
\begin{equation}
\liminf_{z\to (u\nabla v)^+(1)+}\frac{(u\nabla v)((u\nabla v)^+(1))-(u\nabla v)(z)}{(u\nabla v)^+(1)-z}
=
\liminf_{z\to (u\nabla v)^+(1)+}
\frac{1-(u\nabla v)(z)}     {(u\nabla v)^+(1)-z}
\geq
0.
\label{infrde}
\end{equation}
Combined with \eqref{suprde} and \eqref{infrde}, we obtain that
$$(u\nabla v)'_+((u\nabla v)^+(1))=0.$$
So statement (\romannumeral2) holds.
\ep

\section{Proof of Theorem \ref{dbzeron}    }

To prove Theorem \ref{dbzeron}, we need
 following lemmas and corollary.
Therein,
it
shows that the derivatives of a fuzzy number can be computed
by calculations which are determined by its   values at  the endpoints of $\alpha$-cuts and strong-$\alpha$-cuts.

\begin{lm} \label{lrde}
Let $u\in \mathcal{F}(\mathbb{R}) $. Set $u(u^-(0))=\alpha_0$ and $u(u^+(0))=\beta_0$.
Then the following statements hold.
\\
(\romannumeral1) \ Given $    \alpha  \in  (\alpha_0,   1 ] $, then
$u'_-(  u^-(\al)    )=  \varphi$
 if and only if
$$\lim_{\gamma \to  \alpha  -}
\frac   { u(u^-(\al))   -  \gamma   }
{  u^-(\al)  -   u^-(\gamma)   }
=
\lim_{\gamma \to  \alpha    -}
\frac   { u(u^-(\al))   -  \gamma   }
{ u^-(\al)    -    u_s^-(\gamma)}
=
\varphi. $$
\\
(\romannumeral2) \ Given $    \beta  \in  (\beta_0,   1 ] $, then
$u'_+(  u^+(\beta)     )=  \psi$
 if and only if
$$\lim_{\lambda  \to  \beta-}
\frac   { u(u^+(\beta))   -  \lambda   }
{  u^+(\beta)  -   u^+(\lambda)   }
=
\lim_{\lambda  \to  \beta-}
\frac   { u(u^+(\beta))   -  \lambda   }
{  u^+(\beta)  -   u_s^+(\lambda)   }
=
\psi.
$$

\end{lm}

\pof \ (\romannumeral1) \ Sufficiency. \
 Suppose that
$$\lim_{\gamma \to  \alpha  -}
\frac   { u(u^-(\al))   -  \gamma   }
{  u^-(\al)  -   u^-(\gamma)   }
=
\lim_{\gamma \to  \alpha    -}
\frac   { u(u^-(\al))   -  \gamma   }
{ u^-(\al)    -    u_s^-(\gamma)}
=
\varphi, $$
we show that
$$u'_-(  u^-(\al)    )  =  \lim_{  z\to u^-(\al)-    }
\frac{   u (  u  ^{-}(\alpha) )    -     u(z)    }
{u ^{-}(\alpha)  -  z} =   \varphi.$$
In fact,
 given
  $z
  \in
  (     u^-(0),    u^-(\alpha)   )$, suppose that $u (z)=\gamma$,
   then
 $ \gamma< \alpha $, and,
    by Proposition \ref{bse},
$$
 u ^{-}(\gamma)
\leq
z
  \leq
   u   _{s}^{-}(\gamma),$$
hence
\begin{equation*}
\frac{ u(u^-(\al))-\gamma}
{u ^{-}(\alpha)  -   u ^{-}(\gamma)}
\leq
\frac{   u (  u  ^{-}(\alpha) )    -     u(z)    }
{u ^{-}(\alpha)  -  z}
\leq
\frac{ u(u^-(\al)) -   \gamma}
{    u  ^{-}(\alpha)  -   u_s^{-}(\gamma)  },
\end{equation*}
note that if $\gamma \to \alpha-$, then $z\to u^-(\al)-$, and thus
\begin{equation*}
\lim_{  \gamma \to \alpha-   }
\frac{ u(u^-(\al))-\gamma}
{u ^{-}(\alpha)  -   u ^{-}(\gamma)}
=
\lim_{  z\to u^-(\al)-    }
\frac{   u (  u  ^{-}(\alpha) )    -     u(z)    }
{u ^{-}(\alpha)  -  z}
=
\lim_{   \gamma \to \alpha-   }
\frac{ u(u^-(\al)) -   \gamma}
{    u  ^{-}(\alpha)  -   u_s^{-}(\gamma)  }
=\varphi,
\end{equation*}
i.e.
 $$   u'_-(  u^-(\al)    )=  \varphi. $$

Necessity. \  Suppose that  $u'_-(  u^-(\al)    )=  \varphi$,
we prove that
$$\lim_{\gamma \to  \alpha  -}
\frac   { u(u^-(\al))   -  \gamma   }
{  u^-(\al)  -   u^-(\gamma)   }
=
\lim_{\gamma \to  \alpha    -}
\frac   { u(u^-(\al))   -  \gamma   }
{ u^-(\al)    -    u_s^-(\gamma)}
=
\varphi. $$

Now
we use a trick
which  also is used in the proof of Lemma \ref{dzeron}.
Given $\varepsilon>0$, since
$u'_-(u^-(\al))=\varphi$,
there is a $\xi>0$, such that
for all
$y\in (u^-(\al)-\xi,   u^-(\al))$,
$$
\varphi-\varepsilon   \leq    \frac{ u(u^-(\al)) -u(y)}{u^-(\al)-y}   \leq    \varphi+\varepsilon,
$$
notice that, for each $\beta$,
$$\lim_{x\to u^-(\beta)-}     u(x)    \leq    \beta  \leq   u(  u^-(   \beta  ) ),$$
and hence if $0< u^-(\al)- u^-(   \gamma   )<\xi$,
then
\begin{gather*}
\varphi+\varepsilon
\geq
 \lim_{   x\to    u^-(\gamma)-    }        \frac{ u(u^-(\al)) - u(x)  }    {u^-(\al)-x}
\geq
   \frac  {u(u^-(\al))- \gamma}    {  u^-(\al)  -   u^-(\gamma)    }
   \geq
\frac{ u(u^-(\al)) -  u(u^-(\gamma))  } {u^-(\al)  -   u^-(\gamma) }
\geq
\varphi-\varepsilon,
\end{gather*}
now let $\gamma  \to \alpha-$,
 we obtain that
\begin{gather*}
\lim_{\gamma  \to \alpha-}
 \frac  {u(u^-(\al))- \gamma}    {  u^-(\al)  -   u^-(\gamma)    }
=\varphi
 =u'_-(u^-(\al)).
\end{gather*}
Similarly, we can prove that
\begin{gather*}
\lim_{\gamma  \to \alpha-}
 \frac  {u(u^-(\al))- \gamma}    {  u^-(\al)  -   u_s^-(\gamma)    }
=\varphi
 =u'_-(u^-(\al)).
\end{gather*}

(\romannumeral2) \ Since $ u'_+(  u^+(\beta)     )=  \psi  $,
we know that
 if $\lambda \to \beta-$, then $u^+(\lambda)   \to    u^+(\beta)+$.
The remainder proof is
similarly to
the proof of statement (\romannumeral1)
\ep

\begin{tl} \label{lrdef}
Let $u\in \mathcal{F}(\mathbb{R}) $. Set $u(u^-(0))=\alpha_0$ and $u(u^+(0))=\beta_0$.
Then the following statements hold.
\\
(\romannumeral1) \ Given $    \alpha  \in  (\alpha_0,   1 ] $, then
$u'_-(  u^-(\al)    )=  \varphi$
 if and only if  $u(  u^-(\al)    )   =    \alpha$,
 and
$$\lim_{\gamma \to  \alpha  -}
\frac   {\alpha   -  \gamma   }
{  u^-(\al)  -   u^-(\gamma)   }
=
\lim_{\gamma \to  \alpha    -}
\frac   {\alpha  -  \gamma   }
{ u^-(\al)    -    u_s^-(\gamma)}
=
\varphi. $$
\\
(\romannumeral2) \ Given $    \beta  \in  (\beta_0,   1 ] $, then
$u'_+(  u^+(\beta)     )=  \psi$
 if and only if
 $u(  u^+(\beta)     )= \beta$,
 and
$$\lim_{\lambda  \to  \beta-}
\frac   {\beta   -  \lambda   }
{  u^+(\beta)  -   u^+(\lambda)   }
=
\lim_{\lambda  \to  \beta-}
\frac   {\beta   -  \lambda   }
{  u^+(\beta)  -   u_s^+(\lambda)   }
=
\psi.
$$

\end{tl}

\pof \ Note that if $u'_-(  u^-(\al)    )=  \varphi$, then $u$ is left-continuous at $u^-(\al)$,
and then
$$\lim_{x\to  u^-(\al)-}   u(x)=u(u^-(\al))=\al.$$
Similarly, if
$u'_+(  u^+(\beta)     )=  \psi$,
then
$$\lim_{x\to  u^+(\beta)+}   u(x)=u(u^+(\beta))=\beta. $$
So the desired results follow immediately from Lemma \ref{lrde}. \ep

\begin{lm}
\label{lrsd}
Let $u\in \mathcal{F}(\mathbb{R}) $.
Then the following statements hold.
\\
(\romannumeral1) \ Given $   u^-(\al)   \in    [  u^-(0), u^-(1)   ) $,
 set $ u( u^-(\al) )=\rho$, then
$u'_+(  u^-(\al)    )=  \phi$
is equivalent to
$$\lim_{\gamma \to \rho  +}
\frac   {\rho   -  \gamma   }
{  u^-(\al)  -   u^-(\gamma)   }
=
\lim_{\gamma \to  \rho   +}
\frac   { \rho  -  \gamma   }
{ u^-(\al)    -    u_s^-(\gamma)}
=
\phi. $$
\\
(\romannumeral2) \ Given  $   u^+(\beta)   \in    (   u^+(1),   u^+(0)  ]  $, set $ u( u^+(\beta) )=\sigma$, then
$u'_-(  u^+(\beta)    )=  \omega$
is equivalent to
$$\lim_{\gamma \to \sigma  +}
\frac   {\sigma   -  \gamma   }
{  u^+(\beta)  -   u^+(\gamma)   }
=
\lim_{\gamma \to    \sigma   +}
\frac   { \sigma  -  \gamma   }
{ u^+\beta)    -    u_s^+(\gamma)}
=
\omega. $$
\end{lm}

\pof \ (\romannumeral1) \ Since $ u( u^-(\al) )=\rho$, we know that $u^-(\al) = u^-(\rho)$ and $\rho<1$.
The proof is divided into two cases.

Case (A) \  $u^-(\al) <  u_s^-(\al) $.

In this case, $u'_+(  u^-(\al)    )=  0$, and
$$\lim_{\gamma \to \rho  +}
\frac   {\rho   -  \gamma   }
{  u^-(\al)  -   u^-(\gamma)   }
=
\lim_{\gamma \to  \rho   +}
\frac   { \rho  -  \gamma   }
{ u^-(\al)    -    u_s^-(\gamma)}
=
\frac   {  0 }
{  u^-(\al)  -   u_s^-(\alpha)   }
=    0,
$$
so the statement holds.

Case (B) \ $u^-(\al) =  u_s^-(\al) $.

 In this cases,
 $ u^-(\gamma) \to   u^-(\al) +$ as $\gamma \to \rho  +$. The rest of the proof
 is similar to   the proof of statement (\romannumeral1) in Lemma \ref{lrde}.

(\romannumeral2) \ The proof is similar to statement (\romannumeral1).   \ep

\textbf{Proof of Theorem \ref{dbzeron}.}
(\romannumeral1)
\ Set
$\alpha_0=u( u^-(0) )$ and $\alpha_1= v (  v^-(0)  )$.
Since
 $u_-'(u^-(\al))=\varphi>0$ and $v_-'(v^-(\al))=\psi>0$,
we know that $  \alpha >   \max   \{  \alpha_0, \alpha_1  \}    $,
and
 \begin{gather}
 \lim_{x\to  u^-(\al)-}   u(x)=u(u^-(\al))=\al,    \nonumber
 \\
  \lim_{y\to  v^-(\al)-}   v(y) =v(v^-(\al))=\al,    \nonumber
  \\
  (u\nabla v)((u\nabla v)^-(\alpha))=\alpha, \label{acp}
  \end{gather}
  thus $\alpha>  (u\nabla v)((u\nabla v)^-(0))=\min\{  \alpha_0, \alpha_1    \}$.
To prove $(u\nabla v)_-'(     (u\nabla v)  ^- (\al)    )=(\varphi^{-1} + \psi^{-1})^{-1}$,
by Corollary \ref{lrdef} and \eqref{acp},
we only need to show that
$$
\lim_{\gamma \to  \alpha  -}
\frac   {\alpha   -  \gamma   }
{ (u\nabla v)^-(\al)  -   (u\nabla v)^-(\gamma)   }
=
\lim_{\gamma \to  \alpha    -}
\frac   {\alpha  -  \gamma   }
{ (u\nabla v)^-(\al)    -    (u\nabla v)_s^-(\gamma)}
=
(\varphi^{-1} + \psi^{-1})^{-1}.
$$
In
fact,
\begin{align}
&\lim_{\gamma \to  \alpha  -}
\frac   {\alpha   -  \gamma   }
{ (u\nabla v)^-(\al)  -   (u\nabla v)^-(\gamma)   } \nonumber
\\
&= \lim_{ \gamma \to  \alpha  -     }
\frac
{  \alpha- \gamma   }
{u ^{-}(\alpha)   + v ^{-}(\alpha ) - (u ^{-}(\gamma)   + v ^{-}(\gamma))}\nonumber
\\
&= \lim_{\gamma \to  \alpha  -  } \frac{1}
{
\frac{u ^{-}(\alpha)- u ^{-}(\gamma)}{ \alpha -  \gamma   }
+
\frac{v ^{-}(\alpha)- v ^{-}(\gamma)}{ \alpha-  \gamma }
}.
\label{detf}
\end{align}
Since $u_-'(u^-(\al))=\varphi>0$ and $v_-'(v^-(\al))=\psi>0$,
by Corollary \ref{lrdef},
we have
\begin{gather*}
  \lim_{\gamma \to  \alpha  -  } \frac{u ^{-}(\alpha)- u ^{-}(\gamma)}{ \alpha -  \gamma   }   =\varphi^{-1},
  \\
   \lim_{\gamma \to  \alpha  -  } \frac{v ^{-}(\alpha)- v ^{-}(\gamma)}{ \alpha -  \gamma   }   =\psi^{-1},
  \end{gather*}
  and thus, combined with \eqref{detf}, we get
    $$
\lim_{\gamma \to  \alpha  -}
\frac   {\alpha   -  \gamma   }
{ (u\nabla v)^-(\al)  -   (u\nabla v)^-(\gamma)   }
=
(\varphi^{-1} + \psi^{-1})^{-1}.
$$
  Similarly, we can obtain
  $$\lim_{\gamma \to  \alpha    -}
\frac   {\alpha  -  \gamma   }
{ (u\nabla v)^-(\al)    -    (u\nabla v)_s^-(\gamma)}
=
(\varphi^{-1} + \psi^{-1})^{-1}.$$
So
  statement (\romannumeral1) is proved.

%
%
%
%
%
%

(\romannumeral2) \ Since $u(u^-(\al))=v(v^-(\al))=\beta$,
we know that
$
u^-(\alpha)=u^-(\beta)
$,
$
v^-(\alpha)=v^-(\beta),
$
and then
$$(u  \nabla  v)    (   (u  \nabla  v)^-(\beta)   )  =   \beta.$$
Note
 that
$u_+'(u^-(\al))=\varphi>0$ and $v_+'(v^-(\al))=\psi>0$,
so
$u^-(\al) <  u^-(1)$
and
$v^-(\al) <  v^-(1)$
, and thus
$$\alpha \leq \beta < 1.$$
To prove that
$$(u\nabla v)_+'(     (u\nabla v)  ^- (\beta)    )=(\varphi^{-1} + \psi^{-1})^{-1},$$
by Lemma \ref{lrsd},
we only need to show that
$$
\lim_{ \gamma \to  \beta+   }
\frac   { \beta -\gamma}
{  (u\nabla v)^-(\beta)  -   (u\nabla v)^-(\gamma)   }
=
\lim_{ \gamma \to  \beta+   }
\frac   { \beta -\gamma}
{  (u\nabla v)^-(\beta)  -   (u\nabla v)_s^-(\gamma)   }
=
(\varphi^{-1} + \psi^{-1})^{-1}. $$
In
fact,
reasoning as in the proof of statement (\romannumeral1),
we obtain
\begin{align*}
&
\lim_{ \gamma \to  \beta+   }
\frac   { \beta -\gamma}
{  (u\nabla v)^-(\beta)  -   (u\nabla v)^-(\gamma)   }       \nonumber
\\
&= \lim_{ \gamma \to  \beta+      }
\frac
{ \beta- \gamma   }
{u ^{-}(\beta)   + v ^{-}(\beta ) - (u ^{-}(\gamma)   + v ^{-}(\gamma))}\nonumber
\\
&= \lim_{\gamma \to  \beta+   }
\frac{1}
{
\frac{u ^{-}(\beta)- u ^{-}(\gamma)}{ \beta -  \gamma   }
+
\frac{v ^{-}(\beta)- v ^{-}(\gamma)}{ \beta-  \gamma }
}                                                        \nonumber
\\
&=(\varphi^{-1} + \psi^{-1})^{-1}.
\end{align*}
Similarly, we can get
$$
\lim_{ \gamma \to  \beta+   }
\frac   { \beta -\gamma}
{  (u\nabla v)^-(\beta)  -   (u\nabla v)_s^-(\gamma)   }
=
(\varphi^{-1} + \psi^{-1})^{-1}.
$$
So
 statement (\romannumeral2) is proved.

(\romannumeral3) \ Since
$u(u^-(\al))=\beta$, and $v(v^-(\al))=\gamma>\beta$,
we know
that
$\alpha\leq \beta < \gamma$,
and
\begin{gather*}
u^-(\alpha)   =   u^-(\beta)<u^-(1),
\\
v^-(\alpha)   =   v^-(\gamma),
\\
(u\nabla v)     (     (u\nabla v)  ^- (\beta)    )   =\beta      <  1.
\end{gather*}
Note that
$v^-(\lambda)  =   v^-(\al)$ for all $\lambda \in [\alpha, \gamma]$,
hence
$$(u\nabla v)^-(\rho) =  u^-(\rho)    +    v^-(\rho)  =    u^-(\rho)   +    v^-(\beta)$$
for each $\rho \in  [\beta, \gamma]$,
and thus
\begin{align}
 & \lim_{ \theta \to  \beta+   }
\frac   { \beta -\theta}
{  (u\nabla v)^-(\beta)  -   (u\nabla v)^-(\theta)   }           \nonumber
\\
&=\lim_{ \theta \to  \beta+   }
\frac   { \beta -\theta}
{  u^-(\beta)    +  v^-(\beta)    -     (  u^-(\theta)  +  v^-(\beta)   )      }   \nonumber
\\
&=\lim_{ \theta \to  \beta+   }
\frac   { \beta -\theta}
{  u^-(\beta)    -      u^-(\theta)       }           \nonumber
\\
&= u'_+(  u^-(\beta)    )   =    \varphi.
\label{acd}
\end{align}
Similarly, from
$v_s^-(\lambda)  =   v^-(\al)$ for all $\lambda \in [\alpha, \gamma)$,
we get
$$(u\nabla v)_s^-(\rho) =  u_s^-(\rho)    +    v_s^-(\rho)  =    u_s^-(\rho)   +    v^-(\beta)$$
for each $\rho \in  [\beta, \gamma)$,
and thus
\begin{align}
& \lim_{ \theta \to  \beta+   }
\frac   { \beta -\theta}
{  (u\nabla v)^-(\beta)  -   (u\nabla v)_s^-(\theta)   }          \nonumber
\\
&=\lim_{ \theta \to  \beta+   }
\frac   { \beta -\theta}
{  u^-(\beta)    +  v^-(\beta)    -     (  u_s^-(\theta)  +  v^-(\beta)   )      }     \nonumber
\\
&=\lim_{ \theta \to  \beta+   }
\frac   { \beta -\theta}
{  u^-(\beta)    -      u_s^-(\theta)       }          \nonumber
\\
&= u'_+(  u^-(\beta)    )  =  \varphi.
\label{sacd}
\end{align}
Now
it follows from Lemma \ref{lrsd}, \eqref{acd} and \eqref{sacd} that
$(u\nabla v)_+'(     (u\nabla v)  ^- (\beta)    )=\varphi$.

(\romannumeral4) \ We assert that $u^-(\alpha)> u^-(0)$.
On the contrary,
 if $u^-(\al) = u^-(0)$,
then
from $u_-'(u^-(\al))=\varphi$,
we know
$u( u^-(\al)  )  =\alpha =  u( u^-(0)  )  = 0$.
Note that
 $\lim_{y\to   v^-(\al)- }  v(y) = \lambda < \alpha$,
this yields that $\lambda<0$,
which is a
contradiction.

Set $u( u^-(0) ) = \alpha_0$, then $u^-(\alpha_0)  =   u^-(0)$,
and hence $u^-(\alpha) > u^-(\al_0)$.
This implies that
$\alpha> \al_0$.
Since $u_-'(u^-(\al))=\varphi$, we know that $u(u^-(\al))=\alpha$,
and thus
$$(u\nabla w) (   (u\nabla w)^-(\al)   )=\alpha.$$
 By
  $\lim_{y\to   v^-(\al)- }  v(y) = \lambda < \alpha$,
we get
that
\begin{gather*}
 v^-(\rho) =  v^-(\alpha)
\end{gather*}
for all $\rho \in (\lambda, \alpha]$,
and hence
$$(u\nabla v)^-(\rho) =  u^-(\rho)    +    v^-(\rho)  =    u^-(\rho)   +    v^-(\alpha)$$
for each $\rho \in  (\lambda, \alpha]$.
Thus
\begin{align}
 &  \lim_{ \rho \to  \al-   }
\frac   {\al -\rho}
{  (u\nabla v)^-(\alpha)  -   (u\nabla v)^-(\rho)   }            \nonumber
\\
&=\lim_{ \rho \to \al-   }
\frac   {\al -\rho}
{  u^-(\alpha)    +  v^-(\alpha)    -     (  u^-(\rho)  +  v^-(\alpha)   )      }   \nonumber
\\
&=\lim_{ \rho \to \al-   }
\frac   {\al -\rho}
{  u^-(\alpha)    -      u^-(\rho)       }           \nonumber
\\
&= u'_-(  u^-(\alpha)    )   =    \varphi.
\label{ecd}
 \end{align}
Similarly, note that
$v_s^-(\rho)  =   v^-(\al)$ for all $\rho \in [\lambda, \alpha)$,
we get
$$(u\nabla v)_s^-(\rho) =  u_s^-(\rho)    +    v_s^-(\rho)  =    u_s^-(\rho)   +    v^-(\alpha)$$
for each $\rho \in  [\lambda, \alpha)$,
and thus
\begin{align}
 \lim_{ \rho \to  \al-   } &
\frac   {\al -\rho}
{  (u\nabla v)^-(\alpha)  -   (u\nabla v)_s^-(\rho)   }  \nonumber
\\
=&\lim_{ \rho \to \al-   }
\frac   {\al -\rho}
{  u^-(\alpha)    +  v^-(\alpha)    -     (  u_s^-(\rho)  +  v^-(\alpha)   )      }  \nonumber
\\
=&\lim_{ \rho \to \al-   }
\frac   {\al -\rho}
{  u^-(\alpha)    -      u_s^-(\rho)       }     \nonumber
\\
=& u'_-(  u^-(\alpha)    )   =    \varphi.
\label{secd}
\end{align}
Now
it follows from Corollary \ref{lrdef}, \eqref{ecd} and \eqref{secd} that
$(u\nabla v)_-'(     (u\nabla v)  ^- (\alpha)    )=\varphi$.
\ep

\section*{Acknowledgements}
This work
was
supported by the National Natural
Science Foundation of China (Grant No. 61103052).
The authors would like to thank the referees for their comments and suggestions
which have been
very helpful in improving this paper.

\end{document}